\documentclass{article}
\usepackage[left=3cm,right=3cm,top=2.6cm,bottom=3.6cm]{geometry}

\usepackage{amsmath,amsthm,amssymb,mathtools}                  
\usepackage{mathptmx,siunitx,isotope}
\usepackage[utf8]{inputenc}
\usepackage{booktabs,url,caption}       
\usepackage[arrow, matrix, curve, 2cell]{xy}  \UseAllTwocells

\usepackage[mathscr]{euscript}
\usepackage{stmaryrd} 
\usepackage{hyperref} 
\usepackage[T1]{fontenc}
\usepackage{cleveref}
\usepackage[dvipsnames]{xcolor}
\usepackage[backend=biber, 
style=alphabetic,
giveninits=true,
isbn=false, 
url=false,
minalphanames=1,
maxalphanames=11,
maxbibnames=99]{biblatex}
\AtEveryBibitem{\clearfield{month}}
\DeclareNameAlias{author}{given-family}
\renewbibmacro{in:}{}
\DeclareFieldFormat{extraalpha}{#1}
\DeclareLabelalphaTemplate{
	\labelelement{
		\field[final]{shorthand}
		\field{label}
		\field[strwidth=3,strside=left,ifnames=1]{labelname}
		\field[strwidth=1,strside=left]{labelname}
	}
}
\usepackage{multicol}
\usepackage{setspace}
\usepackage{tikz}
\usetikzlibrary{decorations.markings}
\usetikzlibrary{cd}
\usepackage{soul}  
\setuldepth{Berlin}

\usepackage{comment}

\usepackage{todonotes}
\setuptodonotes{size=\footnotesize, color=orange!40}

\numberwithin{equation}{section}

\usetikzlibrary{math} 
\usetikzlibrary{decorations.markings}

\newcommand{\Daggercat}{{\operatorname{Cat}^{\dagger}}} 
\newcommand{\Oonevolcat}{{\operatorname{Cat}^{\Oone\text{-vol}}}} 
\newcommand{\opp}{\operatorname{op}} 
\newcommand{\oneop}{\operatorname{1-op}} 
\newcommand{\twoop}{\operatorname{2-op}} 
\newcommand{\onetwoop}{\operatorname{(1,2)-op}} 
\newcommand{\symmoncat}{\text{SymMon}\mathscr{C}\text{at}} 
\newcommand{\symmoncatd}{\text{SymMon}\mathscr{C}\text{at}^{\mathrm{d}}} 
\newcommand{\efix}{\exists\operatorname{fix}} 
\newcommand{\pos}{\operatorname{pos}} 
\newcommand{\Ide}{I\hspace{-0.02cm}} 
\newcommand{\Idewo}{I} 
\newcommand{\rSpintwovoltwocat}{{\text{2Cat}^{\rSpintwo\text{-vol}}}} 
\newcommand{\strictrSpintwovoltwocat}{{\text{2Cat}^{\rSpintwo\text{-}\dagger}}} 
\newcommand{\symmontwocat}{\text{SymMon}\BiCatadj} 
\newcommand{\Oonevolbicat}{{\text{2Cat}^{\Oone\text{-vol}}}} 
\newcommand{\strictOonevolbicat}{{\text{2Cat}^{\Oone\text{-}\dagger}}} 
\newcommand{\Otwovolbicat}{{\text{2Cat}^{\Otwo\text{-vol}}}} 
\newcommand{\strictOtwovolbicat}{{\text{2Cat}^{\Otwo\text{-}\dagger}}} 
\newcommand{\strictbivolbicat}{{\text{2Cat}^{\text{bi}\text{-}\dagger}}} 
\newcommand{\eu}{\operatorname{eu}} 
\newcommand{\DD}{\operatorname{DD}} 
\newcommand{\Vectbdl}{\operatorname{VectBdl}} 
\newcommand{\Linebdl}{\operatorname{LineBdl}} 

\newcommand{\Z}{\mathbb{Z}} 
\newcommand{\C}{\mathscr{C}} 
\newcommand{\D}{\mathscr{D}} 
\newcommand{\B}{\mathscr{B}} 
\newcommand{\id}{\operatorname{id}} 
\newcommand{\Id}{\operatorname{Id}} 
\newcommand{\orb}{\operatorname{orb}} 
\newcommand{\BiCat}{2\textrm{Cat}} 
\newcommand{\TriCat}{3\textrm{Cat}} 
\newcommand{\BiCatadj}{2\textrm{Cat}_{\text{adj}}} 
\newcommand{\Cat}{\mathscr{C}\text{at}} 
\newcommand{\Hom}{\mathscr{H}\text{om}} 
\newcommand{\Nat}{\mathscr{N}\text{at}} 
\newcommand{\Bord}{\text{Bord}} 
\newcommand{\ev}{\operatorname{ev}} 
\newcommand{\coev}{\operatorname{coev}} 
\newcommand{\Sothree}{\operatorname{SO}(3)} 
\newcommand{\Sotwo}{\operatorname{SO}(2)} 
\newcommand{\Soone}{\operatorname{SO}(1)} 
\newcommand{\rSpintwo}{\operatorname{Spin}(2)^r} 
\newcommand{\Spintwo}{\operatorname{Spin}(2)} 
\newcommand{\Onn}{\operatorname{O}(n)} 
\newcommand{\Otwo}{\operatorname{O}(2)} 
\newcommand{\Oone}{\operatorname{O}(1)} 
\newcommand{\Grb}{\operatorname{Grb}} 
\newcommand{\vect}{\operatorname{vect}_{\Bbbk}} 
\newcommand{\vectc}{\operatorname{vect}_{\mathbb{C}}} 
\newcommand{\LG}{\mathscr{L}\mathscr{G}_{\Bbbk}} 
\newcommand{\RW}{\mathscr{C}^{\operatorname{aRW}}} 
\newcommand{\hmf}{\operatorname{hmf}} 
\newcommand{\MF}{\mathcal{M}\hspace*{-0.07cm}\mathcal{F}} 
\newcommand{\Herm}{\operatorname{herm}} 
\newcommand{\BB}{\operatorname{B}\!} 

\newcommand{\xu}{\text{\ul{$x$}}}
\newcommand{\yu}{\text{\ul{$y$}}}
\newcommand{\zu}{\text{\ul{$z$}}}

\let\emph\undefined
\newcommand{\emph}[1]{\textsl{#1}}

\tikzcdset{
	arrow style=tikz,
	diagrams={>={Straight Barb[scale=2.7]}}
}

\tikzset{
	string/.style={draw=#1, postaction={decorate}, decoration={markings,mark=at position .51 with {\arrow[draw=#1]{>}}}},
	costring/.style={draw=#1, postaction={decorate}, decoration={markings,mark=at position .51 with {\arrow[draw=#1]{<}}}},
	ostring/.style={draw=#1, postaction={decorate}, decoration={markings,mark=at position .47 with {\arrow[draw=#1]{>}}}},
	ustring/.style={draw=#1, postaction={decorate}, decoration={markings,mark=at position .56 with {\arrow[draw=#1]{>}}}},
	oostring/.style={draw=#1, postaction={decorate}, decoration={markings,mark=at position .43 with {\arrow[draw=#1]{>}}}},
	uustring/.style={draw=#1, postaction={decorate}, decoration={markings,mark=at position .59 with {\arrow[draw=#1]{>}}}},
	directed/.style={string=blue!50!black}, 
	odirected/.style={ostring=blue!50!black}, 
	udirected/.style={ustring=blue!50!black}, 
	oodirected/.style={oostring=blue!50!black}, 
	uudirected/.style={uustring=blue!50!black},     
	redirected/.style={costring= blue!50!black},
	redirectedgreen/.style={costring= green!50!black},
	directedgreen/.style={string= green!50!black},
}

\tikzset{-dot-/.style={decoration={
			markings,
			mark=at position 0.5 with {\fill circle (1.875pt);}},postaction={decorate}}}

\tikzset{
	Fdot/.style={circle, draw, fill, inner sep=0pt}, 
	Odot/.style={circle, draw, inner sep=0.1pt, minimum size=0.1cm}
}

\theoremstyle{theorem} 
\newtheorem{theorem}{Theorem}[section] 
\theoremstyle{definition}
\newtheorem{definition}[theorem]{Definition}
\newtheorem{remark}[theorem]{Remark}

\newtheorem{lemma}[theorem]{Lemma}
\newtheorem{corollary}[theorem]{Corollary}
\newtheorem{proposition}[theorem]{Proposition}

\newtheorem{example}[theorem]{Example}

\newcounter{comments}
\newenvironment{commentfigure}{\begin{comment}}{\end{comment}}
\newenvironment{sidewayscommentfigure}{\begin{minipage}}{\end{minipage}}
\newenvironment{displaycomment}{\begin{list}{}{\rightmargin=1cm\leftmargin=1cm}\item\sf\begin{small}\color{gray}}{\end{small}\end{list}}
\def\comments{
	\setcounter{comments}{1}
	\renewenvironment{comment}{\begin{displaycomment}}{\end{displaycomment}}
	
	}

\makeatletter%
\newcommand{\xRrightarrow}[2][]{\ext@arrow 0359\Rrightarrowfill@{#1}{#2}}
\newcommand{\Rrightarrowfill@}{\arrowfill@\equiv\equiv\Rrightarrow}
\newcommand{\xLleftarrow}[2][]{\ext@arrow 3095\Lleftarrowfill@{#1}{#2}}
\newcommand{\Lleftarrowfill@}{\arrowfill@\Lleftarrow\equiv\equiv}
\makeatother
\newcommand{\RRightarrow}{\xRrightarrow{\phantom{--}}} 

\author{Nils Carqueville, Tim Lüders}
\title{\vspace{-1cm}Orbifolds, higher dagger structures, and idempotents}
\date{}

\comments 

\begin{document}
\maketitle
\begin{abstract}
    \begin{spacing}{1.125}
    The orbifold/condensation completion procedure of defect topological quantum field theories can be seen as carrying out a lattice or state sum model construction internal to an ambient theory. 
    In this paper, we propose a conceptual algebraic description of orbifolds/condensations for arbitrary tangential structures in terms of higher dagger structures and higher idempotents. 
    In particular, we obtain (oriented) orbifold completion from (framed) condensation completion by using a general 
    strictification procedure for higher dagger structures which we describe explicitly in low dimensions; we also discuss the spin and unoriented case. We provide several examples 
    of higher dagger categories, such as those associated to state sum models, (orbifolds of) Landau--Ginzburg models, and truncated affine Rozansky--Witten models. 
    We also explain how their higher dagger structures are naturally induced from rigid symmetric monoidal structures, recontextualizing and extending 
    results from the literature. 
    \end{spacing}
\end{abstract}

\setcounter{tocdepth}{2}
\tableofcontents
\section{Introduction}
\label{sec:Introduction}

The study of topological quantum field theories has a plethora of applications in pure mathematics and theoretical physics, and it has become a crucible for collaborations and cross-pollination. 
The simplest non-trivial class of examples are lattice or state sum models, which are generally expected to be fully extended TQFTs whose targets are vanilla categorifications of the category $\textrm{Vect}$ of vector spaces. 
In two dimensions this is a theorem for oriented theories: the non-extended state sum models of \cite{Fukuma1994, LAUDA2007} lift to extended TQFTs with values in the 2-category of algebras, bimodules and intertwiners, see \cite{26, 82}. 
In three dimensions, the results of \cite{DSPS} strongly suggest that Turaev--Viro--Barrett--Westbury models lift to fully extended TQFTs, with defects as described in \cite{meusburger2023statesummodelsdefects, carqueville2023orbifold}. 
It is natural to ask how this generalizes in various directions, including to tangential structures other than orientations. 

The standard state sum constructions can be formulated iteratively in the dimension, as explained for example in \cite{24} or \cite[Section~1]{carqueville2023orbifold}. 
This can naturally be thought of as internal to a specific defect TQFT which is trivial in the top dimension, and determined by lower-dimensional state sum models on proper defects. 
We refer to this defect TQFT as ``trivial''. 
For example in dimension two, line defects of the trivial defect TQFT may be arbitrary finite-dimensional vector spaces, and point defects are linear maps. 
In general dimension, the oriented state sum construction proceeds in four steps for every bordism~$M$: 
(1) choose a decomposition of~$M$ into submanifolds (such as the dual of a triangulation), 
(2) choose a defect label~$A_j$ for every $j$-dimensional submanifold that is compatible with the chosen decomposition, 
(3) evaluate the thus-obtained defect bordism with the trivial defect TQFT, and 
(4) take the colimit over all decomposition choices. 
The last step imposes a finite number of defining conditions on the data~$A_j$. 

In dimension one (where orientations and framings are equivalent), this amounts to specifying an ordinary idempotent (in $\textrm{Vect}$) to describe one-dimensional \textsl{closed} state sum models, and to the idempotent completion of $\textrm{Vect}$ (which is $\textrm{Vect}$ itself) to describe one-dimensional \textsl{defect} state sum models. 
Similarly, in dimension two, the procedure leads to the familiar algebraic input data of ($\Delta$-)separable symmetric Frobenius algebras for closed oriented state sum models, and to their pivotal Morita 2-category to describe defect state sum models. 

Depending on the relevant tangential structure and/or semantic conventions, the generalization of the state sum construction to arbitrary defect TQFTs (as introduced in \cite{6}, and not just the above ``trivial'' ones) goes by the name of ``orbifold completion'' \cite{Carqueville2012, 6} or ``condensation completion'' \cite{24, Johnson_Freyd_2022}. 
The former has been developed for \textsl{oriented} non-extended TQFTs in arbitrary dimension (see \cite{carqueville2023orbifoldstopologicalquantumfield} for a survey), while the latter is motivated by \textsl{framed} fully extended TQFTs.
In either case, the construction naturally amounts to a categorification of idempotent completion, namely of $n$-categories (with extra structure), where~$n$ is the dimension of the ambient defect TQFT. 
The objects of these completed $n$-categories are categorifications of idempotents themselves (for $n=2$, these are $\Delta$-separable Frobenius algebras in the framed case, and $\Delta$-separable symmetric Frobenius algebras in the oriented case), while morphisms are appropriate bimodules and their (higher) maps. 
In \Cref{subsec:IdempotentCompletion} we review more details for $n=1$ and $n=2$. 

It is natural to ask for a variant of the orbifold or condensation completion construction for arbitrary tangential structures. 
One option to address this question is topological and combinatorial in nature, namely to systematically study the interactions between the allowed cell decompositions (such as dual triangulations) and the chosen tangential structure. 
We do not pursue this direction here. 
Another option is purely algebraic and suggested by the cobordism hypothesis of \cite{lurie2009classificationtopologicalfieldtheories}: the categorical algebra pertaining to the framed case is in some sense fundamental, and for any other tangential structure one should take appropriate homotopy fixed points. 
This together with the work \cite{ferrer2024dagger} on ``$G$-volutive $n$-categories'' (the name draws from the fact that for $G=\Oone$ they are closely related to involutive structures, see below), and \cite{walker2021universalstatesum} suggesting a relation between state sum models of tangential structure $G$ and ``$G$-pivotal structures'', are three of our four main inspirations for the present paper. 

Our fourth inspiration is the fact that to an $n$-dimensional oriented defect TQFT one expects to be able to associate a ``pivotal $n$-category''. 
This is known to be true for $n=2$ \cite{55} and $n=3$ \cite{40}, where in the latter case we take a pivotal 3-category to be something that is equivalent to a Gray category with duals as introduced in \cite{Barrett2012}. 
On the other hand, pivotal 2-categories can be understood as a strict version of ``$\Sotwo$-volutive'' 2-categories, or $\Sotwo$-dagger 2-categories, in a way made precise in \Cref{Subsection: rspintwovolutive 2-categories} below, see also \cite{StehouwerMueller2024} and \cite[Example 5.10]{ferrer2024dagger}. 
More generally, a ``pivotal $n$-category'' is expected to be a strict version of an $\textrm{SO}(n)$-volutive $n$-category. 

The cobordism hypothesis with singularities \cite[Section 4.3]{lurie2009classificationtopologicalfieldtheories} suggests the existence of a canonical action $\Onn \to \operatorname{Aut}(n\textrm{Cat}_{\operatorname{adj}})$ on the category of $n$-categories
with adjoints, see \cite[Conjecture 5.3]{ferrer2024dagger}. Given a map $G \to \Onn$, a $G$-volutive $n$-category is defined to be a homotopy fixed point with respect 
to the induced $G$-action on $n\textrm{Cat}_{\operatorname{adj}}$. At least in low dimensions, it turns out that there is a natural notion of higher non-invertible morphisms between $G$-volutive $n$-categories, leading to an $(n+1)$-category $n\textrm{Cat}^{G\text{-vol}}$. 
We expect $n$-dimensional $G$-structured defect TQFTs (which have not yet been rigorously defined in general) to give rise to \textsl{strict}\footnote{%
What we call $G$-dagger is also called flagged $G$-dagger \cite{ferrer2024dagger} or $G$-pivotal \cite{walker2021universalstatesum}.
Moreover, $\Oone$-dagger categories are classically called dagger categories, $\Oone$-volutive categories are called anti-involutive \cite{Stehouwer2023DaggerCV}, 
$\Sotwo$-dagger 2-categories are called pivotal, and we expect $\Sothree$-dagger 3-categories to be what is called Gray categories with duals \cite{Barrett2012}. \label{fnt1}} 
$G$-volutive $n$-categories, or \textsl{$G$-dagger $n$-categories}, which in turn are the objects of the $(n+1)$-category $n\textrm{Cat}^{G\text{-}\dagger}$ whose rigorous (and more subtle) construction is outlined in \cite[Definition 2.6]{ferrer2024dagger}. Instead of trying to directly construct a completion procedure for $n\textrm{Cat}^{G\text{-}\dagger}$, 
it seems more natural to work with $n\textrm{Cat}^{G\text{-vol}}$ instead, provided we are able to compare the two $(n+1)$-categories $n\textrm{Cat}^{G\text{-}\dagger}$ and $n\textrm{Cat}^{G\text{-vol}}$. 

Thus we arrive at a program for how to systematically and algebraically approach orbifold/condensation completion associated to a given tangential structure of type~$G$. 
We shall first give a broad sketch of this program, and then summarize to what extent we implement it in the present paper.
First, we expect there to be ``trivial re-interpretation'' and ``strictification\footnote{As explained in \Cref{sec:DaggerStructures}, $G$-dagger higher categories have less coherent structure than $G$-volutive higher categories. We understand~$S_G$ to be a strictification in this sense, but we do not expect $T_G$ and $S_G$ to be fully faithful in general.}'' $(n+1)$-functors  
\begin{equation}
	\label{eq:SGTG}
	T_G \colon n\textrm{Cat}^{G\text{-}\dagger} 
		\xymatrix{\ar[r]&}
		n\textrm{Cat}^{G\text{-vol}}
		\, , \quad 
		S_G \colon n\textrm{Cat}^{G\text{-vol}}
		\xymatrix{\ar[r]&}
		n\textrm{Cat}^{G\text{-}\dagger}
\end{equation}
such that $S_G$ is right adjoint to $T_G$, $T_G \dashv S_G$. 
Next, we may ask for a ``$G$-idempotent completion'' $(n+1)$-functor 
\begin{equation}
	I_G \colon n\textrm{Cat}^{G\text{-vol}} 
	\xymatrix{\ar[r]&}
	n\textrm{Cat}^{G\text{-vol}}
\end{equation}
which endows the (framed!) idempotent completion of the underlying $n$-category of a $G$-volutive $n$-category with a compatible $G$-volutive structure. 
This should itself satisfy a universal property, in particular we should have $I_G \circ I_G \cong I_G$. 
Finally, we may go full circle and consider the composition $S_G\circ I_G \circ T_G$, viewing it as a construction on the $G$-dagger $n$-categories $\mathscr{D}_\mathscr{Z}$ associated to $n$-dimensional $G$-structured defect TQFTs~$\mathscr{Z}$. 
In the oriented case with $n\leqslant 3$, both the TQFTs $\mathscr{Z}$ and their $n$-categories $\mathscr{D}_{\mathscr{Z}}$ can be ``orbifold completed'' as well as ``Euler completed'' (see \Cref{Definition: Euler completion}), where the latter amounts to twisting the pivotal structure of~$\mathscr{D}_\mathscr{Z}$. These are compatible in the sense that
${\mathscr{D}_{{\mathscr{Z}}_{\orb}}} \cong ({\mathscr{D}_{\mathscr{Z}}})_{\orb}$ and $\mathscr{D}_{{\mathscr{Z}}_{\textrm{eu}}} \cong ({\mathscr{D}_{\mathscr{Z}}})_{\textrm{eu}}$ holds. 
It is thus natural to ask how these completions relate to $S_G,I_G$, and $T_G$. 
In summary, we expect: 
\begin{equation}
	\label{eq:TheDiagram}
		\begin{tikzcd}[column sep=3em, >=stealth]
			\begin{matrix} \text{$n$-dim.\ $G$-structured} \\ \text{defect TQFTs} \end{matrix}
				\;\;\;\;\;\; 
				\arrow[rr, dashed, "\text{$\mathscr{Z} \longmapsto \mathscr{D}_{\mathscr{Z}}$}"] 
				\ar[loop, ->, out=120, in=60, looseness=3, "\text{orbifold/condensation completion}"]
				\ar[loop, ->, in=-120, out=-60,  looseness=3, "\text{Euler completion}"]
			&& 
			\;\;\;\;\;\;  
			{n\,\textrm{Cat}^{G\text{-}\dagger}} 
				\arrow[rr, out=30, in=150, "T_{G}"] 
			& \perp  & 
			{n\,\textrm{Cat}^{G\text{-vol}}}\vphantom{\int}
				\arrow[ll, out=-150, in=-30, "S_{G}"] 
				\ar[loop, ->, in=-30, out=30, looseness=5, "I_G \;\cong\; I_G\circ I_G"]
		\end{tikzcd}
\end{equation}

In the present paper we study this proposal mostly for $n=2$, and for~$G$ being one of the groups $\Oone$, $\Sotwo$, $\Otwo$, $\textrm{Spin}(2)$, and more generally the $r$-spin group $\rSpintwo$ (which is the $r$-fold cover of $\Sotwo$, thus recovering $\Sotwo$ and $\textrm{Spin}(2)$ for $r=1$ and $r=2$, respectively). 

Before outlining our results, we summarize what is known about the case $n=1$ and $G=\Oone$. 
An \textsl{$\Oone$-volutive} structure on a category~$\C$ is a functor $d\colon \C \to \C^{\textrm{op}}$ together with a natural isomorphism $\eta \colon d^{\opp} \circ d \Rightarrow \id_{\C}$ such that $d(\eta_a) = \eta_{d(a)}^{-1}$ for all $a \in \C$. 
The strict version is known as a \textsl{dagger category}, where~$d$ is the identity on objects, and~$\eta$ is equal to the identity. 
The 2-functors $S_{\Oone}, T_{\Oone}$ of~\eqref{eq:SGTG} were in this case constructed and shown to participate in an adjunction $T_{\Oone} \dashv S_{\Oone}$ in \cite{Stehouwer2023DaggerCV}. 
A natural completion of dagger categories $(\C, d)$ was introduced as the \textsl{$d$-Karoubi envelope} $(\overline{\C}^d, d)$ in \cite{SELINGER2008107}, and may be identified with the completion situated on the left-hand side of~\eqref{eq:TheDiagram}. 
In \Cref{Subsection: Oone orbifold completion} we construct a 2-functor $I_{\Oone}\colon \textrm{Cat}^{\Oone\text{-vol}} \to \textrm{Cat}^{\Oone\text{-vol}}$ which basically lifts an $\Oone$-volutive structure on a 1-category~$\C$ to its (ordinary) idempotent completion $I\C$. 
The comparison is then made precise by showing (in \Cref{prop:FFdagger}) that there is a fully faithful dagger functor 
\begin{equation}
	\big(\overline{\C}^d,d\big) \xymatrix{\ar[r]&} \big(S_{\Oone} \circ I_{\Oone} \circ T_{\Oone})(\C,d\big)\,. 
\end{equation}
These results suggests that, while the $d$-Karoubi envelope is (in the sense of a universal property) the ``minimal'' dagger idempotent completion, our construction 
usually leads to a larger object; this is motivated by our interest in Euler and orbifold completion, cf.~\eqref{eq:SO2completions}. 
In other settings it is of interest to impose additional conditions or structure, such as unitarity in the context of higher Hilbert spaces, see e.g.\ \cite{chen2024manifestlyunitaryhigherhilbert}. 
We also remark that $S_{\Oone} \circ I_{\Oone} \circ T_{\Oone}$ naturally has the structure of a monad in $\Hom_{\BiCat}(\Daggercat,\Daggercat)$, but its multiplication need not be an equivalence. 

For $n=2$ and~$G$ one of the aforementioned groups we aim to obtain similarly clear pictures based on the basic idea summarized in~\eqref{eq:TheDiagram}, which we achieve to varying levels of completeness. 
For $G = \rSpintwo$ for $r\in\mathbb Z_+$, we propose complete definitions of the 3-categories $2\textrm{Cat}^{\textrm{Spin}(2)^r\text{-}\dagger}$ and $2\textrm{Cat}^{\textrm{Spin}(2)^r\text{-vol}}$ as well as the 3-functors $S_{\textrm{Spin}(2)^r}, T_{\textrm{Spin}(2)^r}$ between them, which we prove (in \Cref{prop:STadjointSpin}) to be adjoint: 
\begin{equation}
	\begin{tikzcd}[column sep=3em, >=stealth]
		{2\textrm{Cat}^{\rSpintwo\text{-}\dagger}}    \arrow[rr, out=30, in=150, "T_{\textrm{Spin}(2)^r}"] 
		& \perp  & 
		{2\textrm{Cat}^{\rSpintwo\text{-vol}}} .   \arrow[ll, out=-150, in=-30, "S_{\textrm{Spin}(2)^r}"] 
	\end{tikzcd}
\end{equation}
We do not construct the 3-functor $I_{\textrm{Spin}(2)^r}$ in general, but we do define it object-wise, i.e.\ to every 2-category~$\B$ with $\textrm{Spin}(2)^r$-volutive structure~$S$ (reviewed in \Cref{Subsection: rspintwovolutive 2-categories}) we construct an associated $\textrm{Spin}(2)^r$-volutive structure~$S'$ on the (framed) idempotent completion $I\B$, and we set $I_{\textrm{Spin}(2)^r}(\B, S) = (I\B, S')$. 
For $r=1$ we then show (in \Cref{Lemma: Euler completion from our perspective} and \Cref{Theorem: comparison of Euler orbifold completion and SIT}) that the approach suggested in~\eqref{eq:TheDiagram} recovers both the Euler completion and the Euler completed orbifold completion in the sense that there are equivalences
\begin{equation}
	\label{eq:SO2completions}
	\big(S_{\Sotwo} \circ T_{\Sotwo}\big)(\B, S) 
		\; \cong \; 
		(\B, S)_{\textrm{eu}}
		\, , \quad 
		\big(S_{\Sotwo} \circ I_{\Sotwo} \circ T_{\Sotwo}\big)(\B, S) 
		\; \cong \; 
		\big( (\B, S)_{\textrm{orb}} \big){}_{\textrm{eu}}
\end{equation}
for every $\Sotwo$-dagger 2-category $(\B, S)$, i.e.\ for every pivotal 2-category. 
Hence $S_{\Sotwo} \circ T_{\Sotwo}$ is an idempotent. 
We expect, but do not prove, that $I_{\textrm{Spin}(2)^r}$ is also an idempotent for any~$r$. Moreover, we expect that analogues of~\eqref{eq:SO2completions}
hold in arbitrary dimension and for arbitrary tangential structures. 

For $n=2$ and $G = \Oone$ or $G = \Otwo$ we similarly propose the main ingredients of the associated 3-categories, and we describe the 3-functors $S_G, T_G$ on (at least) object level. 
We also prove (in \Cref{Lemma: Oone volution on two idempotent completion} and \Cref{Proposition: Otwo volutive structures extend to the idempotent completion}) that $\Oone$- and $\Otwo$-volutive structures on 2-categories lift to their (framed) idempotent completions, thus defining the 3-functors $I_{\Oone}$ and $I_{\Otwo}$ on objects (which again we expect to be weak idempotents). 
While~\eqref{eq:SO2completions} arguably shows that $S_{\Sotwo} \circ I_{\Sotwo} \circ T_{\Sotwo}$ is the construction most relevant from a TQFT perspective, we do not expect $S_G\circ I_G\circ T_G$ itself to be an idempotent in general, see e.g.\ \Cref{rem:SITnotStrongIdempotent}. 

Along the way of establishing the above results, we show how to produce $G$-volutive 2-categories from symmetric monoidal 2-categories with duals and adjoints. 
It was more generally conjectured in \cite[Statement 6.1]{ferrer2024dagger} that any symmetric monoidal $n$-category with duals on all levels admits a canonical $\Onn$-volutive structure. 
The canonical $\Oone$-volutive structure based on the dualization functor was described before (see \Cref{Example: symmetric monoidal categories with duals are Oone volutive});
the canonical $\Sotwo$-volutive structure given by the Serre automorphism is described in \Cref{Example: Symmetric monoidal 2-categories with duals and adjoints are rSpintwo volutive} 
and \cite{StehouwerMueller2024}. We further corroborate the general picture~\eqref{eq:TheDiagram} by applying it to several examples (obtained in this way), where we find agreement with results obtained in the literature by other means.

The first examples we consider in \Cref{subsec:StateSumModels} are two-dimensional state sum models, where the underlying 2-category~$\B$ is the delooping of the category of finite-dimensional vector spaces. 
For~$G$ one of the groups $\Sotwo, \rSpintwo, \Otwo$ we construct the $G$-dagger structures and compute the action of the operation $S_G\circ I_G \circ T_G$. 
In all cases we find precise agreement with the literature, and a slight generalization in the case of $\Otwo$.  

The other examples we consider in \Cref{sec:ExamplesApplications} are the 2-categories of bundles gerbes, of Landau--Ginzburg models, and of truncated affine Rozansky--Witten models. 
In each case we construct $G$-dagger structures directly from the general theory. 
In the case of $G = \Sotwo$ for Landau--Ginzburg and Rozansky--Witten models, this recovers the previously established pivotal structures, while in the case of $G = \Otwo$ for Landau--Ginzburg models we identify a 2-categorical interpretation of the results of \cite{Hori2008} on orientifold branes. 

\medskip  

\noindent
\textbf{Acknowledgements.} 
We especially thank Lukas Müller and Luuk Stehouwer for many helpful discussion on higher dagger structures, as well as valuable feedback on an earlier draft of this manuscript. 
We are also grateful to Lukas for collaboration on the relationship between Serre automorphisms and pivotal structures, summarized in \Cref{Example: Symmetric monoidal 2-categories with duals and adjoints are rSpintwo volutive}, and to Aleksandar Ivanov for insightful discussions and comments.
Tim Lüders is supported by a scholarship of the Studienstiftung des deutschen Volkes e.V., while Nils Carqueville acknowledges partial support from the German Science Fund DFG (Heisenberg Programme) and the Austrian Science Fund FWF (project no.\ P\,37046).

\newpage 

\noindent
\textbf{Conventions.} 
We refer to bicategories, pseudofunctors, etc.\ as 2-categories, 2-functors etc., and we usually refer to equivalences in higher categories as isomorphisms. 
Occasionally we call (monoidal, higher) categories with fixed choices of duals and adjoints on all levels rigid. 
When filling pasting diagrams, we allow some flexibility in the labeling, see \Cref{Remark: Labeling of pasting diagrams}. 
For conventions on higher dagger structures see \Cref{fnt1}.

\section{Dagger structures and their generalizations}
\label{sec:DaggerStructures}

In this section we review and introduce higher categories of 1- and 2-categories that come from coherent actions of the groups $\Oone$, $\Sotwo$, $\Otwo$, as well as spin groups. 
In particular, we present (or: give a detailed sketch of) the 3-categories of ``$G$-volutive 2-categories'' and their stricter/dagger variants, where~$G$ is one of the aforementioned groups. 
We also describe (in part conjectural) strictification 3-functors as right adjoints to forgetful 3-functors. 
Most examples are deferred to \Cref{sec:ExamplesApplications}.

\subsection[Dagger categories]{Dagger categories}
\label{Subsection: Oonevolutive categories}

In this section we review the 2-categories of dagger categories and their non-strict variants called $\Oone$-volutive categories. 
Following \cite{Stehouwer2023DaggerCV}, we also exhibit an adjunction between them, given by trivially interpreting dagger categories as $\Oone$-volutive categories, and ``strictifying'' the latter to the former. 

\medskip 

Let $\C$ be a category. We denote by $\C^{\opp}$ the \emph{opposite} category with the same objects and morphisms given by 
$\hom_{\C^{\opp}}(a,b) = \hom_{\C}(b,a)$ for $a,b \in \C$. 
We identify $(\C^{\opp})^{\opp} = \C$. If $\B$ is a 2-category, 
we denote by $\B^{\oneop},\B^{\twoop},$ and $\B^{\onetwoop}$ the 2-categories with the same objects and only 1-morphism, only 2-morphisms, and both 1- and 2-morphisms reversed, 
respectively. 
Taking the opposite category then extends to a 2-functor $\Cat \to \Cat^{\twoop}$ where $\Cat$ denotes the (strict) 2-category of categories, functors, and natural 
transformations. 
For background on 2-categories we refer to \cite{johnson20202dimensional,Leinster1998BasicB}, and for 3-categories to \cite{39}.

\begin{definition}
	\label{def:DaggerCategory}
    Let $\C$ be a category. 
    A \emph{dagger structure} on $\C$ is a functor $d \colon \C \to \C^{\opp}$ such that $d^{\opp} \circ d = \id_{\C}$
    and $d(c) = c$ for all objects $c \in \C$. 
    The pair $(\C,d)$ is called a \emph{dagger category}.
\end{definition}

The second condition of the functor $d$ is often referred to as being \emph{identity-on-objects}. 
A morphism $u \colon a \to b$ in a dagger category $(\C,d)$ is an \emph{isometry} if $d(u) \circ u = \id_a$. 
If moreover $u \circ d(u) = \id_b$ holds, we call~$u$ \emph{unitary}. 
We also recall that, by definition, an automorphism $f \colon a \to a$ in a dagger category $(\C,d)$ is \emph{self-adjoint} if $f = d(f)$. 
A self-adjoint automorphism $f$ is \emph{positive} if there is another object $b \in \C$ and an isomorphism $g \colon a \to b$ such that $f = d(g) \circ g$. 
This terminology comes from the following standard example.

\begin{example} \label{Example: The category of hermitian vector spaces}
    Let $\Herm_\mathbb{C}$ be the category of finite-dimensional hermitian vector spaces and linear maps. In particular, any object in $\Herm_\mathbb{C}$
    comes with a non-degenerate hermitian sesquilinear form $\langle - , -\rangle$. Taking the adjoint with respect to the pairings $\langle - , -\rangle$ equips $\Herm_\mathbb{C}$ with a dagger structure.
\end{example}

\begin{definition}
    Let $(\C,d)$ and $(\C',d')$ be dagger categories. A \emph{dagger functor} $(\C,d) \to (\C',d')$ is a functor $F \colon \C \to \C'$ that satisfies 
    $F^{\opp} \circ d = d' \circ F$. 
\end{definition}

\begin{definition}
    Let $F_1, F_2 \colon (\C,d) \to (\C',d')$ be dagger functors. A \emph{dagger natural transformation} $F_1 \to F_2$ is a natural transformation 
    $F_1 \Rightarrow F_2$ whose components are isometries.
\end{definition}

We denote the (strict) 2-category of dagger categories, dagger functors, and dagger natural transformations by $\Daggercat$. 

\begin{remark}
    As remarked in \cite{Stehouwer2023DaggerCV}, the above notion of a dagger natural transformation recovers the appropriate notion
    of equivalence between dagger categories \cite{Vicary2011}. 
    Namely, a dagger functor $F$ is part of an equivalence of dagger categories if 
    and only if it is fully faithful and surjective on objects up to unitaries.
\end{remark}

\medskip 

From a categorical point of view, dagger categories are ill-behaved in the sense that, given a dagger category $(\C,d)$ and an equivalence of categories $\C \cong \C'$, the equivalence does in general not induce a natural dagger structure on $\C'$. 
There is a more coherent version of dagger categories that adresses this problem, which we will now recall; see \cite{Stehouwer2023DaggerCV} for more discussion.

\begin{definition}
	\label{def:O1volutiveCategory}
    Let $\C$ be a category. 
    An \emph{$\Oone$-volution} on $\C$ consists of a functor $d \colon \C \to \C^{\opp}$ and a natural isomorphism $\eta \colon d^{\opp} \circ d \Rightarrow \id_{\C}$ such that $d(\eta_a) = \eta_{d(a)}^{-1}$ for all $a \in \C$. 
    A category together with an $\Oone$-volution is called an \emph{$\Oone$-volutive category}.
\end{definition}

\begin{remark}
    Another name for $\Oone$-volutive is \emph{anti-involutive}. 
    We prefer ``$\Oone$-volutive'' to emphasize parallels between the coherent version of dagger structures discussed here and their (higher) variants in the subsequent subsections -- as well as their homotopy-theoretic origins, see \Cref{Remark: Oone volutive categories as homotopy fixed points}.
\end{remark}

\begin{definition}
    Let $(\C,d,\eta)$ and $(\C',d',\eta')$ be $\Oone$-volutive categories. An \emph{$\Oone$-volutive functor} $(\C,d,\eta) \to (\C',d',\eta')$ 
    consists of a functor $F \colon \C \to \C'$ and a natural isomorphism $\alpha \colon F^{\opp} \circ d \Rightarrow d' \circ F$ such that the following diagram 
    commutes: 
    \begin{equation}
        \begin{aligned}
        \xymatrix{
            (d')^{\opp} \circ F^{\opp} \circ d \ar@2{->}[rr]^-{\alpha^{\opp} \circ \id} \ar@2{->}[d]_-{\id \circ \alpha} && F \circ d^{\opp} \circ d \ar@2{->}[d]^-{\id \circ \eta} \\
            (d')^{\opp} \circ d' \circ F \ar@2{->}[rr]_-{\eta' \circ \id} && F
        }
        \end{aligned}
    \end{equation}
\end{definition}
\begin{definition}
    Let $(F_1,\alpha_1),(F_2,\alpha_2) \colon (\C,d,\eta) \to (\C',d',\eta')$ be $\Oone$-volutive functors. An \emph{$\Oone$-volutive natural transformation} 
    is a natural transformation $\varphi \colon F_1 \Rightarrow F_2$ such that the following diagram commutes: 
    \begin{equation}
        \begin{aligned}
        \xymatrix{
            F_1^{\opp} \circ d \ar@2{->}[rr]^-{\alpha_1} && d' \circ F_1 \ar@2{->}[d]^-{\id \circ \varphi} \\
            F_2^{\opp} \circ d \ar@2{->}[rr]_-{\alpha_2} \ar@2{->}[u]^-{\varphi^{\opp} \circ \id} && d' \circ F_2 
        }
    \end{aligned}
    \end{equation}
\end{definition}
$\Oone$-volutive categories, $\Oone$-volutive functors, and $\Oone$-volutive natural transformations form a 2-category \cite{Stehouwer2023DaggerCV},
which we denote by $\Oonevolcat$.

\begin{remark} 
	\label{Remark: Oone volutive categories as homotopy fixed points}
    Upon forgetting non-invertible 2-morphisms, that is, by passing to the underlying (2,1)-category, 
    $\Oonevolcat$ has a natural interpretation in terms of homotopy fixed points.
    To review this, let $G$ be a topological group and 
    $\B$ a 2-category. 
    A \emph{$G$-action} on $\B$ is a 3-functor 
    \begin{equation}
        \rho \colon \textrm{B}\Pi_2(G) \xymatrix{\ar[r]&} \BiCat \qquad \text{such that} \qquad 
        \rho(\star) = \B 
    \end{equation}
    where $\BiCat$ is the 3-category of 2-categories, 2-functors, 2-transformations, and modifications, and $\textrm{B}\Pi_2(G)$ is the delooping of the fundamental 
    2-groupoid of $G$ whose unique object is denoted $\star$. 
    The 2-category of \emph{homotopy fixed points} with respect to the $G$-action $\rho$ is then defined to be 
    the 3-limit of $\rho$, which is expected to be $\Nat(\Delta,\rho)$ where $\Delta$ is the constant 3-functor $\textrm{B}\Pi_2(G) \to \BiCat$ whose image is the trivial 
    2-category, and $\Nat$ is the 2-category of 3-transformations between two 3-functors. See \cite[Section 2.2]{82} for more details.
    
    In our present setting, we have $G = \Oone$ acting on the (2,1)-category\footnote{Note that it is necessary to restrict to the (2,1)-category 
    because taking the opposite is contravariant on 2-morphism level.} of $\Cat$ via $\C \mapsto \C^{\opp}$. 
	As (2,1)-categories, the (2,1)-category of homotopy fixed points with 
    respect to this action is precisely $\Oonevolcat$. 
    Similarly, as the group $\Soone$ is trivial, one could define \emph{$\Soone$-volutive} categories to be ordinary categories.
\end{remark}
\begin{example}\label{Example: symmetric monoidal categories with duals are Oone volutive}
    Any symmetric monoidal category $\C$ with duals admits the structure of an $\Oone$-volutive category. Indeed, recall that any choice of dualization data $(a^*,\widetilde{\ev}_a, \widetilde{\coev}_a)$ for objects $a\in \C$
    gives rise to a dualization functor $(-)^* \colon \C \to \C^{\opp}$ which assigns~$a$ to its chosen dual $a^*$, and morphisms $X \colon a \to b$ to 
    \begin{equation}
        \xymatrix{
            X^* \;=\;  
            \Big(b^* \ar[rr]^-{\widetilde{\coev}_a \otimes \id} && a^* \otimes a \otimes b^* \ar[rr]^{\id \otimes X \otimes \id} && a^* \otimes b \otimes b^* \ar[rr]^-{\id \otimes \widetilde{\ev}_b} && a^* \Big).
        }
    \end{equation}
    Recalling that any two choices of dualization data induce a natural isomorphism between the associated dualization functors, the symmetric braiding $\beta$ on $\C$ induces a natural isomorphism $\chi \colon ((-)^*)^{\opp} \circ (-)^* \Rightarrow \id_{\C}$. 
    One then checks that $((-)^*,\chi)$ is an $\Oone$-volutive structure on $\C$. 
\end{example}

As shown in \cite[Theorem 2.2.3]{Stehouwer2024Unitaryfermionic}, the above construction is functorial: 

\begin{proposition}
	\label{Proposition: From symmetric monoidal categories to Oone volutive categories}
    The assignment described in \Cref{Example: symmetric monoidal categories with duals are Oone volutive} extends to a 2-functor 
    \begin{equation}
        \symmoncat^{\textrm{d}} \xymatrix{\ar[r]&} \Oonevolcat
    \end{equation}
    from the 2-category of symmetric monoidal categories with duals to the 2-category of $\Oone$-volutive categories
\end{proposition}

This provides a large class of examples of $\Oone$-volutive categories, $\Oone$-volutive functors, and $\Oone$-volutive natural transformations between them. The following is a subexample. 

\begin{example} 
	\label{Example: vect is Oone volutive}
    Let $\vect$ be the symmetric monoidal category of finite-dimensional vector spaces over some field $\Bbbk$. 
    By \Cref{Example: symmetric monoidal categories with duals are Oone volutive}, it admits an $\Oone$-volutive structure whose underlying functor 
    assigns each vector space $V$ to its dual $V^* = \hom_{\Bbbk}(V,\Bbbk)$ and each linear map $f \colon V \to W$ to the linear map 
    $f^* \colon W^* \to V^*$, $\phi \mapsto \phi \circ f$. The component of the natural isomorphism $\eta$ at $V$ is given by 
    the canonical isomorphism $V^{**} \cong V$ between the double dual of~$V$ and itself.
\end{example}

In general, a category may admit several (inequivalent) $\Oone$-volutive structures. 
One way of constructing new $\Oone$-volutive structures from old ones is by composition with compatible (covariant) involutions, as the following demonstrates, cf.\ \cite[Example 3.2]{Stehouwer2023DaggerCV}. 

\begin{example}
	\label{Example: complex vector spaces with complex conjugation and dualization}
    Consider $\vectc$ and the complex conjugation functor $\overline{(-)} \colon \vectc \to \vectc$. The functor $(-)^* \circ \overline{(-)} \colon \vectc \to \vectc^{\opp}$ 
    together with the natural isomorphism $((-)^* \circ \overline{(-)})^{\opp} \circ ((-)^* \circ \overline{(-)}) \Rightarrow \id$ whose component at $V$ is given by 
    \begin{equation}
    	\xymatrix{
    		\overline{\overline{V}^*}^* 
    		\ar[r]^\cong & 
    		\overline{\overline{V}}^{**} 
    		\ar[r]^\cong &
    		\overline{\overline{V}} 
    		\ar[r]^\cong &
    		V,
    		}
    \end{equation}
    constitutes an $\Oone$-volutive structure on $\vectc$.
\end{example}

\medskip 

In the remainder of this subsection, we briefly review the theory of \cite{Stehouwer2023DaggerCV} relating the 2-categories $\Daggercat$ and $\Oonevolcat$. First, 
any dagger category $(\C,d)$ can be \emph{trivially} reinterpreted as an $\Oone$-volutive category $(\C,d,\id)$.
This assignment is in fact functorial, i.e.\ it extends to a 2-functor 
\begin{equation}
	\label{eq:TO1}
    T_{\Oone} \colon \Daggercat \xymatrix{\ar[r]&} \Oonevolcat.
\end{equation}
Conversely, any $\Oone$-volutive category $(\C,d,\eta)$ gives rise to a dagger category (which loosely speaking can be viewed as a ``strictification'' of $(\C,d,\eta)$) as follows. Define a category $S_{\Oone}\C$ whose objects are 
pairs $(a,\theta_a)$ consisting of an object $a \in \C$ and an isomorphism $\theta_a \colon a \to d(a)$ satisfying 
\begin{equation}
	\label{eq:CoherenceO1}
	\eta_a \circ d(\theta_a)^{-1} \circ \theta_a = \id_a, 
\end{equation}
and whose morphisms $(a,\theta_a) \to (b,\theta_b)$ are morphisms $a\to b$ in $\C$ without further properties. We define 
a dagger structure on $S_{\Oone}\C$ by assigning a morphism $X \colon (a,\theta_a) \to (b,\theta_b)$ to the morphism $(b,\theta_b) \to (a,\theta_a)$ given by 
\begin{equation}
	\label{eq:SO1onMorphisms}
    \xymatrix{
        b \ar[r]^-{\theta_b} & d(b) \ar[r]^-{d(X)} & d(a) \ar[r]^-{\theta_a^{-1}} & a.
    }
\end{equation}
One checks that this really defines a dagger structure by using the naturality of $\eta$ and our 
assumptions on $\theta_a,\theta_b$. This \emph{strictification\footnote{Note again that $T_{\Oone}$ and $S_{\Oone}$ are in general not fully faithful. We do not stress this point further.}} procedure in fact also extends to a 2-functor \cite[Lemma 3.13]{Stehouwer2023DaggerCV}
\begin{equation}
	\label{eq:SO1}
    S_{\Oone} \colon \Oonevolcat \xymatrix{\ar[r]&} \Daggercat.
\end{equation}

\begin{example}
	\label{exa:HermVect}
    The image of the $\Oone$-volutive category $(\vectc,\overline{(-)}^*)$ under the 2-functor $S_{\Oone}$ is the dagger category $\Herm_\mathbb{C}$ described in 
    \Cref{Example: The category of hermitian vector spaces}. The $\Oone$-volutive category obtained from the dagger category $\Herm_\mathbb{C}$ by applying $T_{\Oone}$ has the same 
    underlying category $\Herm_\mathbb{C}$ and the $\Oone$-volutive structure whose functor is given by the dagger structure on  $\Herm_\mathbb{C}$ and whose natural isomorphism 
    is trivial.
\end{example}

\begin{remark}
    The image of $T_{\Oone}$ lies in the full sub-2-category $\Oonevolcat^{\efix} \subset \Oonevolcat$ of those $\Oone$-volutive categories in which every object 
    $a$ admits a \emph{fixed point structure} as above, i.e.\ an isomorphism $\theta_a\colon a \to d(a)$ satisfying the coherence property~\eqref{eq:CoherenceO1}. 
    The image of $S_{\Oone}$ lies in the full sub-2-category $\Daggercat^{\pos} \subset \Daggercat$ 
    of those dagger categories, in which every self-adjoint automorphism is positive. 
\end{remark}

\begin{theorem}[{\cite[Theorem 4.9]{Stehouwer2023DaggerCV}}]\label{Theorem: The equivalence of strict and non-strict Oone volutive categories with restrictions}
    $S_{\Oone}$ is right adjoint to $T_{\Oone}$ in the 3-category of 2-categories. Moreover, the respective restrictions of $S_{\Oone}$ and $T_{\Oone}$ induce an equivalence 
    of 2-categories 
    \begin{equation}
    	\Oonevolcat^{\efix} \cong \Daggercat^{\pos}.
    \end{equation}
\end{theorem}
\begin{example}
	\label{Example: ST Herm = Herm}
    As a consequence of the theorem, we have $(S_{\Oone} \circ T_{\Oone})(\Herm_\mathbb{C}) \cong \Herm_\mathbb{C}$.
\end{example}

\begin{remark}
    Combining \Cref{Proposition: From symmetric monoidal categories to Oone volutive categories} with the 2-functor $S_{\Oone}$, we obtain a 2-functor 
    \begin{equation}
        \symmoncat^{\textrm{d}} \xymatrix{\ar[r]&} \Daggercat
    \end{equation} 
    assigning to each symmetric monoidal category with duals a dagger category.
\end{remark}

\subsection[$\rSpintwo$-dagger 2-categories]{$\textrm{\textbf{Spin}}\pmb{(2)^r}$-dagger 2-categories}
\label{Subsection: rspintwovolutive 2-categories}

In this section we follow a similar program to the one in \Cref{Subsection: Oonevolutive categories} with 1-categories replaced by 2-categories and 
$\Oone$ replaced by $\rSpintwo$. More precisely, we introduce 3-categories of (both strict and non-strict versions of) ``$\rSpintwo$-volutive 2-categories''.
We also describe an adjunction between both 3-categories, given by ``trivial re-interpretation'' and ``strictification''. 
Specific examples of $\rSpintwo$-dagger 2-categories are discussed in \Cref{sec:ExamplesApplications}.

\subsubsection{Volutive case}

$\Oone$-volutive categories can be interpreted as homotopy fixed points of an $\Oone$-action, cf.\ \Cref{Remark: Oone volutive categories as homotopy fixed points}. 
Analogously, there should be a $\rSpintwo$-action on the 3-category $\BiCatadj$ whose homotopy fixed points are $\rSpintwo$-volutive 2-categories. 
Recalling that $\BB\Pi_3(\rSpintwo)$ is equivalent to the 4-groupoid with precisely one object~$\star$, only identity 1-morphisms, and a $\Z$'s worth of 2-morphisms, 
to define such an action one in particular has to specify the image of a generator of $\Z$.%
\footnote{Recalling that $\BB^2 \mathbb{Z} \cong \BB \Sotwo \cong \mathbb{C}\mathbb{P}^{\infty}$,  in order to define 
a 4-functor $\rho \colon \BB\Pi_3(\rSpintwo) \to \TriCat$ with $\rho(\star) = \BiCatadj$ one also needs to specify the image of the 4-cell of $\mathbb{C}\mathbb{P}^{\infty}$. If such a 
4-functor is given, one may consider its homotopy fixed points, i.e.\ the 4-limit of~$\rho$ which, similar to the discussion in \Cref{Remark: Oone volutive categories as homotopy fixed points}, is expected to be $\Nat(\Delta,\rho)$ where $\Delta$ is the constant 4-functor $\BB\Pi_3(\rSpintwo) \to \TriCat$ whose image is the trivial 3-category, and $\Nat$ is the 3-category of 4-transformations 
between two 4-functors.} 
This image should be an invertible 3-transformation from the identity 3-functor on $\BiCatadj$ to itself. 
The 1-morphism component of this 3-transformation is generally expected to be the $2r$-fold power of the 2-functor $(-)^\vee\colon \B \to \B^{\onetwoop}$ which takes right adjoints of (2- and) 1-morphisms of $\B\in \BiCatadj$, see also \cite{StehouwerMueller2024}.

As the theory of 4-categories and 4-limits is not very well developed at the time of this writing, the above comments are more of a guiding principle, rather than a recipe for constructing 
the 3-category of $\rSpintwo$-volutive 2-categories. 
Motivated by examples, such as those in \Cref{sec:ExamplesApplications}, as well as we the work of \cite{StehouwerMueller2024, ferrer2024dagger}, we arrive at the following.
\begin{definition}
    Let $\B$ be a 2-category with right adjoints, and let $r \in \Z_{\geqslant 1}$. 
    A \emph{$\rSpintwo$-volutive} structure on $\B$ is an invertible 2-transformation 
    $S \colon \id_{\B} \Rightarrow ((-)^{\vee \vee})^r$. We call the pair $(\B,S)$ a \emph{$\rSpintwo$-volutive 2-category}.
\end{definition}

\begin{remark}
    Since for $r=1$ we have $\rSpintwo = \Sotwo$, we also speak of \textsl{$\Sotwo$-volutive} structures in this case. Similarly, for $r=2$, we may speak of 
    \textsl{$\Spintwo$-volutive} structures.
\end{remark}

\begin{definition}
    Let $(\B,S)$ and $(\B',S')$ be $\rSpintwo$-volutive 2-categories. A \emph{$\rSpintwo$-volutive 2-functor} $(F,\eta) \colon (\B,S) \to (\B',S')$ consists of a 2-functor 
    $F \colon \B \to \B'$ and an invertible modification 
    \begin{equation}
        \begin{aligned}
            \eta \colon 
            \begin{matrix}
            \xymatrix{
                & \ar@{}[d]^(.3){}="a"^(1){}="b" \ar@{=>}^{S} "a";"b" &  \\
                B \ar[rr]_-{((-)^{\vee \vee})^r} \ar[d]_-{F} \ar@/^2pc/[rr]^-{\id_{\B}} && B \ar[d]^-{F} \ar@{=>}[dll]^-{\cong} 
                \\
                B' \ar[rr]_-{((-)^{\vee \vee})^r} && B'
            }
            \end{matrix}
            \;\RRightarrow\;
            \begin{matrix}
                \xymatrix{
                B \ar[rr]^-{\id_{\B}} \ar[d]_-{F} && B \ar[d]^-{F} \ar[d]^-{F} \ar@{=>}[dll]_-{\cong}\\
                B' \ar[rr]^-{\id_{\B'}} \ar@/_2pc/[rr]_-{((-)^{\vee \vee})^r} & \ar@{}[d]^(0){}="a"^(0.7){}="b" \ar@{=>}^{S'} "a";"b" & B'\\
                & & 
            }
            \end{matrix}
        	.
        \end{aligned}
    \end{equation}
    Here we used the fact that every 2-functor automatically commutes with the adjunction 2-functor up to an invertible 2-transformation.
\end{definition}

\begin{definition}
    Let $(F,\eta), (\tilde{F},\tilde{\eta}) \colon (\B,S) \to (\B',S')$ be $\rSpintwo$-volutive 2-functors. A \emph{$\rSpintwo$-volutive 2-transformation} 
    $\alpha \colon (F,\eta) \to (\tilde{F},\tilde{\eta})$ is a 2-transformation $\alpha \colon F \Rightarrow \tilde{F}$ such that 
    \begin{equation}
        \begin{aligned}
            \begin{matrix}
                \xymatrix{
                    & F \circ ((-)^{\vee \vee})^r \ar@{=>}[dr]^-{\cong}  \ar@{}[d]^(.3){}="a"^(0.8){}="b" \ar@3{->}^{\eta} "a";"b" & \\
                    F \circ \id_{\B} \ar@{=>}[d]_-{\alpha \circ \id} \ar@{=>}[r]_-{\cong} \ar@{=>}[ur]^-{\id \circ S} & \id_{\B'} \circ F \ar@{=>}[r]_-{S' \circ \id} \ar@{=>}[d]_-{\id \circ \alpha} \ar@{}[dl]^(.35){}="a"^(0.75){}="b" \ar@3{->}^{\cong} "a";"b" & ((-)^{\vee \vee})^r \circ F \ar@{=>}[d]^-{\id \circ \alpha} \ar@{}[dl]^(.35){}="a"^(0.75){}="b" \ar@3{->}^{\cong} "a";"b"\\
                    \tilde{F} \circ \id_{\B} \ar@{=>}[r]_-{\cong} & \id_{\B'} \circ \tilde{F} \ar@{=>}[r]_-{S' \circ \id} & ((-)^{\vee \vee})^r \circ \tilde{F}
                }
            \end{matrix} 
            \;=\; 
            \begin{matrix}
                \xymatrix{
                    F \circ \id_{\B} \ar@{=>}[d]_-{\alpha \circ \id} \ar@{=>}[r]^-{\id \circ S} & F \circ ((-)^{\vee \vee})^r \ar@{=>}[r]^-{\cong} \ar@{=>}[d]_-{\alpha \circ \id} \ar@{}[dl]^(.35){}="a"^(0.75){}="b" \ar@3{->}_{\cong} "a";"b" & ((-)^{\vee \vee})^r \circ F \ar@{=>}[d]^-{\id \circ \alpha} \ar@{}[dl]^(.35){}="a"^(0.75){}="b" \ar@3{->}_{\cong} "a";"b"\\
                    \tilde{F} \circ \id_{\B} \ar@{=>}[r]^-{\id \circ S} \ar@{=>}[dr]_-{\cong} & \tilde{F} \circ ((-)^{\vee \vee})^r \ar@{=>}[r]^-{\cong} \ar@{}[d]^(.3){}="a"^(0.8){}="b" \ar@3{->}^{\tilde{\eta}} "a";"b"  & ((-)^{\vee \vee})^r \circ \tilde{F} \\
                    & \id_{\B'} \circ \tilde{F} \ar@{=>}[ur]_-{S' \circ \id}  & 
                }
            \end{matrix}
        	.
        \end{aligned}
    \end{equation}
    Here we used the interchange law as well as the naturality of the 2-transformations expressing compatibility of arbitrary 2-functors with the adjunction 2-functor.
\end{definition}

\begin{definition}
    Let $\alpha,\beta \colon (F,\eta) \to (\tilde{F},\tilde{\eta})$ be $\rSpintwo$-volutive 2-transformations. A \emph{$\rSpintwo$-volutive modification} 
    $\xi \colon \alpha \Rrightarrow \beta$ is a modification $\xi \colon \alpha \Rrightarrow \beta$.
\end{definition}
Let $(F,\eta) \colon (\B_1,S_1) \to (\B_2,S_2)$ and  $(F',\eta') \colon (\B_2,S_2) \to (\B_3,S_3)$ be $\rSpintwo$-volutive 2-functors. The \emph{composition} 
$(F',\eta') \circ (F,\eta)$ consists of the 2-functor $F' \circ F \colon \B_1 \to \B_3$ and the invertible modification 
\begin{equation} 
	\label{Equation: rspintwo chapter the composed invertible modification}
    \begin{aligned}
        \xymatrix{
            & F' \circ F \circ ((-)^{\vee \vee})^r \ar@{=>}[r]^-{\cong}  \ar@{}[dd]^(.3){}="a"^(0.7){}="b" \ar@3{->}_{\id_{F'} \circ \eta} "a";"b" & F' \circ ((-)^{\vee \vee})^r \circ F \ar@{=>}[dr]^-{\cong} \ar@{}[dd]^(.3){}="a"^(0.7){}="b" \ar@3{->}^{\eta' \circ \Id_F} "a";"b"  & \\
            F' \circ F \circ \id_{\B_1} \ar@{=>}[dr]_-{\cong} \ar@{=>}[ur]^-{\id \circ \id \circ S_1} &&& ((-)^{\vee \vee})^r \circ F' \circ F. \\
            &F' \circ \id_{\B_2} \circ F \ar@{=>}[r]_-{\cong} \ar@{=>}[uur]^-{\id \circ S_2 \circ \id} & \id_{\B_3} \circ F' \circ F \ar@{=>}[ur]_-{S_3 \circ \id \circ \id}&
        }
    \end{aligned}
\end{equation}
The composition of $\rSpintwo$-volutive 2-transformations is the usual composition of 2-transformations and similar for $\rSpintwo$-volutive modifications. 
It is straightforward to check that the above structures comprise a strict 3-category (see e.g.\ \cite{39} for background on 3-categories).
\begin{definition}
    The 3-category of $\rSpintwo$-volutive 2-categories, $\rSpintwo$-volutive 2-functors, $\rSpintwo$-volutive 2-transformations, $\rSpintwo$-volutive modifications and compositions as above is denoted $\rSpintwovoltwocat$.
\end{definition}
\begin{remark}
    We stress again that $\rSpintwovoltwocat$ is only conjecturally the 3-category of homotopy fixed points of a $\rSpintwo$-action on the 3-category $\BiCatadj$. We do not 
    claim completeness of the definitions above.
\end{remark}
\begin{example} \label{Example: Symmetric monoidal 2-categories with duals and adjoints are rSpintwo volutive}
    Let $\B$ be a symmetric monoidal 2-category with duals $(-)^*$ and adjoints~$(-)^\vee$. 
    For symmetric monoidal structures on 2-categories we refer to \cite{26,Coalgebroids2000BalancedCP}, while for 
    duals and adjoints in monoidal 2-categories we refer to \cite{84} for more details. The \emph{Serre automorphism} of an object $a \in \B$ is the 1-automorphism 
    \begin{equation}
    	\label{eq:SerreAutomorphismSa}
        S_a := (\id \otimes \ev_a) \circ (\beta_{a,a} \otimes \id) \circ (\id \otimes \ev^\vee_a) \colon a \xymatrix{\ar[r]&} a
    \end{equation}
	where~$\beta$ denotes the braiding. 
    Extending results of \cite{lurie2009classificationtopologicalfieldtheories} (restricting to the maximal subgroupoid of $\B$) and \cite{75} (assuming the presence 
    of a pivotal structure), there is an invertible 2-transformation $\mathbb{S} \colon \id_{\B} \Rightarrow (-)^{\vee \vee}$ whose 1-morphism components are given by the Serre automorphisms~\eqref{eq:SerreAutomorphismSa}. 
    The construction of the 2-morphism components of $\mathbb{S}$ is formally analogous to the graphical one in \cite[Proposition 3.2]{75} if one omits the pivotal structure or the one given in \cite[Proposition 4.9]{82} if one omits the invertibility assumption and keeps track of adjoints.
    Concretely, for every 1-morphism $X\colon a\to a'$ the associated 2-morphism component of~$\mathbb{S}$ is 
    \begin{equation}
    	S_X = 
    	\begin{tikzpicture}[thick,scale=1.0,color=black, baseline=2.7cm]
    		\coordinate (p1) at (0,0);
    		\coordinate (p2) at (2,-0.5);
    		\coordinate (p3) at (2,0.5);
    		\coordinate (p4) at (4,0);
    		\coordinate (u1) at (0,4.5);
    		\coordinate (u2) at (2,4);
    		\coordinate (u3) at (2,5);
    		\coordinate (u4) at (4,4.5);
    		\coordinate (ld) at (-2,1);
    		\coordinate (lu) at (-2,5.5);
    		\coordinate (rd) at (6,1);
    		\coordinate (ru) at (6,5.5);
    		%
    		\fill [orange!20!white, opacity=0.8] (p3) -- (ld) -- (lu) -- (u3); 
    		\fill [orange!20!white, opacity=0.8] (p3) -- (rd) -- (ru) -- (u3); 
    		%
    		\fill [orange!25!white, opacity=0.8] 
    		(p1) .. controls +(0,0.25) and +(-1,-0.2) ..  (p3)
    		-- (p3) .. controls +(1,-0.2) and +(0,0.25) ..  (p4)
    		-- (p4) --  (u4)
    		-- (u4) .. controls +(0,0.25) and +(1,-0.2) ..  (u3)
    		-- (u3) .. controls +(-1,-0.2) and +(0,0.25) ..  (u1)
    		;
    		%
    		\coordinate (x2) at (2,1.59);
    		\draw[ultra thick] (p3) --  (u3); 
    		\fill (2.5,1.35) circle (0pt) node {{\small $\beta_{1_a,X}$}};
    		\fill[color=blue!50!black] (x2) circle (2.5pt) node {};
    		\coordinate (x2s) at (2,3.28);
    		\fill[color=blue!50!black] (x2s) circle (2.5pt) node {};
    		\fill (2.6,3) circle (0pt) node {{\small $\beta_{{}^{\vee\vee}\!X,1_{a'}}$}};
    		%
    		\coordinate (x1) at (-1.5,0.93);
    		\coordinate (x3) at (4,2);
    		\coordinate (x4) at (0,2.7);
    		\coordinate (x5) at (5.5,5.44);
    		\draw[ultra thick, blue!50!black] (x1) .. controls +(0,0.7) and +(0,-0.5) .. (x3); 
    		\draw[ultra thick, blue!50!black] (x4) .. controls +(0,0.75) and +(0,-2.5) .. (x5);
    		\fill[color=blue!50!black] (-1.6,1.2) circle (0pt) node {{\small $X$}};
    		\fill[color=blue!50!black] (5.1,5.1) circle (0pt) node {{\small ${}^{\vee\vee}\!X$}};
    		%
    		\draw[very thick, red!80!black] (p1) .. controls +(0,0.25) and +(-1,-0.2) ..  (p3); 
    		\draw[very thick, red!80!black] (p4) .. controls +(0,0.25) and +(1,-0.2) ..  (p3); 
    		\draw[very thick, red!80!black] (p3) -- (ld);
    		\draw[very thick, red!80!black] (u3) -- (lu); 
    		\draw[very thick, red!80!black] (p3) -- (rd);
    		\draw[very thick, red!80!black] (ru) -- (u3); 
    		\draw[thin] (lu) --  (ld); 
    		\draw[thin] (ru) --  (rd); 
    		%
    		\fill [orange!30!white, opacity=0.8] 
    		(p1) .. controls +(0,-0.25) and +(-1,0) ..  (p2)
    		-- (p2) .. controls +(1,0) and +(0,-0.25) ..  (p4)
    		-- (p4) --  (u4)
    		-- (u4) .. controls +(0,-0.25) and +(1,0) ..  (u2)
    		-- (u2) .. controls +(-1,0) and +(0,-0.25) ..  (u1)
    		;
    		\draw[thin] (p1) --  (u1); 
    		\draw[thin] (p4) --  (u4); 
    		%
    		\draw[ultra thick, blue!50!black] (x3) .. controls +(0,0.5) and +(0,-0.5) .. (x4); 
    		%
    		\draw[very thick, red!80!black] (p1) .. controls +(0,-0.25) and +(-1,0) ..  (p2)
    		--(p2) .. controls +(1,0) and +(0,-0.25) .. (p4); 
    		\draw[very thick, red!80!black] (u1) .. controls +(0,-0.25) and +(-1,0) ..  (u2) 
    		--(u2) .. controls +(1,0) and +(0,-0.25) .. (u4); 
    		\draw[very thick, red!80!black] (u4) .. controls +(0,0.25) and +(1,-0.2) ..  (u3); 
    		\draw[very thick, red!80!black] (u1) .. controls +(0,0.25) and +(-1,-0.2) .. (u3); 
    		%
    		\fill (2,0) circle (0pt) node {{\small $a^*$}};
    		\fill (1.3,4.4) circle (0pt) node {{\small $a'$}};
    		\fill (3.6,2.8) circle (0pt) node {{\small $a'^*$}};
    		\fill (5.5,1.5) circle (0pt) node {{\small $a$}};
    		\fill (-0.5,1.1) circle (0pt) node {{\small $a$}};
    		\fill (4.3,4.9) circle (0pt) node {{\small $a'$}};
    		\fill (-1.5,4) circle (0pt) node {{\small $a'$}};
    		\fill[color=blue!50!black] (1.11,2.1) circle (0pt) node {{\small $({}^{\vee\vee}\!X)^*$}};
    		\fill[color=blue!50!black] (4.6,1.95) circle (0pt) node {{\small $\widetilde\Omega^\vee_{({}^{\vee\vee}\!X)^\vee}$}};
    		\fill[color=blue!50!black] (-0.4,2.7) circle (0pt) node {{\small $\widetilde\Omega_{{}^{\vee\vee}\!X}$}};
    		%
    		%
    	\end{tikzpicture}
    	\,\colon X \circ S_a \xymatrix{\ar[r]&} S_{a'} \circ {}^{\vee\vee}\!X
    \end{equation}
	where $\widetilde\Omega_X\colon 1_a \otimes X^* \to X \otimes 1_{a'^*}$ is the 2-isomorphism determined by the cusp isomorphism which is part of the duality data of $a'\in\B$, and we identify left and right duals thanks to the braided structure. 
	Here we employ the graphical calculus of \cite{Barrett2012} with the conventions of \cite{Brunner2023TruncatedAR}; in particular, we have $\widetilde\Omega_X = \Omega_{X^*}$ for~$\Omega$ as in \cite[Equation~(2.24)]{Brunner2023TruncatedAR}. 
	Thus~$S_X$ is manifestly invertible. 
    
    The (appropriately defined) $r$-th power of $\mathbb{S}$ gives rise to an invertible 2-transformation $\mathbb{S}^r \colon \id_{\B} \Rightarrow ((-)^{\vee \vee})^r$ whose 1-morphism components are given by the $r$-th powers of the Serre automorphisms. 
\end{example}

The following is a subexample of \Cref{Example: Symmetric monoidal 2-categories with duals and adjoints are rSpintwo volutive}.

\begin{example} 
	\label{Example: The delooping of vector spaces is strict rSpintwo volutive}
    Recall that the category $\vect$ is symmetric monoidal and rigid. In consequence, $\mathrm{B}\vect$ also carries a symmetric monoidal structure; moreover, 
    $\mathrm{B}\vect$ has duals and adjoints. Here the former are trivial while the latter are given by the dual vector space construction. Computing the 
    invertible 2-transformation~$\mathbb{S}$ described in \Cref{Example: Symmetric monoidal 2-categories with duals and adjoints are rSpintwo volutive} yields the following. 
    Since there is only one object in $\mathrm{B}\vect$, the 1-morphism component of $\mathbb{S}$ is trivial. The 2-morphism component of $\mathbb{S}$ at a 1-morphism 
    $V \colon \star \to \star$ is given by the canonical isomorphism $V \to V^{**}$ from $V$ to its double dual. 
    Similarly, the invertible 2-transformation $\mathbb{S}^r$ 
    has trivial 1-morphism components and 2-morphism components given by iterated trivialization of the double dual. We note that the $\rSpintwo$-volutive 
    structure on $\mathrm{B}\vect$ described here is in fact strict in the sense of \Cref{Definition: strict rSpintwo volutive 2-category}. 
\end{example}

\begin{remark} 
	\label{Remark: The 3-functor from symmontwocats to tspintwovoltwocats}
    We expect that \Cref{Example: Symmetric monoidal 2-categories with duals and adjoints are rSpintwo volutive} can be extended to a 3-functor 
    \begin{equation}
        \symmontwocat^{\textrm{d}} \xymatrix{\ar[r]&} \rSpintwovoltwocat
    \end{equation}
    from the 3-category of symmetric monoidal 2-categories with duals and adjoints to the 3-category of $\rSpintwo$-volutive 2-categories. The compatibility of the 
    2-transformation $\mathbb{S}$ with symmetric monoidal 2-functors is proved in \cite[Proposition B.12]{müller2023reflection} for the maximal subgroupoid. An extension of 
    this argument should then yield the assignment of the 3-functor on 1-morphism level, while for the assignment on 2-morphism level one only has to check a condition, and 
    on 3-morphism level there is nothing to show.
\end{remark}

\subsubsection{Dagger case}

In the following, we will introduce the appropriate strict versions of $\rSpintwo$-volutive 2-categories and their higher morphisms, and relate them to pivotal structures. 

\begin{definition} \label{Definition: strict rSpintwo volutive 2-category}
    Let $\B$ be a 2-category with right adjoints, and let $r \in \Z_{\geqslant 1}$. A \emph{$\rSpintwo$-dagger} structure on $\B$ is an invertible 2-transformation 
    $S \colon \id_{\B} \Rightarrow ((-)^{\vee \vee})^r$ whose 1-morphism components are identities, i.e.\ $S_a = \id_a \colon a \to a$. We call the pair $(\B,S)$ a \emph{$\rSpintwo$-dagger 2-category}.
\end{definition}
Note that the 2-morphism components of a $\rSpintwo$-dagger structure need not be identities. Any $\rSpintwo$-dagger 2-category $(\B,S)$ can be trivially 
reinterpreted as a $\rSpintwo$-volutive 2-category $(\B,S)$. We will come back to this after defining appropriate notions of (higher) morphisms between $\rSpintwo$-dagger 2-categories. 

\begin{remark} \label{Remark: Equivalent ways of talking about pivotal structures}
    There is another, equivalent way of thinking about $\rSpintwo$-dagger structures on a 2-category with right adjoints, namely, in terms of 
    left adjoints. Let $X \colon a \to b$ be a 1-morphism in a $\rSpintwo$-dagger 2-category $(\B,S)$. We can define the left adjoint of $X$ to be 
    $^\vee{X} := ((-)^{\vee})^{2r-1}(X) $ with left adjunction 2-morphisms
    \begin{align}
        \ev_X \;=\;
        \Big(
        \xymatrix{
            ^\vee{X} \circ X \ar[rr]^-{\id \circ S_X^{-1}} &&  ((-)^{\vee})^{2r-1}(X) \circ  ((-)^{\vee})^{2r}(X) 
            \ar[rr]^-{\widetilde{\ev}_{((-)^{\vee})^{2r-1}(X)}} && \id_a
        }
    	\Big)
    \end{align}
    and  
    \begin{align}
    	\coev_X \;=\;
    	\Big( 
        \xymatrix{
            \id_b \ar[rr]^-{\widetilde{\coev}_{((-)^{\vee})^{2r-1}(X)}} && ((-)^{\vee})^{2r}(X) \circ  ((-)^{\vee})^{2r-1}(X) \ar[rr]^-{S_X \circ \id} && X \circ ^\vee{X}
        }
    	\Big)
    \end{align}
    where $\widetilde{\ev}$ and $\widetilde{\coev}$ denote right evaluation and coevaluation maps, respectively. One verifies that $\ev_X$ and $\coev_X$ 
    satisfy the Zorro moves. 
    In the case $r=1$ we have $^\vee{X} = X^\vee$, hence an $\Sotwo$-dagger structure is precisely the same as a pivotal structure on a 
    2-category with right adjoints, see e.g.\ \cite[Section 2.1]{carqueville2023orbifold} and \cite[Section 2.3]{Carqueville2012}.
\end{remark}

\begin{definition}
	\label{def:StrictSpinVolutove2Functor}
    Let $(\B,S)$ and $(\B',S')$ be $\rSpintwo$-dagger 2-categories. A \emph{$\rSpintwo$-dagger 2-functor} $F \colon (\B,S) \to (\B',S')$ consists of a 2-functor 
    $F \colon \B \to \B'$ such that
    \begin{equation}
        \begin{aligned}
            \begin{matrix}
            \xymatrix{
                & \ar@{}[d]^(.3){}="a"^(1){}="b" \ar@{=>}^{S} "a";"b" &  \\
                B \ar[rr]_-{((-)^{\vee \vee})^r} \ar[d]_-{F} \ar@/^2pc/[rr]^-{\id_{\B}} && B \ar[d]^-{F} \ar@{=>}[dll]^-{\cong} \\
                B' \ar[rr]_-{((-)^{\vee \vee})^r} && B'
            }
            \end{matrix}
            \;\;\;=\;\;\;
            \begin{matrix}
                \xymatrix{
                B \ar[rr]^-{\id_{\B}} \ar[d]_-{F} && B \ar[d]^-{F} \ar[d]^-{F} \ar@{=>}[dll]_-{\cong}\\
                B' \ar[rr]^-{\id_{\B'}} \ar@/_2pc/[rr]_-{((-)^{\vee \vee})^r} & \ar@{}[d]^(0){}="a"^(0.7){}="b" \ar@{=>}^{S'} "a";"b" & B'\\
                & & 
            }
            \end{matrix}
        	.
        \end{aligned}
    \end{equation}
\end{definition}
\begin{definition}
    Let $F,\tilde{F}\colon (\B,S) \to (\B',S')$ be $\rSpintwo$-dagger 2-functors. A \emph{$\rSpintwo$-dagger 2-transformation} 
    $\alpha \colon F \to \tilde{F}$ is a 2-transformation $\alpha \colon F \Rightarrow \tilde{F}$ such that 
    \begin{equation}
        \begin{aligned}
            \begin{matrix}
                \xymatrix{
                    F \circ \id_{\B} \ar@{=>}[d]_-{\alpha \circ \id} \ar@{=>}[r]^-{\cong} & \id_{\B'} \circ F \ar@{=>}[r]^-{S' \circ \id} \ar@{=>}[d]_-{\id \circ \alpha} \ar@{}[dl]^(.35){}="a"^(0.75){}="b" \ar@3{->}_{\cong} "a";"b" & ((-)^{\vee \vee})^r \circ F \ar@{=>}[d]^-{\id \circ \alpha} \ar@{}[dl]^(.35){}="a"^(0.75){}="b" \ar@3{->}_{\cong} "a";"b"\\
                    \tilde{F} \circ \id_{\B} \ar@{=>}[r]_-{\cong} & \id_{\B'} \circ \tilde{F} \ar@{=>}[r]_-{S' \circ \id} & ((-)^{\vee \vee})^r \circ \tilde{F}
                }
            \end{matrix} 
            \;\;\; =\;\;\; 
            \begin{matrix}
                \xymatrix{
                    F \circ \id_{\B} \ar@{=>}[d]_-{\alpha \circ \id} \ar@{=>}[r]^-{\id \circ S} & F \circ ((-)^{\vee \vee})^r \ar@{=>}[r]^-{\cong} \ar@{=>}[d]_-{\alpha \circ \id} \ar@{}[dl]^(.35){}="a"^(0.75){}="b" \ar@3{->}_{\cong} "a";"b" & ((-)^{\vee \vee})^r \circ F \ar@{=>}[d]^-{\id \circ \alpha} \ar@{}[dl]^(.35){}="a"^(0.75){}="b" \ar@3{->}_{\cong} "a";"b"\\
                    \tilde{F} \circ \id_{\B} \ar@{=>}[r]_-{\id \circ S} & \tilde{F} \circ ((-)^{\vee \vee})^r \ar@{=>}[r]_-{\cong} & ((-)^{\vee \vee})^r \circ \tilde{F} \\
                }
            \end{matrix}
        	.
        \end{aligned}
    \end{equation}
\end{definition}
\begin{definition}
    Let $\alpha,\beta \colon F \to \tilde{F}$ be $\rSpintwo$-dagger 2-transformations. A \emph{$\rSpintwo$-dagger modification} 
    $\xi \colon \alpha \to \beta$ is a modification $\xi \colon \alpha \Rrightarrow \beta$.
\end{definition}
One checks that the above structures together with ordinary compositions comprise a 3-category:
\begin{definition}
    The 3-category of $\rSpintwo$-dagger 2-categories, $\rSpintwo$-dagger 2-functors, $\rSpintwo$-dagger 2-transformations, 
    and $\rSpintwo$-dagger modifications is denoted $\strictrSpintwovoltwocat$.
\end{definition}

\subsubsection{Strictification}

Clearly, any $\rSpintwo$-dagger 2-category $(\B,S)$ can be \emph{trivially} reinterpreted as a $\rSpintwo$-volutive 2-category $(\B,S)$. Similarly, any 
$\rSpintwo$-dagger 2-functor $F$ between $\rSpintwo$-dagger 2-categories can be reinterpreted as a $\rSpintwo$-volutive 2-functor $(F,\id)$ in which 
the modification $\eta$ is trivial. Extending this argument and recalling our definition of the composition of $\rSpintwo$-volutive 2-functors, 
one finds that this assignment is in fact functorial, i.e.\ it extends to a 3-functor
\begin{equation}
    T_{\rSpintwo} \colon \strictrSpintwovoltwocat \xymatrix{\ar[r]&} \rSpintwovoltwocat.
\end{equation}
Conversely, any $\rSpintwo$-volutive 2-category $(\B,S)$ gives rise to a $\rSpintwo$-dagger 2-category as follows. Define a 2-category $S_{\rSpintwo}\B$ whose objects 
are pairs $(a,\lambda_a)$ where $a \in \B$ and $\lambda_a \colon S_a \to \id_a$ is a 2-isomorphism, whose 1-morphisms $X \colon (a,\lambda_a) \to (b,\lambda_b)$ are 1-morphisms 
$X \colon a \to b$ in $\B$, and whose 2-morphisms $f \colon X \to Y$ are 2-morphisms in $\B$. We also define a $\rSpintwo$-dagger structure $S_{\rSpintwo}S$ on $S_{\rSpintwo}\B$ 
as follows. The 2-morphism component of $S_{\rSpintwo}S$ at $X \colon (a,\lambda_a) \to (b,\lambda_b)$ is given by the composition
\begin{equation} \label{Equation: strictifying a strict rspintwo volution}
    \xymatrix{
        X^{\vee \vee} \cong X^{\vee \vee} \circ \id_a \ar[r]^-{\id \circ \lambda_a^{-1}} &X^{\vee \vee} \circ S_a \ar[r]^-{S_X} & S_b \circ X \ar[r]^-{\lambda_b \circ \id} & \id_b \circ X \cong X. 
    }
\end{equation}
One checks that this defines a $\rSpintwo$-dagger structure on $S_{\rSpintwo}\B$. 
This is analogous to the discussion leading up to~\eqref{eq:SO1onMorphisms}, and analogously we have: 

\begin{lemma}
    The assignment $(\B,S) \mapsto (S_{\rSpintwo}\B,S_{\rSpintwo}S)$ extends to a 3-functor 
    \begin{equation}
        S_{\rSpintwo} \colon \rSpintwovoltwocat \xymatrix{\ar[r]&} \strictrSpintwovoltwocat.
    \end{equation}
    between the 3-category of $\rSpintwo$-volutive 2-categories and the 3-category of $\rSpintwo$-dagger 2-categories.
\end{lemma}

\begin{proof}
    Let $(F,\eta) \colon (\B,S) \to (\B',S')$ be a $\rSpintwo$-volutive 2-functor. 
    We construct a $\rSpintwo$-dagger 2-functor
    $S_{\rSpintwo}F \colon (S_{\rSpintwo}\B,S_{\rSpintwo}S) \to (S_{\rSpintwo}\B',S_{\rSpintwo}S')$ by sending an object $(a,\lambda_a)$ to the object consisting 
    of $F(a) \in \B'$ and the 2-isomorphism 
    \begin{equation}
        \xymatrix{
            S'_{F(a)} \ar[r]^-{\eta_a^{-1}} & F(S_a) \ar[r]^-{F(\lambda_a)} & F(\id_a) \cong \id_{F(a)}
        }
    \end{equation}
    and 1- and 2-morphisms to their images under $F$. 
    To check that this is functorial, let 
    $(F,\eta) \colon (\B_1,S_1) \to (\B_2,S_2)$ and $(F',\eta') \colon (\B_2,S_2) \to (\B_3,S_3)$ be $\rSpintwo$-volutive 2-functors. 
    Noting that we only have 
    to consider the respective assignments on object level, we compute 
    \begin{equation}
        (S_{\rSpintwo}F' \circ S_{\rSpintwo}F)(a,\lambda_a) = S_{\rSpintwo}F'(F(a), F(\lambda_a) \circ \eta_a^{-1}) = (F'F(a),F'(F(\lambda_a) \circ \eta_a^{-1}) \circ {\eta'}_{F(a)}^{-1})
    \end{equation}
    and recognize that $F'(\eta_a^{-1}) \circ {\eta'}_{F(a)}^{-1}$ is precisely the component of the modification defined in~\eqref{Equation: rspintwo chapter the composed invertible modification} 
    at $a$, which yields the claim. 
    Describing the assignment on 2- and 3-morphism level and proving functoriality there is simpler, and we omit the details.
\end{proof}

Recall that in order to describe an adjunction in a higher category, one has to give unit and counit, as well as higher isomorphisms that trivialize the Zorro compositions of unit and counit. 
In the following, the latter will be particularly simple, which allows us to carry out the argument rigorously. 
(For adjunctions in higher categories we refer to \cite[Chapter~2]{RiehlVerity2022}.)

\begin{proposition}
	\label{prop:STadjointSpin}
    $S_{\rSpintwo}$ is right adjoint to $T_{\rSpintwo}$ in the 4-category of 3-categories. 
\end{proposition}

\begin{proof}
	We divide the argument into four parts. 
	
	\noindent 
    Step I: We start by constructing the unit of the adjunction, which is a 3-transformation 
    \begin{equation}
        U \colon \id_{\strictrSpintwovoltwocat} \xymatrix{\ar[r]&} S_{\rSpintwo} \circ T_{\rSpintwo}.
    \end{equation}
    As a 3-transformation, $U$ has three layers of structure; we start by describing the 1-morphism component of $U$ at some strict 
    $\rSpintwo$-volutive 2-category $(\B,S)$, which is a $\rSpintwo$-dagger 2-functor $(\B,S) \to (S_{\rSpintwo} \circ T_{\rSpintwo})(\B,S)$. 
    Consider the 2-functor $U_{(\B,S)} \colon \B \to (S_{\rSpintwo} \circ T_{\rSpintwo})\B$ with $a \mapsto (a,\id_{\id_a})$, $X \mapsto X$ and $f \mapsto f$ on objects, 1- and 2-morphisms, respectively. We claim that this 
    2-functor is $\rSpintwo$-dagger. Indeed, this is true because the $\rSpintwo$-dagger structure on $(S_{\rSpintwo} \circ T_{\rSpintwo})\B$ 
    has components given by~\eqref{Equation: strictifying a strict rspintwo volution} which simplifies on the image of our 2-functor to the $\rSpintwo$-dagger
    structure we started with.
    
    Next, we want to describe the 2-morphism component of $U$ at some $\rSpintwo$-dagger 2-functor $F \colon (\B,S) \to (\B',S')$, which is a 
    $\rSpintwo$-dagger 2-transformation $(S_{\rSpintwo} \circ T_{\rSpintwo})(F) \circ U_{(\B,S)} \to U_{(\B',S')} \circ F$. To do so, we spell out the
    assignments of the source and the target of our 2-transformation on object level. Explicitly, we have 
    \begin{equation}
        (((S_{\rSpintwo} \circ T_{\rSpintwo})(F)) \circ U_{(\B,S)})(a) = ((S_{\rSpintwo} \circ T_{\rSpintwo})(F)) (a,\id_{\id_a}) = (F(a),\id_{\id_{F(a)}})
    \end{equation}
    and 
    \begin{equation}
        (U_{(\B',S')} \circ F)(a) = U_{(\B',S')}(F(a)) = (F(a),\id_{\id_{F(a)}}).
    \end{equation}
    This together with the fact that the assignments of the source and target on 1- and 2-morphism level coincide by definition allows us to define the 
    2-transformation $U_{(\B,S)}$ to be trivial, hence $\rSpintwo$-dagger. By a similar argument, we may set the 3-morphism component of $U$ also to be trivial. 
    This finishes the construction of the unit of the adjunction. 
    
    \medskip 

	\noindent 
    Step II: The next piece of the adjunction is the counit, which is a 3-transformation 
    \begin{equation}
        K \colon T_{\rSpintwo} \circ S_{\rSpintwo} \xymatrix{\ar[r]&} \id_{\rSpintwovoltwocat}.
    \end{equation}
    Again, $K$ has three levels of structure; we start by describing the 1-morphism component of $K$ at some $\rSpintwo$-volutive 2-category $(\B,S)$, 
    which is a $\rSpintwo$-volutive 2-functor $(T_{\rSpintwo} \circ S_{\rSpintwo})(\B,S) \to (\B,S)$. Consider the 2-functor 
    $K_{(\B,S)} \colon (T_{\rSpintwo} \circ S_{\rSpintwo})(\B) \to \B$ with $(a,\lambda_a) \mapsto a$, $X \mapsto X$ and $f \mapsto f$. Spelling out the $\rSpintwo$-volutive structure 
    of $(T_{\rSpintwo} \circ S_{\rSpintwo})(\B,S)$, which has identity 1-morphism components, we find that, in order to enhance $K_{(\B,S)}$ 
    to a $\rSpintwo$-volutive 2-functor, we need to provide natural 2-isomorphisms $\id_a \to S_a$ for each $(a,\lambda_a)\in (T_{\rSpintwo} \circ S_{\rSpintwo})(\B)$. The evident choice for such 2-isomorphisms 
    is $\lambda_a^{-1} \colon \id_a \to S_a$, which is in fact natural. The functor $K_{(\B,S)}$ together with the invertible modification 
    $\Lambda_{(\B,S)}$ whose 2-morphism component at $(a,\lambda_a)$ is $\lambda_a^{-1}$ then forms a $\rSpintwo$-volutive 2-functor. 
    
    By a similar argument as in the construction of the unit of the adjunction, we find that the 2- and 3-morphism component of $K$ can be chosen to be trivial. 
    This finishes the construction of the counit of the adjunction. 
	
	\medskip 
	
	\noindent 
    Step III: We compute the Zorro compositions of counit and unit. Consider first the composition 
    \begin{equation}
        \xymatrix{
            T_{\rSpintwo} \ar[r]^-{\id \circ U} & T_{\rSpintwo} \circ S_{\rSpintwo} \circ T_{\rSpintwo} \ar[r]^-{K \circ \id} & T_{\rSpintwo}
        }
    \end{equation}
    which is again a 3-transformation. As argued before, only the 1-morphism components of this 3-transformation are of interest here. In particular, we 
    wish to compute the 1-morphism component of $(K \circ \id) \bullet (\id \circ U)$ at some $\rSpintwo$-dagger 2-category $(B,S)$, where ``$\bullet$'' denotes the horizontal composition of 3-transformations. This requires us to study the 
    composition of $\rSpintwo$-dagger 2-functors 
    \begin{equation}
        \xymatrix{
            T_{\rSpintwo}(\B,S) \ar[rr]^-{T(U_{(\B,S)})} && (T_{\rSpintwo} \circ S_{\rSpintwo} \circ T_{\rSpintwo})(\B,S) \ar[rr]^-{K_{T_{\rSpintwo}(B,S)}} && T_{\rSpintwo}(\B,S).
        }
    \end{equation}
    On object level, we have $(K_{T_{\rSpintwo}(\B,S)} \circ T(U_{(\B,S)}))(a) = K_{T_{\rSpintwo}(\B,S)} (a,\id_{\id_a}) = a$, while on 2- and 3-morphism level the assignment 
    is again trivial. This shows that the 1-morphism component of the composition $(K \circ \id) \bullet (\id \circ U)$ and hence the whole 3-transformation is trivial. 
    
    Consider now the second Zorro composition 
    \begin{equation}
        \xymatrix{
            S_{\rSpintwo} \ar[r]^-{U \circ \id} & S_{\rSpintwo} \circ T_{\rSpintwo} \circ S_{\rSpintwo} \ar[r]^-{\id \circ K} & S_{\rSpintwo}
        }
    \end{equation}
    of which we wish to compute the 1-morphism component at some $\rSpintwo$-volutive 2-category $(\B,S)$. On object level, we have
    $(S_{\rSpintwo}(K_{(\B,S)}) \circ U_{S_{\rSpintwo}(\B,S)}) (a,\lambda_a) = (S_{\rSpintwo}(K_{(\B,S)}))(a,\lambda_a,\id_{\id_{(a,\lambda_a)}}) = (a,\lambda_a)$. 
    The second step requires slightly more elaboration. 
    Namely, we note that the assignment of the 1-morphism component of $\id \circ K$ on object level is given by $(a,\lambda_a,\psi_a) \mapsto (a,\psi_a \circ \lambda_a)$ where $\lambda_a \colon S_a \to \id_a$ and $\psi_a \colon \id_a = \id_{(a,\lambda_a)} \to \id_{(a,\lambda_a)} = \id_a$ 
    are 2-isomorphisms. 
    We conclude that the 2-functor underlying the 1-morphism component of $(\id \circ K)\bullet (U \circ \id)$ is trivial. It remains to check that the 
    invertible modification underlying the 1-morphism component of $(\id \circ K)\bullet (U \circ \id)$ at $(\B,S)$ is also trivial. 
    This is also true, by definition of the invertible modification underlying the 1-morphism component of $(\id \circ K)_{(\B,S)}$ restricted to the image of $U_{S_{\rSpintwo}(\B,S)}$, as well as the triviality of the invertible modification underlying $(U \circ \id)_{(\B,S)}$. 
    This shows that the whole 3-transformation is trivial. 
	
	\medskip 
	
	\noindent 
    Step IV: Noting that both Zorro compositions are strictly trivial, we may complete the given data to an adjunction in the 4-category of 3-category 
    by chosing trivial higher coherence morphisms. This finishes the proof.
\end{proof}
\begin{remark}
    Under the assumption that there is a 3-functor $\symmontwocat^{\textrm{d}} \to \rSpintwovoltwocat$ as in \Cref{Remark: The 3-functor from symmontwocats to tspintwovoltwocats}, post-composition with the 3-functor $S_{\rSpintwo}$
    yields a 3-functor 
    \begin{equation}
        \symmontwocat^{\textrm{d}} \xymatrix{\ar[r]&} \strictrSpintwovoltwocat
    \end{equation}
    from the 3-category of symmetric monoidal 2-categories with duals and adjoints to the 3-category of $\rSpintwo$-dagger 2-categories. 
\end{remark}
\begin{remark}
    There is a natural notion of the full sub-3-category 
    \begin{equation}
    	\rSpintwovoltwocat^{\efix} \; \subset \; \rSpintwovoltwocat
    \end{equation}
    of those $\rSpintwo$-volutive 2-categories $(\B,S)$ which admit 2-isomorphisms 
    $\lambda_a \colon S_a \to \id_a$ for all objects $a \in \B$. Clearly, the essential image of $T_{\rSpintwo}$ is contained in $\rSpintwovoltwocat^{\efix}$. A 
    notion of positivity in $\strictrSpintwovoltwocat$ analogous to the $\Oone$-case is less obvious. 
\end{remark}

\subsection[$\Oone$-dagger 2-categories]{$\textrm{\textbf{O}}\pmb{(1)}$-dagger 2-categories}
\label{Subsection: Oonevolutive 2-categories}

In this section we follow a similar program to the one in \Cref{Subsection: rspintwovolutive 2-categories} with $\rSpintwo$ replaced by $\Oone$. 
We give some (but not all) details on the associated 3-category of ``$\Oone$-volutive 2-categories'', as well as the 3-category of $\Oone$-dagger 2-categories, 
and on a strictification 3-functor between them; see also \cite{Stehouwer2023DaggerCV,ferrer2024dagger}.

\subsubsection{Volutive case}

\begin{definition}
	\label{def:O1volutive2Category}
    Let $\B$ be a 2-category. 
    An \emph{$\Oone$-volution} on $\B$ is a triple $(d,\eta,\tau)$ consisting of a 2-functor $d \colon \B \to \B^{\oneop}$, 
    an invertible 2-transformation $\eta \colon d^{\oneop} \circ d \Rightarrow \id_{\B}$, and an invertible modification $\tau \colon (\eta^{\oneop})^{-1} \circ \id_d \Rrightarrow \id_{d} \circ \eta$ 
    such that the following equality of modifications holds:
    \begin{equation} \label{anti-involution bicategory coherence diagram}
        \begin{matrix}
            \xymatrix{
                & d^{\oneop} \circ d \ar@2{->}[dr]^-{\id} \ar@2{->}@/^1.5pc/[drr]^-{\eta} \ar@2{->}@{}[d]^(.25){}="a"^(.75){}="b" \ar@3{->}_{(\tau^{\oneop})^{-1} \circ \id} "a";"b" & \ar@2{->}@{}[d]^(.25){}="a"^(.75){}="b" \ar@3{->}^{\cong} "a";"b" & \\
                d^{\oneop} \circ d \circ d^{\oneop} \circ d \ar@2{->}[rr]^-{\id \circ (\eta^{\oneop})^{-1} \circ \id} \ar@2{->}@/_1.25pc/[dr]_-{\id \circ \id \circ \eta} \ar@2{->}@/^1.25pc/[ur]^-{\eta \circ \id \circ \id} & \ar@2{->}@{}[d]^(.25){}="a"^(.75){}="b" \ar@3{->}_{\id \circ \tau} "a";"b" & d^{\oneop} \circ d \ar@2{->}[r]^-{\eta} \ar@2{->}@{}[d]^(.25){}="a"^(.75){}="b" \ar@3{->}^{\cong} "a";"b" & \id_{\B}\\
                & d^{\oneop} \circ d \ar@2{->}[ur]_-{\id} \ar@2{->}@/^-1.5pc/[urr]_-{\eta} &&
        }
        \end{matrix}
        = 
        \begin{matrix}
            \xymatrix{
                & d^{\oneop} \circ d  \ar@2{->}[dr]^-{\eta} \ar@2{->}@{}[dd]^(.25){}="a"^(.75){}="b" \ar@3{->}_{\cong} "a";"b"  & \\
                d^{\oneop} \circ d \circ d^{\oneop} \circ d \ar@2{->}[dr]_-{\id \circ \id \circ \eta} \ar@2{->}[ur]^-{\eta \circ \id \circ \id} && \id_{\B}\\
                & d^{\oneop} \circ d \ar@2{->}[ur]_-{\eta} &
        }
        \end{matrix}
    \end{equation}
    A 2-category together with an $\Oone$-volution is called an \emph{$\Oone$-volutive 2-category}.
    \end{definition}
\begin{remark}\label{Remark: Labeling of pasting diagrams}
    When filling pasting diagrams, we allow for some flexibility in the labeling, for instance, let $X,Y \colon a \to b$ be invertible 1-morphisms and $f \colon X \Rightarrow Y$ a 
    2-morphism. Then, we denote the following two 2-morphisms by the same symbol:
    \begin{equation}
        \xymatrix@C=2pc@R=1.5pc{ 
        a 
        \ar@/^1.5pc/[rr]_{\quad}^{X}="1" 
        \ar@/_1.5pc/[rr]_{Y}="2" 
        && b
        \ar@{}"1";"2"|(0.135){\,}="7" 
        \ar@{}"1";"2"|(0.875){\,}="8" 
        \ar@{=>}^{f}"7" ;"8"
        }   
        \hspace{1cm}
        \text{and}
        \hspace{1cm}
        \xymatrix@C=2pc@R=1.5pc{ 
        b 
        \ar@/^1.5pc/[rr]_{\quad}^{Y^{-1}}="1" 
        \ar@/_1.5pc/[rr]_{X^{-1}}="2" 
        && a .
        \ar@{}"1";"2"|(0.135){\,}="7" 
        \ar@{}"1";"2"|(0.875){\,}="8" 
        \ar@{=>}^{f}"7" ;"8"
        }   
    \end{equation}
    Of course, by the right-hand side diagram we mean the 2-morphism 
    \begin{equation}
        \xymatrix@C=3pc{
            Y^{-1} \ar[r]^-{\cong} & Y^{-1} \circ X \circ X^{-1} \ar[r]^-{\id \circ f \circ \id} & Y^{-1} \circ Y \circ X^{-1} \ar[r]^-{\cong} & X^{-1}.
        }
    \end{equation}
\end{remark}

\begin{example} 
	\label{Example: Any Oone volutive category is trivially a Oone volutive 2-category}
    Let $(\C,d,\eta)$ be an $\Oone$-volutive 1-category as in \Cref{def:O1volutiveCategory}. 
    Considering $\C$ as a 2-category with only identity 2-morphisms, $d$ as the induced 2-functor, $\eta$ as the induced 
    2-transformation, and the trivial modification yields an $\Oone$-volutive 2-category $(\C,d,\eta,\id)$.
\end{example}
\begin{example} \label{Example: symmetric monoidal 2-categories with duals are Oone volutive}
    Let $\B$ be a symmetric monoidal 2-category with duals and adjoints. Given any choice of dualization data, consider the dualization 2-functor 
    $(-)^* \colon \B \to \B^{\oneop}$ whose assignment on object and 1-morphism level is described in \Cref{Example: symmetric monoidal categories with duals are Oone volutive}, 
    and which assigns to a 2-morphism $f \colon X \to Y$ the 2-morphism $f^* \colon X^* \to Y^*$ defined via 
    \begin{equation} 
        \xymatrix{
            b^* \ar[rr]^-{\widetilde{\coev}_a \otimes \id} \ar[d]^{\id} && a^* \otimes a \otimes b^* \ar[rr]^{\id \otimes X \otimes \id} \ar[d]^{\id} \ar@{=>}[dll]^-{\id} && a^* \otimes b \otimes b^* \ar[rr]^-{\id \circ \otimes \widetilde{\ev}_b} \ar[d]^{\id} \ar@{=>}[dll]^-{\id \otimes f \otimes \id} && a^* \ar[d]^{\id} \ar@{=>}[dll]^-{\id} \\
            b^* \ar[rr]_-{\widetilde{\coev}_a \otimes \id} && a^* \otimes a \otimes b^* \ar[rr]_{\id \otimes Y \otimes \id} && a^* \otimes b \otimes b^* \ar[rr]_-{\id \circ \otimes \widetilde{\ev}_b} && a^*.
        }
    \end{equation} 
    The rest of the data of the 2-functor $(-)^*$ can also be constructed from the dualization data. As expected from the 1-categorical case discussed in 
    \Cref{Example: symmetric monoidal categories with duals are Oone volutive}, one may complete the functor $(-)^*$ to an $\Oone$-volutive structure on $\B$. In particular, 
    the invertible 2-transformation $((-)^*)^{\oneop} \circ (-)^* \to \id_{\B}$ is constructed from the symmetric braiding.
\end{example}
\begin{definition}
    Let $(\B,d,\eta,\tau)$ and $(\B',d',\eta',\tau')$ be $\Oone$-volutive 2-categories. An \emph{$\Oone$-volutive 2-functor} $(\B,d,\eta,\tau) \to (\B',d',\eta',\tau')$ consists 
    of a 2-functor $F \colon \B \to \B'$, an
    invertible 2-transformation $\alpha \colon F^{\oneop} \circ d \Rightarrow d' \circ F$, and an invertible modification~$\Xi$ fitting into the diagram 
    \begin{equation} 
    	\label{Functor of anti-involutive bicategories diagram}
        \xymatrix{
            (d')^{\oneop} \circ F^{\oneop} \circ d \ar@2{->}[d]_-{\id \circ \alpha} \ar@2{->}[rr]^-{\alpha^{\oneop} \circ \id} && F \circ d^{\oneop} \circ d \ar@2{->}[d]^-{\id \circ \eta} \ar@3{->}[dll]^-{\Xi}  \\
            (d')^{\oneop} \circ d' \circ F \ar@2{->}[rr]_-{\eta' \circ \id} && F
        }
    \end{equation}
    which satisfies a compatiblity condititon with $\tau,\tau'$ that involves the threefold product of $d$, namely, 
    \begin{align*}
        &\begin{matrix}
            \xymatrix{
                F d^{\oneop} d d^{\oneop} \ar@{=>}@/^-2.25pc/[drr]_-{\id \circ \eta \circ \id} \ar@{=>}[drr]^-{\id \circ \id \circ (\eta^{\oneop})^{-1}} \ar@{=>}[rr]^-{(\alpha^{\oneop})^{-1} \circ \id} &  \ar@{=>}@{}[d]^(.55){}="a"^(1.2){}="b" \ar@3{->}^{\tau^{\oneop}} "a";"b" & \ar@{=>}@{}[d]^(.25){}="a"^(.75){}="b" \ar@3{->}^{\cong} "a";"b" d'^{\oneop} F^{\oneop} d d^{\oneop} \ar@{=>}[drr]^-{\id \circ \id \circ (\eta^{\oneop})^{-1}} \ar@{=>}[rr]^-{\id \circ \alpha \circ \id} && \ar@{=>}@{}[d]^(.25){}="a"^(.75){}="b" \ar@3{->}^{\Xi^{\oneop}} "a";"b"  d'^{\oneop} d' F d^{\oneop} \ar@{=>}[rr]^-{\id \circ (\alpha^{\oneop})^{-1}} && d'^{\oneop} d' d'^{\oneop} F^{\oneop} \\
                && F d^{\oneop} \ar@{=>}[rr]_-{(\alpha^{\oneop})^{-1}} && d'^{\oneop} F^{\oneop} \ar@{=>}[urr]_-{\id \circ \eta'^{\oneop} \circ \id} && \\
            }
        \end{matrix}
        \\
        &=
        \begin{matrix}
            \xymatrix{
                F d^{\oneop} d d^{\oneop} \ar@{=>}[drr]_-{\id \circ \eta \circ \id} \ar@{=>}[rr]^-{(\alpha^{\oneop})^{-1} \circ \id} && \ar@{=>}@{}[d]^(.25){}="a"^(.75){}="b" \ar@3{->}_{\Xi^{-1}} "a";"b" d'^{\oneop} F^{\oneop} d d^{\oneop} \ar@{=>}[rr]^-{\id \circ \alpha \circ \id} && \ar@{=>}@{}[d]^(.25){}="a"^(.75){}="b" \ar@3{->}_{\cong} "a";"b" d'^{\oneop} d' F d^{\oneop} \ar@{=>}[rr]^-{\id \circ (\alpha^{\oneop})^{-1}} &\ar@{=>}@{}[d]^(.55){}="a"^(1.2){}="b" \ar@3{->}^{\tau'^{\oneop}} "a";"b" &  d'^{\oneop} d' d'^{\oneop} F^{\oneop} \\
                && F d^{\oneop}  \ar@{=>}[urr]^-{\eta'^{-1} \circ \id \circ \id}  \ar@{=>}[rr]_-{(\alpha^{\oneop})^{-1}} && d'^{\oneop} F^{\oneop} \ar@{=>}[urr]^-{\eta'^{-1} \circ \id \circ \id} \ar@{=>}@/^-2.25pc/[urr]_-{\id \circ \eta'^{\oneop} \circ \id} && \\
            }
        \end{matrix}
    	.
    \end{align*}
\end{definition}

\begin{definition}
    Let $(F,\alpha,\Xi), (\tilde{F},\tilde{\alpha},\tilde{\Xi}) \colon (\B,d,\eta,\tau) \to (\B',d',\eta',\tau')$ be $\Oone$-volutive 2-functors. An 
    \emph{$\Oone$-volutive 2-transformation} $(F,\alpha,\Xi) \to (\tilde{F},\tilde{\alpha},\tilde{\Xi})$ consists of a 2-transformation 
    $\beta \colon F \Rightarrow \tilde{F}$ and an invertible modification $\Omega$ fitting into the diagram 
    \begin{equation*}
        \xymatrix{
            F^{\oneop} \circ d  \ar@2{->}[rr]^-{\alpha}  \ar@3{->}[drr]^-{\Omega} && d' \circ F  \ar@2{->}[d]^-{\id \circ \beta}\\
            \tilde{F}^{\oneop} \circ d  \ar@2{->}[rr]_-{\tilde{\alpha}}  \ar@2{->}[u]^-{\beta^{\oneop} \circ \id} &&  d' \circ \tilde{F}
        }
    \end{equation*}
    which satisfies a compatibility condition with $\Xi, \tilde{\Xi}$ that involves the twofold product of $d$, namely
    \begin{align*}
        &\begin{matrix}
            \xymatrix{
                & F d^{\oneop} d \ar@2{->}@{}[d]^(0.75){}="a"^(1.25){}="b" \ar@3{->}_{\Xi^{-1}} "a";"b"  \ar@2{->}[dr]^-{(\alpha^{\oneop})^{-1} \circ \id}\ar@2{->}[rrrr]^-{\beta \circ \id} && \ar@2{->}@{}[d]^(.15){}="a"^(.65){}="b" \ar@3{->}_{\Omega^{\oneop}} "a";"b" && \tilde{F} d^{\oneop} d \ar@2{->}[dr]^-{\id \circ \eta} \ar@2{->}[dl]_-{(\tilde{\alpha}^{\oneop})^{-1} \circ \id} \ar@2{->}@{}[d]^(0.75){}="a"^(1.25){}="b" \ar@3{->}_{\tilde{\Xi}} "a";"b"   &\\
            F \ar@2{->}[ur]^-{\id \circ \eta^{-1}} \ar@2{->}[dr]_-{\eta'^{-1} \circ \id} &&  \ar@2{->}[dl]_-{\id \circ \alpha} d'^{\oneop} F^{\oneop}d & \ar@2{->}@{}[d]^(.25){}="a"^(.75){}="b" \ar@3{->}_{\Omega} "a";"b"  & d'^{\oneop} \tilde{F}^{\oneop} d \ar@2{->}[ll]_-{\id \circ \beta^{\oneop} \circ \id} \ar@2{->}[dr]^-{\id \circ \tilde{\alpha}} && \tilde{F} \\
                & d'^{\oneop} d' F  \ar@2{->}[rrrr]_-{\id \circ \beta} &&&& d'^{\oneop} d' \tilde{F} \ar@2{->}[ur]_-{\eta'\circ \id} &
            }
        \end{matrix}
        \\
        &=\begin{matrix}
            \xymatrix{
                & F d^{\oneop} d \ar@2{->}[rrrr]^-{\beta \circ \id} && \ar@2{->}@{}[d]^(.25){}="a"^(.75){}="b" \ar@3{->}_{\cong} "a";"b"&& \tilde{F} d^{\oneop} d \ar@2{->}[dr]^-{\id \circ \eta} &\\
            F \ar@2{->}[rrrrrr]^-{\beta} \ar@2{->}[ur]^-{\id \circ \eta^{-1}} \ar@2{->}[dr]_-{\eta'^{-1} \circ \id} &&& \ar@2{->}@{}[d]^(.25){}="a"^(.75){}="b" \ar@3{->}_{\cong} "a";"b" &&& \tilde{F} \\
                & d'^{\oneop} d' F \ar@2{->}[rrrr]_-{\id \circ \beta} &&&& d'^{\oneop} d' \tilde{F} \ar@2{->}[ur]_-{\eta'\circ \id} &
            }
        \end{matrix}
    	.
    \end{align*}
\end{definition}

\begin{definition}
    Let $(\beta,\Omega),(\hat{\beta},\hat{\Omega}) \colon (F,\alpha,\Xi) \to (\tilde{F},\tilde{\alpha},\tilde{\Xi})$ be $\Oone$-volutive 2-transformations. An 
    \emph{$\Oone$-volutive modification} $(\beta,\Omega) \to (\hat{\beta},\hat{\Omega})$ is a modification 
    $\Lambda \colon \beta \Rrightarrow \hat{\beta}$ satisfying a compatibility condition with $\Omega,\hat{\Omega}$, namely, 
    \begin{equation}
        \begin{matrix}
            \xymatrix{
                F_1^{\oneop} \circ d  \ar@2{->}[rr]^-{\alpha_1}  \ar@3{->}[ddrr]^-{\Omega} && d' \circ F_1  \ar@2{->}[dd]_-{\id \circ \beta} \ar@2{->}@/_-3pc/[dd]^-{\id \circ \hat{\beta}}\\
                &\ar@2{->}@{}[r]^(1){}="a"^(1.8){}="b" \ar@3{->}_{\id \circ \Lambda} "a";"b" &  \\
                F_2^{\oneop} \circ d  \ar@2{->}[rr]^-{\alpha_2}  \ar@2{->}[uu]^-{\beta^{\oneop} \circ \id} &&  d' \circ F_2
        }
        \end{matrix}
        \;=\; 
        \begin{matrix}
            \xymatrix{
                F_1^{\oneop} \circ d  \ar@2{->}[rr]^-{\alpha_1}  \ar@3{->}[ddrr]^-{\hat{\Omega}} && d' \circ F_1  \ar@2{->}[dd]^-{\id \circ \hat{\beta}}\\
                &\ar@2{->}@{}[l]^(1){}="a"^(1.8){}="b" \ar@3{<-}_{\Lambda^{\oneop} \circ \id} "a";"b" &  \\
                F_2^{\oneop} \circ d  \ar@2{->}[rr]^-{\alpha_2}  \ar@2{->}[uu]_-{\hat{\beta}^{\oneop} \circ \id} \ar@2{->}@/_-3pc/[uu]^-{\beta^{\oneop} \circ \id} &&  d' \circ F_2
        }
        \end{matrix}
    	.
    \end{equation}
\end{definition}

Let $(F,\alpha,\Xi) \colon (\B_1,d_1,\eta_1,\tau_1) \to (\B_2,d_2,\eta_2,\tau_2)$ and $(F',\alpha',\Xi') \colon 
(\B_2,d_2,\eta_2,\tau_2) \to (\B_3,d_3,\eta_3,\tau_3)$ be $\Oone$-volutive 2-functors. The \emph{composition} 
$(F',\alpha',\Xi') \circ (F,\alpha,\Xi) := (\underline{F},\underline{\alpha},\underline{\Xi})$ consists of the 2-functor $\underline{F} := F' \circ F$, the 2-transformation  
\begin{equation}
	\label{eq:ComposedAlpha}
	\underline{\alpha} := 
	\Big( 
    \xymatrix{
        {F'}^{\oneop} \circ F^{\oneop} \circ d_1 \ar@2{->}[r]^-{\id \circ \alpha} & {F'}^{\oneop} \circ d_2 \circ F \ar@2{->}[r]^-{\alpha' \circ \id} & d_3 \circ F' \circ F 
    }
	\Big) 
\end{equation}
and the modification $\underline{\Xi}$ obtained by filling the diagram 
\begin{equation}
	\label{eq:ComposedXi}
    \begin{aligned}
        \xymatrix{
            F'F d_1^{\oneop} d_1 \ar@2{->}[ddd]_-{\id \circ \eta_1} & && & d_3^{\oneop} {F'}^{\oneop}F^{\oneop} d_1 \ar@2{->}[ddd]^-{\id \circ \underline{\alpha}} \ar@2{->}[llll]_-{\underline{\alpha}^{\oneop} \circ \id} \ar@2{->}[ddl]^-{\id \circ \id \circ \alpha} \ar@2{->}[dlll]^-{{\alpha'}^{\oneop} \circ \id \circ \id}\\
            & F'd_2^{\oneop} F^{\oneop} d_1 \ar@2{->}[d]_{\id \circ  \id \circ \alpha} \ar@2{->}[ul]_-{\id \circ \alpha^{\oneop} \circ \id}&&& \\ 
            &F' d_2^{\oneop} d_2 F \ar@2{->}[dl]_-{\id \circ \eta_2 \circ \id} && d_3^{\oneop} {F'}^{\oneop} d_2 F \ar@2{->}[ll]_-{{\alpha'}^{\oneop} \circ \id \circ \id} \ar@2{->}[dr]^-{\id \circ \alpha' \circ \id}& \\
            F'F&&&& d_3^{\oneop}d_3 F' F \ar@2{->}[llll]^-{\eta_3 \circ \id}
        }
    \end{aligned}
	.
\end{equation}
Here, the left square is filled by $\Xi$, the lower square is filled by $\Xi'$, the upper and right triangles commute by definition, and the middle square commutes by the 
interchange law. 
One checks that the coherence properties of $\Xi$ and $\Xi'$ imply the correct coherence property for $\underline{\Xi}$ in the sense that we have an  
$\Oone$-volutive 2-functor $(\underline{F},\underline{\alpha},\underline{\Xi}) \colon  (\B_1,d_1,\eta_1,\tau_1) \to (\B_3,d_3,\eta_3,\tau_3)$. The composition 
laws for $\Oone$-volutive 2-transformations and $\Oone$-volutive modifications, respectively, are simpler and we omit the details. 

\begin{remark} 
	\label{Remark: Expectation that Oone volutive 2-categories form a 3-category}
    We expect that $\Oone$-volutive 2-categories, $\Oone$-volutive 2-functors, $\Oone$-volutive 2-transformations, and $\Oone$-volutive modifications form a 
    3-category $\Oonevolbicat$. To prove this, one would have to check the necessary coherences a 3-category satisfies, which is rather tedious. 
    Extending \Cref{Example: Any Oone volutive category is trivially a Oone volutive 2-category}, we then expect there to be an 
    inclusion 3-functor 
    \begin{equation}
        \underline{\Oonevolcat} \xymatrix{\ar[r]&} \Oonevolbicat 
    \end{equation}
    where $\underline{\Oonevolcat}$ is $\Oonevolcat$ thought of as a 3-category with only identity 3-morphisms, 
    which assigns to each $\Oone$-volutive category the $\Oone$-volutive 2-category described in \Cref{Example: Any Oone volutive category is trivially a Oone volutive 2-category}.
\end{remark}

\begin{remark}
    Under the assumption that there is a 3-category $\Oonevolbicat$, we expect there to be a 3-functor 
    \begin{equation}
        \symmontwocat^{\mathrm{d}} \xymatrix{\ar[r]&} \Oonevolbicat
    \end{equation}
    assigning to each symmetric monoidal 2-category with duals and adjoints the $\Oone$-volutive 2-category with the same underlying 2-category and the 
    $\Oone$-volutive structure indicated in \Cref{Example: symmetric monoidal 2-categories with duals are Oone volutive}. 
\end{remark}

\subsubsection{Dagger case}

\begin{definition}
    Let $\B$ be a 2-category. 
    An \emph{$\Oone$-dagger structure} on $\B$ is a pair $(d,\eta)$ consisting of an identity-on-objects 2-functor $d \colon \B \to \B^{\oneop}$ and
    an invertible 2-transformation $\eta \colon d^{\oneop} \circ d \Rightarrow \id_{\B}$ whose 1-morphism components are identities, and $(\eta^{\oneop})^{-1} \circ \id_d = \id_{d} \circ \eta$ holds. 
    A 2-category together with an $\Oone$-dagger structure is called an \emph{$\Oone$-dagger 2-category}.
\end{definition}

In the context of Q-system completions for $\textrm{C}^*$-categories, $\Oone$-dagger 2-categories are called ``$\dagger$ 2-categories'', see e.g.\ \cite{chen2021qsystemcompletion3functor}. 
Any $\Oone$-dagger 2-category $(\B,d,\eta)$ can be trivially reinterpreted as an $\Oone$-volutive 2-category $(\B,d,\eta,\id)$. 
We come back to this after introducing 
(higher) morphisms between $\Oone$-dagger 2-categories. 

\begin{definition}
    Let $(\B,d,\eta)$ and $(\B',d',\eta')$ be $\Oone$-dagger 2-categories. An \emph{$\Oone$-dagger 2-functor} $(\B,d,\eta) \to (\B',d',\eta')$ consists 
    of a 2-functor $F \colon \B \to \B'$ and an invertible 2-transformation $\alpha \colon F^{\oneop} \circ d \Rightarrow d' \circ F$ whose 1-morphism components are identities, and $(\id \circ \eta) \bullet (\alpha^{\oneop} \circ \id) = (\eta' \circ \id) \bullet (\id \circ \alpha)$ holds.
\end{definition}

\begin{definition}
    Let $(F,\alpha), (\tilde{F},\tilde{\alpha}) \colon (\B,d,\eta) \to (\B',d',\eta')$ be $\Oone$-dagger 2-functors. A 
    \emph{$\Oone$-dagger 2-transformation} $(F,\alpha) \to (\tilde{F},\tilde{\alpha})$ is a 2-transformation 
    $\beta \colon F \Rightarrow \tilde{F}$ such that $(\id \circ \beta) \bullet \alpha \bullet (\beta^{\oneop} \circ \id) = \tilde{\alpha}$ holds.
\end{definition}

\begin{definition}
    Let $\beta,\hat{\beta} \colon (F,\alpha) \to (\tilde{F},\tilde{\alpha})$ be $\Oone$-dagger 2-transformations. A 
    \emph{$\Oone$-dagger modification} $\beta \to \hat{\beta}$ is a modification 
    $\Lambda \colon \beta \Rrightarrow \hat{\beta}$ satisfying $\id \circ \Lambda = \Lambda^{\oneop} \circ \id$.
\end{definition}

\begin{remark}
    We expect $\Oone$-dagger 2-categories, $\Oone$-dagger 2-functors, $\Oone$-dagger 2-transformations, and $\Oone$-dagger modifications 
    to form a 3-category $\strictOonevolbicat$. To prove this, one needs to describe the various compositions and coherences in this 3-category, which is expected to be 
    slightly easier than the corresponding discussion in the volutive case.
\end{remark}

It is clear that any $\Oone$-dagger 2-category $(\B,d,\eta)$ can be \emph{trivially} reinterpreted as an $\Oone$-volutive 2-category $(\B,d,\eta,\id)$, and 
we expect that this assignment is functorial in the sense that there is a 3-functor 
\begin{equation}
    T_{\Oone} \colon \strictOonevolbicat \xymatrix{\ar[r]&} \Oonevolbicat.
\end{equation}

\subsubsection{Strictification}

Conversely, any $\Oone$-volutive 2-category $(\B,d,\eta,\tau)$ gives rise to an $\Oone$-dagger 2-category, which we will describe in the following.  
Let $S_{\Oone}\B$ be the 2-category whose objects are tuples $(c,\theta_c,\Pi_c)$ where $c \in \B$ is an object, $\theta_c \colon c \to d(c)$ is a 1-isomorphism, 
and $\Pi_c \colon \eta_c \circ d(\theta_c)^{-1} \circ \theta_c \Rightarrow \id_c$ is a 2-isomorphism satisfying 
\begin{equation} \label{Equation: Coherence condition for strictification of Oone volutive 2-categories}
    \begin{matrix}
        \xymatrix{
            \ar@{}[ddrr]^(.375){}="a"^(.625){}="b" \ar@{=>}_{\Pi_c} "a";"b"  & d^2(c) \ar[dr]^-{\eta_c} \ar[rr]^-{d^2(\theta_c)} && d^3(c) \ar[dr]^-{\eta_{d(c)}} \ar@{}[dl]^(.25){}="a"^(.75){}="b" \ar@{=>}_{\cong} "a";"b" \\
            d(c) \ar[ur]^-{d(\theta_c)^{-1}} && c \ar[rr]_-{\theta_c} & \ar@{}[dl]^(.25){}="a"^(.75){}="b" \ar@{=>}_{\cong} "a";"b"& d(c) \ar[dl]^-{\id_{d(c)}} \\
            & \ar[ul]^-{\theta_c} \ar[ur]_-{\id_c} c \ar[rr]_-{\theta_c} && d(c) 
        }
    \end{matrix}
    = 
    \begin{matrix}
        \xymatrix{
            & d^2(c)  \ar[rr]^-{d^2(\theta_c)} &  \ar@{}[dl]^(.25){}="a"^(.75){}="b" \ar@{=>}_{d(\Pi_c)^{-1}} "a";"b"& d^3(c) \ar[dr]^-{\eta_{d(c)}} \ar[dl]^-{d(\eta_c)^{-1}} \\
            d(c) \ar[rr]_-{\id_{d(c)}} \ar[ur]^-{d(\theta_c)^{-1}} &  \ar@{}[dr]^(.25){}="a"^(.75){}="b" \ar@{=>}_{\cong} "a";"b"& d(c) \ar[dr]^-{\id_{d(c)}} && d(c) \ar[dl]^-{\id_{d(c)}}  \ar@{}[ll]^(.375){}="a"^(.625){}="b" \ar@{=>}_{\tau_c} "a";"b"\\
            & \ar[ul]^-{\theta_c} c \ar[rr]_-{\theta_c} && d(c) 
        }
    \end{matrix}
.
\end{equation}
A 1-morphism $(c,\theta_c,\Pi_c) \to (c',\theta_{c'}, \Pi_{c'})$ in $S_{\Oone}\B$ is a 1-morphism $X \colon c \to c'$ in $\B$, and a 2-morphism $X \to Y$ $S_{\Oone}\B$ 
is a 2-morphism $f \colon X \to Y$ in $\B$. It is clear that $S_{\Oone}\B$ defined in this way becomes a 2-category. 

Next, we define an $\Oone$-dagger structure on $S_{\Oone}\B$. 
The identity-on-objects 2-functor $S_{\Oone} d$ assigns to a 1-morphism 
$X \colon (a,\theta_a,\Pi_a) \to (b,\theta_b,\Pi_b)$ the 1-morphism $ (b,\theta_b,\Pi_b) \to (a,\theta_a,\Pi_a)$ given by 
\begin{equation}
    \xymatrix{
        b \ar[r]^-{\theta_b} & d(b) \ar[r]^-{d(X)} & d(a) \ar[r]^-{\theta_a^{-1}} & a
    }
\end{equation}
and to a 2-morphism $f \colon X \to Y$ the 2-morphism $S_{\Oone} d(X) \to S_{\Oone} d(Y)$ defined by
\begin{equation}
    \xymatrix{
        b \ar[rr]^-{\theta_b} \ar[d]_-{\id} && d(b) \ar[rr]^-{d(X)} \ar[d]^-{\id} \ar@{=>}[dll]^-{\id_{\theta_b}} && d(a) \ar[rr]^-{\theta_a^{-1}} \ar[d]^-{\id} \ar@{=>}[dll]^-{d(f)} && a \ar[d]^-{\id} \ar@{=>}[dll]^-{\id_{\theta_a^{-1}}}\\
        b \ar[rr]_-{\theta_b} && d(b) \ar[rr]_-{d(Y)} && d(a) \ar[rr]_-{\theta_a^{-1}} && a.
    }
\end{equation}
One routinely checks that this defines a 2-functor by using the respective properties of $d$. 

Next, we define the 2-transformation $S_{\Oone} \eta \colon S_{\Oone} d^{\oneop} \circ S_{\Oone} d \Rightarrow \id_{S_{\Oone}\B}$. 
Its 1-morphism components are identities. 
The 2-morphism component of this 2-transformation at a 1-morphism $X \colon (a,\theta_a,\Pi_a) \to (b,\theta_b,\Pi_b)$ is given by
\begin{equation}
    \xymatrix{
        X \cong \id_b  X  \id_a \ar[rr]^-{\Pi_b \circ \id_X \circ \Pi_a^{-1}} && \theta_b^{-1}  d(\theta_b) \eta_b^{-1}  X  \eta_a  d(\theta_a)^{-1}  \theta_a
        \ar[rr]^-{\id \circ \eta_X^{-1} \circ \id} && \theta_b^{-1}  d(\theta_b) d^2(X)  d(\theta_a)^{-1}  \theta_a .
    }
\end{equation} 
One checks that this defines a 2-transformation by using the respective properties of $\eta$. 
It is also clear by construction that $\eta$ is an invertible 2-transformation. 
One deduces from~\eqref{Equation: Coherence condition for strictification of Oone volutive 2-categories} that $(S_{\Oone}\eta^{\oneop})^{-1} \circ \id_{S_{\Oone} d} = 
\id_{S_{\Oone} d} \circ S_{\Oone}\eta$ holds. 

\begin{lemma}
    $(S_{\Oone}\B,S_{\Oone}d,S_{\Oone}\eta)$ is an $\Oone$-dagger 2-category.
\end{lemma}

We expect that the assignment $(\B,d,\eta,\tau) \mapsto (S_{\Oone}\B,S_{\Oone}d,S_{\Oone}\eta)$ extends to a 3-functor 
\begin{equation}
    S_{\Oone} \colon \Oonevolbicat \xymatrix{\ar[r]&} \strictOonevolbicat.
\end{equation}
In the following, we will describe the action of the conjectural 3-functor $ S_{\Oone}$ on 1-morphism level. Let  $(\B,d,\eta,\tau)$ and $(\B',d',\eta',\tau')$ be 
$\Oone$-volutive 2-categories and let $(F,\alpha,\Xi) \colon (\B,d,\eta,\tau) \to (\B',d',\eta',\tau')$ be an $\Oone$-volutive 2-functor. We construct an 
$\Oone$-dagger 2-functor $(S_{\Oone}F, S_{\Oone}\alpha) \colon S_{\Oone}(\B,d,\eta,\tau) \to S_{\Oone}(\B',d',\eta',\tau')$ as follows. The 2-functor $S_{\Oone}\B \to S_{\Oone}\B'$ 
assigns an object $(c,\theta_c,\Pi_c)$ to the tuple consisting of the object $F(c)$, the 1-isomorphism $\tilde{\theta}_{F(c)}$ defined as
\begin{equation}
	\label{eq:thetaTilde}
	\tilde{\theta}_{F(c)} = 
	\Big(
    \xymatrix{
        F(c) \ar[r]^-{F(\theta_c)} & F(d(c)) \ar[r]^-{\alpha_c} & d'(F(c))
    }
	\Big)
\end{equation}
and the 2-isomorphism $\tilde{\Pi}_{F(c)}$ obtained by
\begin{equation}
    \xymatrix{
        && d'(F(c)) \ar[rr]^-{d'(\tilde{\theta}_{F(c)})^{-1}} \ar[dr]^-{d'(F(\theta_c))^{-1}} \ar@{}[d]^(.25){}="a"^(1.25){}="b" \ar@{=>}^{\alpha_{\theta_c^{-1}}} "a";"b" && {d'}^2(F(c)) \ar[drrr]^-{\eta'_{F(c)}} \ar[dd]^-{\alpha_{d(c)}^{-1} \bullet d'(\alpha_c)} &\ar@{}[d]^(.55){}="a"^(1.45){}="b" \ar@{=>}^{\Xi_c} "a";"b" \\
        F(c) \ar@/^1.5pc/[urr]^-{\tilde{\theta}_{F(c)}} \ar[r]^-{F(\theta_c)} \ar@/^-7pc/[rrrrrrr]_-{\id_{F(c)}} & F(d(c)) \ar[ur]^-{\alpha_c} \ar[drrr]_-{F(d(\theta_c))^{-1}} && d'(F(d(c))) \ar[ur]^-{d'(\alpha_c)^{-1}} \ar[dr]^-{\alpha_{d(c)}^{-1}} \ar@{}[d]^(1){}="a"^(1.75){}="b" \ar@{=>}_{F(\Pi_c)} "a";"b" &&&& F(c)\\
        &&&&F(d^2(c)) \ar[urrr]_-{F(\eta_c)}&&&
    }
\end{equation}
where the unfilled diagrams commute by definition. One checks that the triple $(F(c), \tilde{\theta}_{F(c)}, \tilde{\Pi}_{F(c)})$ defines an object in 
$S_{\Oone}(\B',d',\eta',\tau')$ by using the coherence condition of $\Pi$. On 1- and 2-morphism level, the assignment of the 2-functor $S_{\Oone}\B \to S_{\Oone}\B'$  
is trivial. 
Next, we describe the 2-transformation $S_{\Oone}\alpha \colon S_{\Oone}F^{\oneop} \circ S_{\Oone}d \to S_{\Oone}d' \circ S_{\Oone}F$. 
Its 1-morphism components are identities. 
Its 2-morphism component at $X \colon (a,\theta_a,\Pi_a) \to (b,\theta_b,\Pi_b)$ is given by the 2-morphism 
\begin{equation}
    \xymatrix{
        \tilde{\theta}_{F(a)}^{-1} d'(F(X)) \tilde{\theta}_{F(b)} \ar[rr]^-{\id \circ \alpha_X \circ \id} && \tilde{\theta}_{F(a)}^{-1} \alpha_a F(d(X)) \alpha_b^{-1} \tilde{\theta}_{F(b)} = F(\theta_a^{-1}) F(d(X)) F(\theta_b) \cong F(\theta_a^{-1}d(X)\theta_b)
    }
\end{equation}
where we again used the definition of $\tilde{\theta}$ in~\eqref{eq:thetaTilde}. 
One checks that this defines an invertible 2-transformation and moreover that the pair $(S_{\Oone}F,S_{\Oone}\alpha)$ 
defines an $\Oone$-dagger 2-functor. Similarly, one describes the action of $S_{\Oone}$ on 2- and 3-morphisms. 

\begin{remark}
    We expect that the 3-functor $S_{\Oone}$ is right adjoint to $T_{\Oone}$ in the 4-category of 3-categories.
\end{remark}

\begin{remark} \label{Remark on coherently self-dual objects}
    Let $\B$ be a symmetric monoidal 2-category with duals and adjoints. Applying the construction described in \Cref{Example: symmetric monoidal 2-categories with duals are Oone volutive}
    yields an $\Oone$-volutive 2-category $\B$, to which we may apply in turn $S_{\Oone}$, resulting in an $\Oone$-dagger 2-category $\B$. The objects in the underlying 
    2-category are by definition objects $a \in \B$ together with 1-isomorphisms $\theta_a \colon a \to a^*$ and 2-isomorphisms 
    $\Pi_a \colon \eta_a \circ (\theta_a^{-1})^* \circ \theta_a \Rightarrow \id_a$ satisfying a coherence property involving the triple dual, cf.~\eqref{Equation: Coherence condition for strictification of Oone volutive 2-categories}.
    We say that an object $a \in \B$ is \emph{coherently self-dual} if it admits such data.
\end{remark}

\subsection[$\Otwo$-dagger 2-categories]{$\textrm{\textbf{O}}\pmb{(2)}$-dagger 2-categories}
\label{Subsection: Otwovolutive 2-categories}

This section is analogous to \Cref{Subsection: Oonevolutive 2-categories}, with $\Oone$ replaced by $\Otwo$.
We propose the basic structure of 3-categories of $\Otwo$-volutive 2-categories, $\Otwo$-dagger 2-categories, and we discuss strictification. 

\subsubsection{Volutive case}

Analogously to the motivation given at the beginning of \Cref{Subsection: rspintwovolutive 2-categories}, there should be an $\Otwo$-action 
on the 3-category $\BiCatadj$ whose homotopy fixed points are $\Otwo$-volutive 2-categories. To construct this action, it seems helpful to realize 
$\Otwo$ as the semidirect product 2-group $\Oone \rtimes \Sotwo$, cf.\ \cite[Appendix B.3]{müller2023reflection}. The problem of describing 
$\Otwo$-actions then decomposes into a study of $\Sotwo$-actions, $\Oone$-actions, and their compatibility. 
Carrying out the details 
of this construction and computing the homotopy fixed points is again difficult, and we propose instead the following definitions, partially following 
\cite{StehouwerMueller2024}. As indicated, an $\Otwo$-volutive structure should consist of an $\Oone$-volutive structure, an $\Sotwo$-volutive structure, and 
a compatibility between them. 

Recall that if a 2-category $\B$ has right adjoints, then $\B^{\oneop}$ and $\B^{\twoop}$ have left adjoints. Moreover, recalling that the adjunction 2-functor 
$(-)^\vee \colon \B^{\oneop} \to \B^{\twoop}$ is an equivalence, one may translate between different conventions on which oppositization is used.
The following definition was originally proposed in \cite{StehouwerMueller2024}.
\begin{definition}
    Let $\B$ be a 2-category with right adjoints. An \emph{$\Otwo$-volutive structure} on $\B$ consists of an $\Oone$-volutive structure $(d,\eta,\tau)$ on $\B$, an 
    $\Sotwo$-volutive structure $S$ on $\B$, and an invertible modification $\Gamma$ fitting into the diagram 
        \begin{equation}
            \xymatrix{
                ^{\vee \vee}(-) \circ d \ar@2{->}[rr]^-{\cong} \ar@2{->}[d]_-{S \circ \Id_d} && d \circ (-)^{\vee \vee} \ar@3{->}[dll]^-{\Gamma} \ar@2{->}[d]^-{\Id_d \circ S^{-1}} \\
                \id_{\B^{\oneop}} \circ d \ar@2{->}[rr]_-{\cong} && d \circ \id_{\B} 
            }
        \end{equation}
    which satisfies a compatibility condition with $\eta$, namely, the equality of modifications
    \begin{equation}\label{Equation: The compatibility condition of Gamma in an Otwo volution}
        \begin{matrix}
        \xymatrix@C=1em{
            & (-)^{\vee \vee} \ar@2{->}[ddr]^-{\eta^{-1} \circ \id} \ar@2{->}[ddl]_-{\id \circ \eta^{-1}} \ar@2{->}[d]^-{S^{-1}} & \\
            & \id \ar@2{->}[d]^-{\eta^{-1}} &\\
            (-)^{\vee \vee} \circ d^{\oneop} \circ d \ar@2{->}[dr]_-{\cong} & d^{\oneop} \circ d \ar@2{->}[d]^-{\id \circ S^{-1} \circ \id} \ar@2{->}[r]^-{\id \circ S} \ar@2{->}[l]_-{S \circ \id} & d^{\oneop} \circ d \circ (-)^{\vee \vee} \ar@2{->}[dl]^-{\cong}\\
            & d^{\oneop} \circ ^{\vee \vee}(-) \circ d &
        }
        \end{matrix}
        = 
        \begin{matrix}
            \xymatrix@C=0.5em{
                & (-)^{\vee \vee} \ar@2{->}[ddr]^-{\eta^{-1} \circ \id} \ar@2{->}[ddl]_-{\id \circ \eta^{-1}} & \\
                &&\\
                (-)^{\vee \vee} \circ d^{\oneop} \circ d \ar@2{->}[dr]_-{\cong} & & d^{\oneop} \circ d \circ (-)^{\vee \vee}\ar@2{->}[dl]^-{\cong}  \\
                & d^{\oneop} \circ ^{\vee \vee}(-) \circ d &
            }
            \end{matrix}
        .
    \end{equation}
    Here on the left-hand side the upper two diagrams are filled by the interchange law and the lower two diagrams are filled employing $\Gamma$, whereas on 
    the right-hand side we use the naturality of the 2-transformation exchanging a 2-functor with the double adjoint 2-functor with respect to $\eta$.
\end{definition}

\begin{example}\label{Example: symmetric monoidal 2-categories with duals are Otwo volutive}
    Let $\B$ be a symmetric monoidal 2-category with duals and adjoints. 
    As explained in Examples~\ref{Example: symmetric monoidal 2-categories with duals are Oone volutive} and~\ref{Example: Symmetric monoidal 2-categories with duals and adjoints are rSpintwo volutive}, there is an $\Oone$-volutive structure on $\B$ given by the dualization 
    2-functor and an $\Sotwo$-structure on $\B$ given by the Serre automorphism. These structures in fact combine to an $\Otwo$-volutive structure. The 2-morphism 
    components of the invertible modification $\Gamma$ are described in \cite{StehouwerMueller2024}, see also \cite[Appendix B.3]{müller2023reflection}. 
    Restricting to the maximal subgroupoid $\B^{\times}$, this $\Otwo$-volutive structure is expected to reduce to the canonical $\Otwo$-action on $\B^{\times}$ 
    induced by the cobordism hypothesis.
\end{example}

The following is a subexample of \Cref{Example: symmetric monoidal 2-categories with duals are Otwo volutive}.

\begin{example} \label{Example: The delooping of vector spaces is strict Otwo volutive}
    We consider again the symmetric monoidal 2-category $\mathrm{B}\vect$, which has duals and adjoints. By 
    \Cref{Example: symmetric monoidal 2-categories with duals are Otwo volutive}, $\mathrm{B}\vect$ carries an $\Otwo$-volutive structure which we now spell out. 
    Since we already described the underlying $\Sotwo$-dagger structure in \Cref{Example: The delooping of vector spaces is strict rSpintwo volutive} (for $r=1$), 
    we now focus on the underlying $\Oone$-volutive structure. 
    Since the single object of $\mathrm{B}\vect$ is self-dual with trivial evaluation and coevaluation 
    morphism, the 2-functor $d \colon \mathrm{B}\vect \to \mathrm{B}\vect^{\oneop}$ is trivial. The higher coherence data $\eta$ and $\tau$ belonging to the 
    $\Oone$-volutive structure are trivial as well. 
    Finally, we compute the modification $\Gamma$ to be trivial. We note that the $\Otwo$-volutive structure described 
    here is strict in the sense of \Cref{Definition: strict Otwo volutive 2-category} below.
\end{example}

Higher maps between between $\Otwo$-volutive 2-categories have not been systematically constructed. 
The following definitions seem to be natural, but we take them to be provisional: 

\begin{definition}
    Let $(\B,d,\eta,\tau,S,\Gamma)$ and $(\B',d',\eta',\tau',S',\Gamma')$ be $\Otwo$-volutive 2-categories. An \emph{$\Otwo$-volutive 2-functor} 
    $(\B,d,\eta,\tau,S,\Gamma) \to (\B',d',\eta',\tau',S',\Gamma')$ consists of a 2-functor $F \colon \B \to \B'$ together with an $\Oone$-volutive structure $(\alpha,\Xi)$
    and an $\Sotwo$-volutive structure $\zeta$.
\end{definition}

\begin{definition} 
    Let $(F,\alpha,\Xi,\zeta),(\tilde{F},\tilde{\alpha},\tilde{\Xi},\tilde{\zeta}) \colon (\B,d,\eta,\tau,S,\Gamma) \to (\B',d',\eta',\tau',S',\Gamma')$ be
    $\Otwo$-volutive 2-functors. An \emph{$\Otwo$-volutive 2-transformation} $(F,\alpha,\Xi,\zeta) \to (\tilde{F},\tilde{\alpha},\tilde{\Xi},\tilde{\zeta})$
    consists of an $\Oone$-volutive 2-transformation $(\beta,\Omega) \colon (F,\alpha,\Xi) \to (\tilde{F},\tilde{\alpha},\tilde{\Xi})$ such that 
    $\beta \colon F \to \tilde{F}$ is also an $\Sotwo$-volutive 2-transformation. 
\end{definition}

\begin{definition}
    Let $(\beta,\Omega),(\hat{\beta},\hat{\Omega}) \colon (F,\alpha,\Xi,\zeta) \to (\tilde{F},\tilde{\alpha},\tilde{\Xi},\tilde{\zeta})$ be $\Otwo$-volutive 2-transformations. 
    An \emph{$\Otwo$-volutive modification} $(\beta,\Omega) \to (\hat{\beta},\hat{\Omega})$ is an $\Oone$-volutive modification.
\end{definition}

Recalling our definition of the composition of $\Oone$-volutive 2-functors (see~\eqref{eq:ComposedAlpha} and~\eqref{eq:ComposedXi}) as well as $\Sotwo$-volutive 2-functors (see~\eqref{Equation: rspintwo chapter the composed invertible modification}), we immediately know how to 
compose $\Otwo$-volutive 2-functors. 
\begin{remark} 
	\label{Remark: The 3-functor from symmetric monoidal 2-categories with duals to Otwo volutive 2-categories}
    We expect that $\Otwo$-volutive 2-categories, $\Otwo$-volutive 2-functors, $\Otwo$-volutive 2-transformations, and $\Otwo$-volutive modification form a 3-category 
    $\Otwovolbicat$. 
    Moreover, we expect to find a 3-functor 
    \begin{equation}
        \symmontwocat^{\mathrm{d}} \xymatrix{\ar[r]&} \Otwovolbicat
    \end{equation}
    assigning to each symmetric monoidal 2-category with duals and adjoints the $\Otwo$-volutive 2-category with the same underlying 2-category and the 
    $\Otwo$-volutive structure indicated in \Cref{Example: symmetric monoidal 2-categories with duals are Otwo volutive}. 
\end{remark}

\subsubsection{Dagger case}

\begin{definition} \label{Definition: strict Otwo volutive 2-category}
    Let $\B$ be a 2-category with right adjoints. 
    An \emph{$\Otwo$-dagger structure} on $\B$ consists of an $\Oone$-dagger structure $(d,\eta)$ and an 
    $\Sotwo$-dagger structure $S$ on $\B$ such that the diagram
    \begin{equation}
        \xymatrix{
            ^{\vee \vee}(-) \circ d \ar@2{->}[rr]^-{\cong} \ar@2{->}[d]_-{S \circ \Id_d} && d \circ (-)^{\vee \vee} \ar@2{->}[d]^-{\Id_d \circ S^{-1}} \\
            \id_{\B^{\oneop}} \circ d \ar@2{->}[rr]_-{\cong} && d \circ \id_{\B} 
        }
    \end{equation}
    commutes. 
\end{definition}

We remark that here the 2-transformation $^{\vee \vee}(-)  \circ d \cong d \circ (-)^{\vee \vee}$ can always be chosen to have identity 1-morphism components, 
which we will assume. We denote the 2-morphism components generically by $\epsilon_X \colon d(X^{\vee \vee}) \to {}^{\vee \vee}d(X)$. 
Any $\Otwo$-dagger 2-category $(\B,d,\eta,S)$ can be trivially reinterpreted as an $\Otwo$-volutive 2-category $(\B,d,\eta,\id,S,\id)$; we come back to this point 
after defining appropriate notions of (higher) morphism between $\Otwo$-dagger 2-categories.

\begin{definition}
    Let $(\B,d,\eta,S)$ and $(\B',d',\eta',S')$ be $\Otwo$-dagger 2-categories. An \emph{$\Otwo$-dagger 2-functor} 
    $(\B,d,\eta,S) \to (\B',d',\eta',S')$ consists of a 2-functor $F \colon \B \to \B'$ together with an $\Oone$-dagger structure $\alpha$ such that 
    $F$ is $\Sotwo$-dagger.
\end{definition}
\begin{definition}
    Let $(F,\alpha),(\tilde{F},\tilde{\alpha}) \colon (\B,d,\eta,S) \to (\B',d',\eta',S')$ be
    $\Otwo$-dagger 2-functors. An \emph{$\Otwo$-dagger 2-transformation} $(F,\alpha) \to (\tilde{F},\tilde{\alpha})$ is a 
    $\Oone$-dagger 2-transformation $\beta \colon (F,\alpha) \to (\tilde{F},\tilde{\alpha})$ such that $\beta$ is 
    $\Sotwo$-dagger. 
\end{definition}
\begin{definition}
    Let $\beta,\hat{\beta} \colon (F,\alpha),(\tilde{F},\tilde{\alpha})$ be $\Otwo$-dagger 2-transformations. A 
    \emph{$\Otwo$-dagger modification} $\beta \to \hat{\beta}$ is an $\Oone$-dagger modification.
\end{definition}
It is expected that $\Otwo$-dagger 2-categories, $\Otwo$-dagger 2-functors, $\Otwo$-dagger 2-transformations, and $\Otwo$-dagger modification
form a 3-category, which we denote by $\strictOtwovolbicat$. 

It is clear that any $\Otwo$-dagger 2-category $(\B,d,\eta,S)$ can be \emph{trivially} reinterpreted as an $\Otwo$-volutive 2-category $(\B,d,\eta,\id,S,\id)$, and 
we expect that this assignment is functorial in the sense that it extends to a 3-functor 
\begin{equation}
	\label{eq:TO2}
    T_{\Otwo} \colon \strictOtwovolbicat \xymatrix{\ar[r]&} \Otwovolbicat.
\end{equation}

\subsubsection{Strictification}

Conversely, any $\Otwo$-volutive 2-category $(\B,d,\eta,\tau,S,\Gamma)$ gives rise to an $\Otwo$-dagger 2-category, which we describe in the following. 
Let $\epsilon$ denote the 2-transformation $^{\vee \vee}(-) \circ d \cong d \circ (-)^{\vee \vee}$, which in general may have non-trivial 1-morphism components.
Let $S_{\Otwo}\B$ be the 2-category whose objects are tuples $(c,\theta_c,\Pi_c,\lambda_c,\phi_c)$ where $(c,\theta_c,\Pi_c)$ is an object in $S_{\Oone}\B$, 
$(c,\lambda_c)$ is an object in $S_{\Sotwo}\B$, and 
\begin{equation}
	\phi_c \colon \theta_c^{-1} \circ \epsilon_c \circ ^{\vee \vee}\theta_c \xymatrix{\ar[r]^{\cong}&} \id_c
\end{equation}
is a 2-isomorphism in $\B$ such that
\begin{equation} 
	\label{Equation: Otwo-volutive strictification condition}
	\begin{aligned}
    \begin{matrix}
        \xymatrix{
        d(c) \ar[d]_-{S_{d(c)}} \ar[rr]^-{\id_{d(c)}} && d(c) \ar[d]^-{\epsilon_c} \ar@{=>}[dll]_-{\Gamma_c} \\
        d(c) \ar[d]_-{\theta_c^{-1}} \ar[rr]^-{d(S_c)} && d(c) \ar[d]^-{\theta_c^{-1}} \ar@{=>}[dll]_-{d(\lambda_c)}\\
        c \ar[rr]_-{\id_c}&& c
        }
    \end{matrix}
    \;=\; 
    \begin{matrix}
        \xymatrix{
            d(c) \ar[rrdd]^-{^{\vee \vee}\theta_c^{-1}}  \ar[d]_-{S_{d(c)}}  \ar[rr]^-{\id_{d(c)}} && d(c) \ar[d]^-{\epsilon_c} \ar@{}[dl]^(.25){}="a"^(.65){}="b" \ar@{=>}_{\phi_c} "a";"b"\\
            d(c) \ar[d]_-{\theta_c^{-1}} & \ar@{}[dl]^(.15){}="a"^(.55){}="b" \ar@{=>}_{S_{\theta_c}} "a";"b" \ar@{}[d]^(0.8){}="a"^(1.2){}="b" \ar@{=>}_{\lambda^{-1}_c} "a";"b" & d(c) \ar[d]^-{\theta_c^{-1}} \\
            c \ar@/^1pc/[rr]^-{S_c^{-1}} \ar@/_1pc/[rr]_-{\id_c} && c
        }
    \end{matrix}
	\end{aligned}.
\end{equation}
We construct an $\Otwo$-dagger structure on $S_{\Otwo}\B$. First, the $\Oone$-dagger and $\Sotwo$-volutive structures $(S_{\Otwo}d,S_{\Otwo}\eta)$ and 
$S_{\Otwo}S$ on $S_{\Otwo}\B$ are constructed as for $S_{\Oone} \B$ and $S_{\Sotwo}\B$, respectively. Next, note that the 2-transformation 
$\epsilon \colon ^{\vee \vee}(-) \circ d \cong d \circ (-)^{\vee \vee}$ induces a 2-transformation 
$\hat{\epsilon} \colon ^{\vee \vee}(-) \circ S_{\Otwo}d \cong S_{\Otwo}d \circ (-)^{\vee \vee}$ with 1-morphism component at $(c,\theta_c,\Pi_c,\lambda_c,\phi_c)$ 
given by $\hat{\epsilon}_c := \theta_c^{-1} \circ \epsilon_c \circ ^{\vee \vee}{\theta_c}$ and 2-morphism component at $X$ given by 
\begin{equation*}
	\begin{aligned}
    \xymatrix{
        b \ar[rr]^-{\theta_b} \ar[d]_-{\hat{\epsilon}_b} && d(b) \ar[rr]^-{d(X^{\vee \vee})} \ar[d]^-{\epsilon_b} \ar@{=>}[dll]^-{=} && d(a) \ar[rr]^-{\theta_a^{-1}} \ar[d]^-{\epsilon_a} \ar@{=>}[dll]^-{\epsilon_X} && a \ar[d]^-{\hat{\epsilon}_a} \ar@{=>}[dll]^-{=}\\
        b \ar[rr]_-{^{\vee \vee}\theta_b} && d(b) \ar[rr]_-{^{\vee \vee}d(X)} && d(a) \ar[rr]_-{^{\vee \vee}\theta_a^{-1}} && a
    }
	\end{aligned}
	.
\end{equation*}
This 2-transformation in turn induces an 2-transformation $S_{\Otwo}\epsilon \colon ^{\vee \vee}(-) \circ S_{\Otwo}d \cong S_{\Otwo}d \circ (-)^{\vee \vee}$ 
whose 1-morphism components are identities, and whose 2-morphism component at $X$ is given by 
\begin{equation}
    \xymatrix{
        \id_a \circ S_{\Otwo}d( X^{\vee \vee}) \ar[r]^-{\phi_a^{-1} \circ \id} & \hat{\epsilon}_a \circ  S_{\Otwo}d( X^{\vee \vee})  \ar[r]^-{\hat{\epsilon}_X} & ^{\vee \vee}(S_{\Otwo}d(X)) \circ \hat{\epsilon}_b \ar[r]^-{\id \circ \phi_b} & ^{\vee \vee}(S_{\Otwo}d(X)) \circ \id_b .
    }
\end{equation}
Again, one shows that this defines a 2-transformation by using the respective properties of $\hat{\epsilon}$. Finally, one deduces 
$(\id_{S_{\Otwo} d} \circ S_{\Otwo} S^{-1}) \bullet S_{\Otwo}\epsilon = S_{\Otwo} S \circ \id_{S_{\Otwo} d}$ from~\eqref{Equation: Otwo-volutive strictification condition}. 

\begin{lemma}
    $(S_{\Otwo}\B,  S_{\Otwo}d,  S_{\Otwo}\eta,  S_{\Otwo}S)$ is an $\Otwo$-dagger 2-category.
\end{lemma}

We expect that the assignment $(\B,d,\eta,\tau,S,\Gamma) \mapsto (S_{\Otwo}\B,  S_{\Otwo}d,  S_{\Otwo}\eta,  S_{\Otwo}S)$ extends to a 3-functor 
\begin{equation}
    S_{\Otwo} \colon \Otwovolbicat \xymatrix{\ar[r]&} \strictOtwovolbicat .
\end{equation}

\begin{remark}
    Assuming that $S_{\Otwo}$ and $T_{\Otwo}$ are indeed 3-functors, we expect that $S_{\Otwo}$ is right adjoint to $T_{\Otwo}$ in the 4-category of 3-categories. 
\end{remark}

\begin{remark}
    As indicated in \Cref{Remark: The 3-functor from symmetric monoidal 2-categories with duals to Otwo volutive 2-categories}, there should be a 3-functor 
    $\symmontwocat^{\mathrm{d}} \to \Otwovolbicat$ whose post-composition with $S_{\Otwo}$ should then yield a 3-functor 
    \begin{equation}
        \symmontwocat^{\mathrm{d}} \xymatrix{\ar[r]&} \strictOtwovolbicat
    \end{equation}
    assigning to each symmetric monoidal 2-category with duals and adjoints an $\Otwo$-dagger 2-category. 
\end{remark}

\subsection[$\Oone$-volutive categories from $\Otwo$-dagger 2-categories]{$\textrm{\textbf{O}}\pmb{(1)}$-volutive categories from $\textrm{\textbf{O}}\pmb{(2)}$-dagger 2-categories}
\label{Subsection: Otwovolutive 2-categories and Oonevolutive 1-categories}

In this section we show that every $\Otwo$-dagger 2-category comes with a canonical and compatible $\Oone$-volutive structure on its Hom categories. 

\medskip

\begin{definition}
    A \emph{bi-$\dagger$} structure $(d_1,d_2,\phi)$ on a 2-category $\B$ consists of two 2-functors $d_1 \colon \B \to \B^{\oneop}$, $d_2 \colon \B \to \B^{\twoop}$
    and an invertible 2-transformation $\phi \colon d_1^{\oneop} \circ d_1 \Rightarrow \id_{\B}$ such that 
    \begin{itemize}
    	\item 
    	$d_2$ is identity on objects and 1-morphisms and strictly squares to the identity, 
    	\item
    	$d_1$ is identity-on-objects, 
    	\item
    	$\phi$ has identity 1-morphism components, 
    	\item
    	there are equalities $(\phi^{\oneop})^{-1} \circ \id_{d_1} = \id_{d_1} \circ \phi$ and $d_2^{\oneop} \circ d_1 = d_1^{\twoop} \circ d_2$,
    	\item
    	$\phi$ is unitary with respect to $d_2$, i.e.\ $d_2(\phi_X) = \phi_X^{-1}$. 
    \end{itemize}
    A 2-category $\B$ together with a bi-$\dagger$ structure 
    is called \emph{bi-$\dagger$ 2-category}.
\end{definition}

We note that the pair $(d_1,\phi)$ in the definition of a bi-$\dagger$ structure is precisely the data of an $\Oone$-dagger structure on the 2-category $\B$. 
In other words, a bi-$\dagger$ 2-category is an $\Oone$-dagger 2-category whose Hom categories carry compatible dagger structures. 

\begin{remark}
    In \cite[Definition 4.1]{ferrer2024dagger}, bi-$\dagger$ 2-categories are called \emph{bi-involutive bicategories}. A more coherent analogue of this definition, 
    which in our terminology would be called a \emph{bi-volutive 2-category}, is called a \emph{fully dagger bicategory} in loc.\ cit., and a strictification procedure akin 
    to the ones we present here is proposed in \cite[Statement 4.5]{ferrer2024dagger}.
\end{remark}

In the following, we will construct bi-$\dagger$ 2-categories from $\Otwo$-dagger 2-categories. Let $(\B,d,\eta,S)$ be an $\Otwo$-dagger
2-category. We claim that the Hom categories of $\B$ carry canonical $\Oone$-volutions. Indeed, let $a,b \in \B$ be objects. 
Consider the functor 
$\hat{d} := d\circ(-)^\vee \colon \Hom_{\B}(a,b) \to \Hom_{\B}(a,b)^{\opp}$ and the natural isomorphism 
$\hat{\eta} \colon \hat{d}^{\opp} \circ \hat{d} \Rightarrow \id$ with component at $X$ given by 
\begin{equation}
	\hat{\eta}_X^{-1} := 
	\Big( 
    \xymatrix{
        d(^{\vee}(d(X^{\vee}))) \ar[r]^-{\cong} & d^2(X^{\vee \vee}) \ar[r]^-{d^2(S_X)} & d^2(X) \ar[r]^-{\eta_X^{-1}} & X
    }
	\Big) .
\end{equation}
One checks that the pair $(\hat{d},\hat{\eta})$ defines an $\Oone$-volution on $\Hom_{\B}(a,b)$. Consider the 2-category $\hat{\B}$ whose objects are those of 
$\B$ and whose Hom category at two objects $a,b \in \hat{B}$ is given by $S_{\Oone}\Hom_{\B}(a,b)$, which by definition comes with a dagger structure $S_{\Oone}\hat{d}$. 
One checks that this procedure, i.e.\ applying the strictification of $\Oone$-volutions locally, indeed results in a 2-category. We further note that the 
dagger structures on the Hom categories of $\hat{B}$ together form an identity-on-objects-and-1-morphisms functor $d_2 \colon \hat{\B} \to \hat{B}^{\twoop}$ which 
strictly squares to the identity. It is straightforward to construct from the $\Oone$-dagger structure $(d,\eta)$ on $\B$ an $\Oone$-dagger structure $(d_1,\phi)$ on $\hat{B}$. 
Moreover, one checks that $d_1$ and $d_2$ are compatible in the sense that $(\hat{\B},d_1,d_2,\phi)$ is a bi-$\dagger$ 2-category. We record this in the following. 

\begin{lemma}\label{Lemma: bi-dagger 2-categories from strict Otwo volutive 2-categories}
    Any $\Otwo$-dagger 2-category $(\B,d,\eta,S)$ gives rise to a bi-$\dagger$ 2-category $(\hat{\B},d_1,d_2,\phi)$ as described in the previous paragraph.
\end{lemma}

\begin{remark}\label{Remark: symmetric monoidal 2-categories with duals and adjoints give bi-dagger 2-categories}
    We expect that this assignment too will turn out to be functorial. Denoting the conjectural 3-category of bi-$\dagger$ 2-categories by $\strictbivolbicat$, 
    we summarize the chain of conjectural 3-categories and 3-functors as follows: 
    \begin{equation}
        \xymatrix{
            \symmontwocat^{\mathrm{d}} \ar[rr]^-{\text{\Cref{Remark: The 3-functor from symmetric monoidal 2-categories with duals to Otwo volutive 2-categories}}} && \Otwovolbicat \ar[rr]^-{S_{\Otwo}} 
            && \strictOtwovolbicat \ar[rr]^-{\text{\Cref{Lemma: bi-dagger 2-categories from strict Otwo volutive 2-categories}}} && \strictbivolbicat .
        }
    \end{equation}
    In particular, from any symmetric monoidal 2-category with duals and adjoints one constructs a bi-$\dagger$ 2-category in a canonical way. 
\end{remark}

\begin{remark}
    One reason to study bi-$\dagger$ 2-categories is their proposed application to the formulation of unitary 2-group representations; see 
    e.g.\ \cite{bartsch2024unitary2groupsymmetries}, where bi-$\dagger$ 2-categories are called \emph{$\dagger$-bicategories}. 
    \Cref{Remark: symmetric monoidal 2-categories with duals and adjoints give bi-dagger 2-categories} then provides a method to study 
    unitary 2-group representations on any symmetric monoidal 2-category with duals and adjoints, or rather its associated bi-$\dagger$ 2-category.
\end{remark}

As explained above, the Hom categories of any $\Otwo$-dagger 2-category carry canonical $\Oone$-volutive structures. The following illustrates this construction. 

\begin{example} 
	\label{Example: comparing Oone volutions on deloopings}
    Let $\C$ be any symmetric monoidal rigid category, for example $\C = \vect$. 
    Then, $\mathrm{B}\C$ is symmetric monoidal with (trivial) duals and adjoints (given by the duals in $\C$) 
    and hence carries an $\Otwo$-volutive structure by \Cref{Remark: The 3-functor from symmetric monoidal 2-categories with duals to Otwo volutive 2-categories}, 
    which one checks to be strict in the sense of \Cref{Definition: strict Otwo volutive 2-category}. The induced $\Oone$-volutive structure on the single Hom category 
    of $\mathrm{B}\C$ coincides with the $\Oone$-volutive structure on $\C$ induced by the symmetric monoidal rigid structure as explained in \Cref{Example: symmetric monoidal categories with duals are Oone volutive}.
\end{example}

\section{Algebraic description of orbifolds}
\label{sec:AlgebraicDescriptionOrbifolds}

In this section we review, reformulate and propose several categorifications of idempotent completion. 
In \Cref{subsec:IdempotentCompletion} we briefly recall the case for arbitrary 1-categories as well as one categorification that applies to arbitrary 2-categories. 
Categorification is not unique, and one expects different variants associated to different tangential structures. 
In \Cref{Subsection: Oone orbifold completion} we review an idempotent completion for $\Oone$-volutive 1-categories, explain that it restricts to the structure-preserving idempotent completion of dagger categories, and show how it interacts with the adjoint 2-functors~$S_{\Oone}$ and~$T_{\Oone}$ of \Cref{Subsection: Oonevolutive categories}. 
Along the same lines we describe idempotent completion of $G$-volutive 2-categories for~$G$ a spin group and~$\Otwo$ in Sections~\ref{Subsection: rSpintwo orbifold completion} and~\ref{Subsection: Otwo orbifold completion}, respectively. 
For $G=\Sotwo$ we explain how this systematic approach recovers (Euler completed) orbifold completion, while for~$\Otwo$ we thus produce a candidate for ``(a generalization of Euler completed) $\Otwo$-idempotent completion''.

\subsection{Idempotent completion}
\label{subsec:IdempotentCompletion}

In Section~\ref{Subsubsection: Idempotent completion of 1-categories} we recall the basic definitions and properties of idempotent completion of 1-categories. 
One particular categorification of idempotency and idempotent completion is reviewed in Section~\ref{Subsubsection: Idempotent completion of 2-categories}, namely the one believed to correspond to 2-framings.

\subsubsection{Idempotent completion of 1-categories}
\label{Subsubsection: Idempotent completion of 1-categories}

Recall that an \emph{idempotent} is an endomorphism that squares to itself. 

\begin{definition}
	An idempotent $e \colon a \to a$ in a given category $\C$ \emph{splits} if there is an object $b \in \C$ together with morphisms $Y \colon a \to b$ and $Z \colon b \to a$ such that $e = Z \circ Y$ and $\id_b = Y \circ Z$. 
	The triple $(b,Y,Z)$ is called a \emph{splitting} of~$e$, with~$b$ the \textsl{image} of~$e$. 
	If every idempotent in~$\C$ splits then~$\C$ is called \emph{idempotent complete}.
\end{definition}

There are several other useful charactizations of whether a splitting exists, see e.g.\ \cite[Proposition 6.5.4]{Borceux1994}: 

\begin{lemma}\label{lemma: Splittings of idempotents are limits}
	Let $e \colon a \to a$ be an idempotent. 
	The following are equivalent: 
	\begin{enumerate}
		\item[(i)] The idempotent $e$ splits as $e = Z\circ Y$ with $Y \colon a \to b$, $Z \colon b \to a$, i.e.\ $Y\circ Z = \id_b$. 
		\item[(ii)] 
		The limit~$b$ of the diagram $e \colon a \to a$ exists, with structure map $Z \colon b \to a$.  
		\item[(iii)] 
		The colimit~$b$ of the diagram $e \colon a \to a$ exists, with structure map $Y \colon a \to b$.  
	\end{enumerate}
	Moreover, under these conditions we have that both the limit and colimit are absolute (i.e., they are preserved by all functors whose domain is the ambient category of~$e$).
\end{lemma}
\begin{definition}
	Let $\C$ be a category. Its \emph{idempotent completion} $\Ide\C$ is the category whose objects are pairs $(a,e)$ where $a \in \C$ and $e \colon a \to a$ is an idempotent 
	in $\C$, and whose morphisms $(a,e) \to (b,f)$ are morphisms $X \colon a \to b$ in~$\C$ such that $f \circ X = X = X \circ e$. The identity morphism on $(a,e)$ is~$e$, and composition in $\Ide\C$ is the one of $\C$.
\end{definition}
\begin{remark}
	$\Ide\C$ is also called the \emph{Karoubi completion} of $\C$. 
\end{remark}
Idempotent completion satisfies the following universal property, see e.g.\ \cite[Appendix A.1]{Dcoppet2022}.
\begin{lemma}\label{Lemma: universal property of the idempotent completion}
	The functor $\C \to \Ide\C$ given by $a \mapsto (a,\id_a)$ on objects and by $X \mapsto X$ on morphisms is fully faithful. 
	Moreover, $\Ide\C$ is idempotent complete, and every other functor 
	$\C \to \D$ with idempotent complete codomain factors (essentially uniquely) through $\C \hookrightarrow \Ide\C$. 
\end{lemma}
\begin{corollary}
	Idempotent completion is a completion, in particular we have $\Ide\Ide\C \cong \Ide\C$ for any category $\C$. 
\end{corollary}
\begin{lemma}\label{Lemma: idempotent completion is a 2-functor}
	The assignment $\C \mapsto \Ide\C$ extends to a 2-functor 
    \begin{equation}
		\Idewo_1 \colon \Cat \xymatrix{\ar[r]&} \Cat^{\text{ic}} \;\subset\; \Cat
	\end{equation}
	from the 2-category of categories to its full sub-2-category of idempotent complete categories.
\end{lemma}
\begin{remark}
	\label{Remark: 1-idempotent completion as adjunction}
	Another way of characterizing the 2-functor $\Ide_1 \hspace{0.05cm} \colon \Cat \to \Cat^{\text{ic}}$ is by identifying it as a left adjoint to the inclusion 2-functor 
	$\iota \colon \Cat^{\text{ic}} \hookrightarrow \Cat$, cf.\ \cite[Proposition A.1.6]{Dcoppet2022}: 
	\begin{equation}
		\begin{tikzcd}[column sep=3em, >=stealth]
			\Cat   \arrow[rr, out=30, in=150, "\Ide_1"] 
			& \perp & 
			\Cat^{\text{ic}}   \arrow[ll, hookrightarrow, out=-150, in=-30, "\iota"] 
		\end{tikzcd}
	\end{equation}
\end{remark}

\begin{lemma}
	\label{Lemma: Idempontent completion of symmetric monoidal categories with duals}
	The 2-functor $\Ide_1 \hspace{0.05cm} \colon \Cat \to \Cat^{\text{ic}}$ extends to a 2-functor 
	\begin{equation}
		\Ide_{1,\otimes} \colon \symmoncat^{\textrm{d}} \xymatrix{\ar[r]&} \symmoncat^{\textrm{d},\text{ic}}
	\end{equation}
	from the 2-category of symmetric monoidal categories with duals to its full sub-2-category of idempotent complete symmetric monoidal categories with duals.
\end{lemma}

\subsubsection{Idempotent completion of 2-categories}
\label{Subsubsection: Idempotent completion of 2-categories}

One dimension higher, a categorification of idempotency and idempotent completion was introduced in \cite{18} and \cite{24}, see also \cite{81} for a detailed analysis.

To prepare for \Cref{Definition: IB} below, we review some standard algebraic notions, see e.g.\ \cite[Section~4]{18} for more details. 
An (associative unital) \textsl{algebra} in a 2-category~$\B$ is a 1-endomorphism $A\colon a \to a$ in~$\B$ together with 2-morphisms $\mu\colon A \otimes A \to A$ and $\eta\colon \id_a \to A$ such that $\mu \circ (\mu \otimes \id_A) = \mu \circ (\id_A \otimes \mu)$ and $\mu\circ(\eta\otimes \id_A) = \id_A = \mu\circ(\id_A\otimes\eta)$. 
A \textsl{right $A$-module} is a 1-morphism $X\colon a\to a'$ for some $a'\in\B$ together with a 2-morphism $\rho_{X,A} \colon X\otimes A \to X$ that is compatible with~$\mu$ and~$\eta$, and a left $A$-module structure $\rho_{A,X}$ is defined analogously. 
If $(a',A',\mu',\eta')$ is also an algebra, then an \textsl{$A'$-$A$-bimodule} is a 1-morphism $X\colon a \to a'$ together with a left $A'$- and a right $A$-module structure that commute with one another. 
If $\widetilde X\colon a\to a'$ also comes with the structure of an $A'$-$A$-bimodule, then an \textsl{intertwiner} $X \to \widetilde X$ is a 2-morphism between the underlying 1-morphisms that commutes with the $A$- and $A'$-actions on~$X$ and~$\widetilde X$. 

Dually, a (coassociative counital) \textsl{coalgebra} structure on~$A$ is given by 2-morphisms $\Delta \colon A \to A \otimes A$ and $\varepsilon\colon A \to \id_a$ satisfying the coassociativity and counitality conditions. 
A \textsl{Frobenius algebra} is an algebra $(a,A,\mu,\eta)$ and simultaneously a coalgebra $(a,A,\Delta,\varepsilon)$ with the same underlying 1-morphism~$A$ such that $(\id_A \otimes \mu) \circ (\Delta \otimes \id_A) = (\mu \otimes \id_A) \circ (\id_A \otimes \Delta)$. 
A Frobenius algebra is \textsl{$\Delta$-separable} if $\mu\circ\Delta = \id_A$. 
We often refer to ($\Delta$-separable) Frobenius algebras $(a,A,\mu,\eta,\Delta,\varepsilon)$ simply as~$A$, and similarly for modules. 

Given algebras $(a,A,\mu,\eta)$, $(a',A',\mu',\eta')$, $(a'',A'',\mu'',\eta'')$ together with an $A'$-$A$-bimodule~$Y$ and an $A''$-$A'$-bimodule~$X$, the \textsl{relative tensor product} $X\otimes_{A'}Y$ by definition is the coequalizer of the two maps $X\otimes A' \otimes Y\to X\otimes Y$ induced by the left $A'$-action on~$X$ and the right $A'$-action on~$Y$. 
If~$A'$ in fact comes with the structure of a $\Delta$-separable Frobenius algebra, then one checks that (where on the right we use the standard graphical calculus)
\begin{equation}
	\label{Equation: Idempotent for relative tensor product}
	e_{X,A',Y} 
	:= 
	(\rho_{X,A'} \otimes \rho_{A',Y})\circ \big( \id_X \otimes ( \Delta \circ \eta ) \otimes \id_Y \big) 
	= 
	\begin{tikzpicture}[very thick,scale=0.75,color=blue!50!black, baseline]	
		\draw (-0.9,-1) node[left] (X) {{\small$X$}};
		\draw (+0.9,-1) node[right] (Xu) {{\small$Y$}};
		\draw[color=green!50!black] (-0.05,0.4) node[right] (Xu) {{\small$A'$}};
		\draw (-1,-1) -- (-1,1); 
		\draw (1,-1) -- (1,1); 
		\fill[color=green!50!black] (-1,0.6) circle (2.5pt) node (meet) {};
		\fill[color=green!50!black] (1,0.6) circle (2.5pt) node (meet) {};
		\draw[-dot-, color=green!50!black] (0.35,-0.0) .. controls +(0,-0.5) and +(0,-0.5) .. (-0.35,-0.0);
		\draw[color=green!50!black] (0.35,-0.0) .. controls +(0,0.25) and +(-0.25,-0.25) .. (1,0.6);
		\draw[color=green!50!black] (-0.35,-0.0) .. controls +(0,0.25) and +(0.25,-0.25) .. (-1,0.6);
		\draw[color=green!50!black] (0,-0.75) node[Odot] (down) {}; 
		\draw[color=green!50!black] (down) -- (0,-0.35); 
		
	\end{tikzpicture} 
\end{equation}
is an idempotent, and the relative tensor product exists iff $e_{X,A',Y}$ splits, in which case $X\otimes_{A'}Y$ is given by the image of $e_{X,A',Y}$. 
The $A$-action on~$Y$ and the $A''$-action on~$X$ induce the structure of an $A''$-$A$-bimodule on $X\otimes_{A'}Y$. 

\medskip 

We are now in a position to consider our first categorification of ordinary idempotent completion: 

\begin{definition}
	\label{Definition: IB}
	Let $\B$ be a 2-category with idempotent complete Hom categories. Its \emph{idempotent completion} $\Ide\B$ is the 2-category consisting of the following data:
	\begin{itemize}
		\item objects in $\Ide\B$ are pairs $(a,A)$ with $a \in \B$ and $A \in \B(a,a)$ with the structure of a $\Delta$-separable Frobenius algebra,
		\item 1-morphisms $(a,A) \to (b,B)$ in $\Ide\B$ are $X \in \B(a,b)$ with the structure of a $B$-$A$-bimodule,
		\item 2-morphisms in $\Ide\B$ are 2-morphisms in $\B$ that are intertwiners,
		\item the composition of 1-morphisms $X \colon (a,A) \to (b,B)$ and  $Y \colon (b,B) \to (c,C)$ 
		is the relative tensor product $Y \otimes_B X \colon (a,A) \to (c,C)$, the composition of 2-morphisms is that of $\B$, and
		\item unitors and associators are induced from the ones in $\B$. 
	\end{itemize}
\end{definition}

\begin{remark}
	\begin{enumerate} 
		\item 
		Note that a (Frobenius) algebra $A\colon a \to a$ is a weaker version of an idempotent in the sense that the self-composition $A\otimes A$ need not be equal to~$A$, but the two are coherently related by the multiplication $A\otimes A \to A$. 
		This is a first indication that idempotent completion of 2-categories is a categorification of the homonymous notion for 1-categories. 
		\item 
		In \cite{18}, $\Ide\B$ is called the \emph{equivariant completion}, while in \cite{24} it is called the \emph{(co)unital condensation completion} of $\B$. 
		While categorification of idempotent completion is non-unique (as illustrated in the subsequent subsections), here we prefer the term ``idempotent completion'' for 2-categories, as it emphasizes what is being categorified, and because it  appears to be the most basic categorification of ordinary idempotent completion (as also illustrated in the subsequent subsections).
		\item 
		The assumption of idempotent complete Hom categories is made in \Cref{Definition: IB} so that the relative tensor product of bimodules exists, i.e.\ so that the idempotents in~\eqref{Equation: Idempotent for relative tensor product} split. 
		One can in fact show that every idempotent in Hom categories is of this form. 
		On the other hand, it is natural to ask for all Hom categories to be idempotent complete; this is part of the induction step in a recursive definition of idempotent completion of $n$-categories for all $n\geqslant 1$. 
		\item 
		Any 2-category $\B$ can \emph{locally} be idempotent completed, i.e.\ applying the 2-functor $\Idewo_1$ from \Cref{Lemma: idempotent completion is a 2-functor} to each 
		Hom category of $\B$ yields another 2-category with idempotent complete Hom categories into which~$\B$ can be included.
	\end{enumerate}
\end{remark} 

\begin{remark}
	\begin{enumerate} 
		\item 
		There is a notion of a ``splitting'' of a $\Delta$-separable Frobenius algebra that satisfies a universal property akin to \Cref{Lemma: universal property of the idempotent completion}; see \cite{Dcoppet2022} and \cite{carqueville2023orbifold} for some details. 
		\item 
		The Hom categories of $\Ide\B$ are again idempotent complete. The universal property of $\Ide\B$ then implies 
		that $\Ide\Ide\B \cong \Ide\B$ for any 2-category $\B$ with idempotent complete Hom categories. 
		\item 
		The 2-category $\Ide\B$ descends from a double category via the construction in \cite[Corollary~6.12]{25}. 
		This fact simplifies proving that $\Ide\B$ inherits symmetric monoidal structures, as well as duals, if~$\B$ is equipped with such structures. 
		\item 
		The assignment $\B \mapsto \Ide\B$ extends to a 3-functor 
		\begin{equation}
			\Idewo_2 \colon \BiCat \xymatrix{\ar[r]&} \BiCat^{\text{ic}}
		\end{equation}
		from the 3-category of 2-categories to the full sub-3-category of idempotent complete 2-categories, cf.\ \cite{Dcoppet2022}.
	\end{enumerate}
\end{remark}

Let $\B$ be a 2-category with right adjoints~$X^\vee$ for all 1-morphisms~$X$, see e.g.\ \cite[Section~2]{18} for details as well as our (standard) conventions. 
In the following, we want to extend the results from the pivotal setting (which by definition involves an isomorphism $(-)^{\vee\vee} \cong \id_{\B}$ of 2-functors) obtained in \cite[Section~4]{18} to the non-pivotal setting. 
Recall that if $A \in \B(a,a)$ is a $\Delta$-separable Frobenius algebra, then~$A^\vee$ canonically also carries the structure of a $\Delta$-separable Frobenius algebra whose structure maps are the adjoints of the structure maps of~$A$. 
Moreover, note that for any $B$-$A$-bimodule~$X$ in~$\B$, the right adjoint $X^\vee$ canonically carries the structure of an $A$-$B^{\vee\vee}$-bimodule, with left and right action given by 
\begin{equation} \label{Equation: Definition of the actions on the adjoint bimodule}
	\begin{tikzpicture}[very thick,scale=0.85,color=blue!50!black, baseline=0cm]
		\draw[line width=0] 
		(-0.5,-1.25) node[line width=0pt, color=green!50!black] (Algebra) {{\small $A\vphantom{X^\vee}$}}
		(0,1.25) node[line width=0pt] (A) {{\small $X^\vee$}}
		(0,-1.25) node[line width=0pt] (A2) {{\small $X^\vee$}};
		\draw (A) -- (A2);
		\draw[color=green!50!black] (0,0) .. controls +(-0.5,-0.25) and +(0,1) .. (Algebra);
		\fill[color=green!50!black] (0,0) circle (2.5pt) node[left] {};
	\end{tikzpicture}
	\!:=\!
	\begin{tikzpicture}[very thick,scale=0.85,color=blue!50!black, baseline=0cm]
		\draw[line width=0] 
		(0.5,-1.25) node[line width=0pt, color=green!50!black] (Algebra) {{\small $A\vphantom{X^\vee}$}}
		(-1,1.25) node[line width=0pt] (A) {{\small $X^\vee$}}
		(1,-1.25) node[line width=0pt] (A2) {{\small $X^\vee$}}; 
		\draw[redirected] (0,0) .. controls +(0,-1) and +(0,-1) .. (-1,0);
		\draw[redirected] (1,0) .. controls +(0,1) and +(0,1) .. (0,0);
		\draw[color=green!50!black] (0,0) .. controls +(0.5,-0.25) and +(0,1) .. (Algebra);
		\fill[color=green!50!black] (0,0) circle (2.5pt) node[left] {};
		\draw (-1,0) -- (A); 
		\draw (1,0) -- (A2); 
	\end{tikzpicture}
	, \qquad 
	\begin{tikzpicture}[very thick,scale=0.85,color=blue!50!black, baseline=0cm]
		\draw[line width=0] 
		(+0.5,-1.25) node[line width=0pt, color=green!50!black] (Algebra) {{\small $B^{\vee\vee}\vphantom{X^\vee}$}}
		(0,1.25) node[line width=0pt] (A) {{\small $X^\vee$}}
		(0,-1.25) node[line width=0pt] (A2) {{\small $X^\vee$}};
		\draw (A) -- (A2);
		\draw[color=green!50!black] (0,0) .. controls +(+0.5,-0.25) and +(0,1) .. (Algebra);
		\fill[color=green!50!black] (0,0) circle (2.5pt) node[left] {};
	\end{tikzpicture}
	\!:=\!
	\begin{tikzpicture}[very thick,scale=0.85,color=blue!50!black, baseline=0cm]
		\draw[line width=0] 
		(1.5,-1.25) node[line width=0pt, color=green!50!black] (Algebra) {{\small $B^{\vee\vee}\vphantom{X^\vee}$}}
		(-1,1.25) node[line width=0pt] (A) {{\small $X^\vee$}}
		(1,-1.25) node[line width=0pt] (A2) {{\small $X^\vee$}}; 
		\draw[redirected] (0,0) .. controls +(0,-1) and +(0,-1) .. (-1,0);
		\draw[redirected] (1,0) .. controls +(0,1) and +(0,1) .. (0,0);
		\draw[color=green!50!black] (0,0) .. controls +(0.15,0.15) and +(0.2,0) .. (-0.25,-0.25);
		\draw[color=green!50!black] (-0.25,-0.25) .. controls +(-0.6,0) and +(-0.85,0) .. (0.5,1.25);
		\draw[color=green!50!black] (0.5,1.25) .. controls +(0.7,0) and +(0,1) .. (1.5,-0.25);
		\draw[<-,color=green!50!black] (-0.21,-0.25) -- (-0.22,-0.25);
		\draw[>-,color=green!50!black] (0.41,1.25) -- (0.42,1.25);
		\draw[color=green!50!black] (1.5,-0.25) -- (Algebra); 
		\fill[color=green!50!black] (0,0) circle (2.5pt) node[left] {};
		\draw (-1,0) -- (A); 
		\draw (1,0) -- (A2); 
	\end{tikzpicture}
	.
\end{equation}

\begin{definition}
	The \emph{(right) Nakayama map} $\widetilde{\gamma}_A \colon A \to A^{\vee \vee}$ of a $\Delta$-separable Frobenius algebra~$A$ is 
	\begin{equation} 
		\label{The Nakyama morphism and its inverse}
		\widetilde{\gamma}_A 
		:= 
		( \widetilde\ev_{A^\vee} \otimes \id_{A^{\vee\vee}} )
		\circ 
		\big(\id_{A^{\vee}} \otimes (\Delta^{\vee\vee}\circ\eta^{\vee\vee}\circ\varepsilon\circ\mu) \big)
		\circ 
		(\widetilde\coev_A \otimes \id_A)
		=
		\begin{tikzpicture}[scale = 0.333, baseline=0.7cm]
			\draw[color=green!50!black, ultra thick, -,] (0,0) arc (0:180:1);
			\draw[color=green!50!black, ultra thick, -] (-1,1) -- (-1,1.9);
			\draw[color=green!50!black, thick] (-1,2) node[Odot] (unit) {}; 
			\filldraw[color=green!50!black] (-1,1) circle (7pt);
			
			\begin{scope}[yscale=-1, yshift =-5cm]
				\draw[color=green!50!black, ultra thick, -,] (0,0) arc (0:180:1);
				\draw[color=green!50!black, ultra thick, -] (-1,1) -- (-1,1.9);
				\filldraw[color=green!50!black] (-1,1) circle (7pt);
				\draw[color=green!50!black, thick] (-1,2) node[Odot] (unit) {}; 
			\end{scope}
			
			\draw[color=green!50!black, ultra thick, -] (0,0) -- (0,-2);
			\draw[color=green!50!black, ultra thick, -] (0,5) -- (0,7);
			
			\begin{scope}[very thick,decoration={
					markings,
					mark=at position 0.5 with {\arrow{>}}}
				] 
				\draw[color=green!50!black, ultra thick, -, postaction={decorate}] (-4,0) arc (-180:0:1);
			\end{scope}
			\begin{scope}[very thick,decoration={
					markings,
					mark=at position 0.5 with {\arrow{<}}}
				] 
				\draw[color=green!50!black, ultra thick, -, postaction={decorate}] (-2,5) arc (0:180:1);
			\end{scope}
			\draw[color=green!50!black, ultra thick, -] (-4,0) -- (-4,5);
			
			\node[color=green!50!black] at (0.5,-2)   {$A$};
			\node[color=green!50!black] at (-4.7,2.5)   {$A^{\vee}$};
			\node[color=green!50!black] at (0.95,7)   {$A^{\vee \vee}$};
		\end{tikzpicture}
	\end{equation}
\end{definition}

We comment on the relation between the above Nakayama map and the standard ``Nakayama automorphism'' in \Cref{Remark: The Nakayama automorphism} below.
Given this connection, the following result comes with little surprise. 

\begin{lemma}
	The Nakayama map $\widetilde{\gamma}_A$ is an isomorphism of Frobenius algebras. 
	Moreover, $\widetilde{\gamma}_A^{\vee \vee} = \widetilde{\gamma}_{A^{\vee \vee}}$.
\end{lemma}
\begin{proof}
	The arguments in the pivotal case, see e.g.\ \cite{88}, carry over to the non-pivotal case if one makes sure that~$A^{\vee\vee}$ is not confused with~$A$. 
\end{proof}

Let $B,B'$ be algebras, and let $\varphi\colon B'\to B$ be an algebra isomorphism. 
Recall that if~$Y$ comes with the structure of a right $B$-module, then the \textsl{$\varphi$-twisted right $B'$-module $Y_\varphi$} has the same underlying 1-morphism~$Y$ and its $B'$-action is the original $B$-action pre-composed with $\id_Y \otimes \varphi$. 

\begin{lemma} \label{Lemma: right adjoints in the idempotent completion}
	Let $\B$ be a 2-category with idempotent complete Hom categories and right adjoints. Then $\Ide\B$ also has right adjoints. The right adjoint of a 1-morphism 
	$X \colon (a,A) \to (b,B)$ in $\Ide\B$ is $X^\star := (X^\vee)_{\widetilde{\gamma}_B}$.
\end{lemma}
\begin{proof}
	The proof is formally analogous to the one in \cite[Proposition~4.7]{18} if one replaces the Nakayama automorphism by the Nakayama map.
\end{proof}
The theory of left adjoints for morphisms in $\Ide \B$ is analogous. We denote the respective left Nakayama map $A \to {}^{\vee \vee}A$ by the same symbol $\widetilde{\gamma}_A$.
It will always we be clear from the context which one is meant.

\subsection[$\Oone$-idempotent completion]{$\textrm{\textbf{O}}\pmb{(1)}$-idempotent completion}
\label{Subsection: Oone orbifold completion}

In this section we recall a variant of idempotent completion for dagger categories, lift it to $\Oone$-volutive categories, and work out how it composes with the 2-functors $S_{\Oone}, T_{\Oone}$ of \Cref{Subsection: Oonevolutive categories}. 

\medskip 

The following definitions were put forward in \cite{SELINGER2008107}. 
Let $(\C,d)$ be a dagger category as in \Cref{def:DaggerCategory}. 
A morphism $e \colon a \to a$ in $\C$ is called a \emph{$d$-idempotent} 
if $e \circ e = e$ and $d(e) = e$ holds. 

\begin{definition}
    Let $(\C,d)$ be a dagger category. Let $\overline{\C}^d \subset \Ide \C$ be the full subcategory whose objects are pairs $(a,e)$ in which 
    $e\colon a\to a$ is a $d$-idempotent. 
    The category~$\overline{\C}^d$ together with the dagger structure (denoted by the same symbol!) $d \colon \overline{\C}^d \to ({\overline{\C}^d})^{\opp}$, $(a,e) \mapsto (a,e)$, $X \mapsto d(X)$ is called the \emph{$d$-Karoubi envelope} of $(\C,d)$.
\end{definition}
\begin{remark}
    The $d$-Karoubi envelope satisfies a universal property analogous to the one in \Cref{Lemma: universal property of the idempotent completion}, where one replaces 
    the notion of splitting of an idempotent by the notion of a splitting of a $d$-idempotent. 
\end{remark}

We note that a dagger structure $d$ on $\C$ does not necessarily induce a dagger structure on $\Ide \C$, which is why one restricts to the full subcategory $\overline{\C}^d$. 
However, any dagger structure can also be reinterpreted as an $\Oone$-volutive structure with trivial higher coherence via the 2-functor $T_{\Oone}$ in~\eqref{eq:TO1}, and any $\Oone$-volutive structure on 
$\C$ canonically induces an $\Oone$-volutive structure on $\Ide \C$, which we  describe in the following. 

\begin{lemma} \label{Lemma: Oone volution on idempotent completion}
    Let $\C$ be a category. 
    Any $\Oone$-volutive structure $(d,\eta)$ on $\C$ induces an $\Oone$-volutive structure $(d',\eta')$ on $\Ide \C$. Moreover, the assignment 
    $(\C,d,\eta) \mapsto (\Ide \C,d',\eta')$ extends to a 2-functor 
    \begin{equation}
        \Ide_{\Oone} \colon \Oonevolcat \xymatrix{\ar[r]&} \Oonevolcat.
    \end{equation}
\end{lemma}
\begin{proof}
    We have the functor $d' \colon \Ide \C \to \Ide \C^{\opp}$, $(a,e) \mapsto (d(a),d(e))$, $X \mapsto d(X)$ which is clearly well-defined. 
    We define the component of the 
    natural isomorphism $\eta' \colon {d'}^{\opp} \circ d' \to \id$ at $(a,e)$ to be the morphism $\eta_{(a,e)}' := e \circ \eta_a \colon (d^2(a),d^2(e)) \to (a,e)$. 
    This is well-defined since by naturality of $\eta$ we have $e \circ \eta_a = \eta_a \circ d^2(e)$ and hence 
    \begin{equation}
        e \circ (e \circ \eta_a) = e \circ \eta_a = \eta_a \circ d^2(e) = \eta_a \circ d^2(e)\circ d^2(e) = (e \circ \eta_a) \circ d^2(e).
    \end{equation}
    Naturality follows from the respective property of $\eta$, and $\eta_{(a,e)}'$ is invertible with inverse given by $\eta_a^{-1} \circ e$. It is also routine to check 
    that the coherence condition of an $\Oone$-volution is satisfied. 
    
    Next, let $(F,\alpha) \colon (\C_1,d_1,\eta_1) \to (\C_2,d_2,\eta_2)$ be an $\Oone$-volutive functor. We construct an $\Oone$-volutive functor 
    $(\Ide \C_1,d_1',\eta_1') \to (\Ide \C_2,d_2',\eta_2')$. As the underlying functor we take $\Ide F \colon \Ide \C_1 \to \Ide \C_2$, which assigns 
    $(a,e) \mapsto (F(a),F(e))$ and $X \mapsto F(X)$. The natural isomorphism $\alpha' \colon \Ide F^{\opp} \circ d_1' \Rightarrow d_2' \circ \Ide F$ 
    has as its 1-morphism component at $(a,e) \in \Ide \C_1$ the isomorphism $\alpha'_{(a,e)} := d(F(e)) \circ \alpha_a$. Finally, to an $\Oone$-volutive natural transformation 
    $\varphi \colon (F,\alpha) \to (\hat{F},\hat{\alpha})$ we assign the  $\Oone$-volutive natural transformation $(\Ide F,\alpha') \to (\Ide \hat{F},\hat{\alpha}')$ 
    given by $\Ide \varphi$. One checks that this assignment indeed defines a 2-functor. 
\end{proof}
The following is a straightforward comparison of induced structures.
\begin{lemma} \label{Lemma: Oone volutive idempotent completion is compatible with symmetric monoidal origin}
    The diagram of 2-functors 
    \begin{equation}
        \xymatrix{
            \symmoncatd \ar[d] \ar[rr]^-{\Ide_{1,\otimes}} &&  \symmoncatd \ar[d] \\
            \Oonevolcat \ar[rr]_-{\Ide_{\Oone}} &&  \Oonevolcat
        }
    \end{equation}
    where the vertical 2-functors are the ones of \Cref{Proposition: From symmetric monoidal categories to Oone volutive categories} and the upper 2-functor 
    is the one obtained (in \Cref{Lemma: Idempontent completion of symmetric monoidal categories with duals}) by extending $\Ide$ to the 2-category of rigid symmetric monoidal categories, is commutative. 
\end{lemma}

\begin{lemma} \label{Lemma: Oone idempotent completion is idempotent}
    There is an invertible 2-transformation $\Ide_{\Oone} \circ \Ide_{\Oone} \cong \Ide_{\Oone}$.
\end{lemma}

\begin{proof}
    We describe the 1-morphism component of the 2-transformation at an $\Oone$-volutive category $(\C,d,\eta)$. Consider the equivalence of categories $I \C \cong II\C$ 
    given on object level by $(a,e) \mapsto (a,e,e)$. We claim that this functor induces an $\Oone$-volutive functor $\Ide_{\Oone}(\C,d,\eta) \to 
    (\Ide_{\Oone}\Ide_{\Oone})(\C,d,\eta)$. Indeed, the natural isomorphism $\alpha$ in the definition of an $\Oone$-volutive functor can here be taken to be trivial. 
    We infer from \cite[Lemma 2.5]{Stehouwer2023DaggerCV} that this $\Oone$-volutive functor is then an equivalence in $\Oonevolcat$. The 2-morphism component of 
    the 2-transformation at an $\Oone$-volutive functor, which consists of an $\Oone$-volutive natural transformation, can also be taken to be trivial.
\end{proof}

Recall the functor~$S_{\Oone}$ from~\eqref{eq:SO1} and consider the composition 
\begin{equation}
    \xymatrix{
        \Daggercat \ar[r]^-{T_{\Oone}} & \Oonevolcat \ar[r]^-{\Ide_{\Oone}} & \Oonevolcat \ar[r]^-{S_{\Oone}} & \Daggercat.
    }
\end{equation}
Given a dagger category $(\C,d)$, we want to give an explicit description of the dagger category $(S_{\Oone} \circ \Ide_{\Oone} \circ T_{\Oone})(\C,d)$. 
First, we recall that the $\Oone$-volutive category $T_{\Oone}(\C,d)$ consists of the category $\C$ together with the $\Oone$-volutive structure $(d,\id)$. 
Next, the $\Oone$-volutive category $(\Ide_{\Oone} \circ T_{\Oone})(\C,d)$ 
consists of the category $\Ide_{\Oone} \C$ together with the $\Oone$-volutive structure whose functor assigns $(a,e) \mapsto (a,d(e))$, $X \mapsto d(X)$ and whose natural isomorphism 
is again trivial, recalling that $\id_{(a,e)} = e$. 
Lastly, the dagger category $(S_{\Oone} \circ \Ide_{\Oone} \circ T_{\Oone})(\C,d)$ has as its underlying category the one whose 
objects are triples $(a,e,\theta_{(a,e)})$ where $a \in \C$ is an object, $e \colon a \to a$ is an idempotent, and $\theta_{(a,e)} \colon (a,e) \to (a,d(e))$ is an isomorphism 
satisfying $d(\theta_{(a,e)})^{-1} \circ \theta_{(a,e)} = \id_{(a,e)}$, and whose morphisms $X \colon (a,e,\theta_{(a,e)}) \to (b,f,\theta_{(b,f)})$ are morphisms $X \colon a \to b$ 
in $\C$ satisfying $f \circ X = X = X \circ e$. 
The associated dagger structure on this category assigns to a morphism $X \colon (a,e,\theta_{(a,e)}) \to (b,f,\theta_{(b,f)})$ 
the morphism 
\begin{equation}
	\theta_{(a,e)}^{-1} \circ d(X) \circ \theta_{(b,f)} \colon (b,f,\theta_{(b,f)}) 
	\xymatrix{\ar[r]&} 
	(a,e,\theta_{(a,e)}).  
\end{equation}
Given this explicit description of the dagger category $(S_{\Oone} \circ \Ide_{\Oone} \circ T_{\Oone})(\C,d)$, we find the following. 

\begin{proposition}
	\label{prop:FFdagger}
    Let $(\C,d)$ be a dagger category. 
    There is a fully faithful dagger functor 
    \begin{equation}
        \psi \colon (\overline{\C}^d,d) \xymatrix{\ar[r]&} (S_{\Oone} \circ I_{\Oone} \circ T_{\Oone})(\C,d).
    \end{equation}
\end{proposition}

\begin{proof}
    We send an object $(a,e) \in \overline{\C}^d$ to the object $(a,e,\id_{(a,e)}) \in (S_{\Oone} \circ \Ide_{\Oone} \circ T_{\Oone})\C$ and a morphism $X \colon (a,e) \to (b,f) \in 
    \overline{\C}^d$ to the morphism $X \colon (a,e,\id_{(a,e)}) \to (b,f,\id_{(b,f)})$. 
    It is clear that this is well-defined, functorial, and fully faithful. 
    Our explicit 
    description of the dagger structure on $(S_{\Oone} \circ \Ide_{\Oone} \circ T_{\Oone})\C$ shows that the functor is a dagger functor. 
\end{proof}

\begin{remark}
	\label{rem:SITnotStrongIdempotent}
    We claim that $S_{\Oone} \circ \Ide_{\Oone} \circ T_{\Oone} \colon \Daggercat \to \Daggercat$ allows for the structure of a monoid in the 2-category 
    $\Hom_{\BiCat}(\Daggercat,\Daggercat)$. Indeed, the component of the unit at a dagger category $(\C,d)$ is the inclusion dagger functor 
    $\iota \colon (\C,d) \to S_{\Oone} \circ \Ide_{\Oone} \circ T_{\Oone}(\C,d)$, and the associative multiplication is given by
    \begin{equation}
        \xymatrix{
            S_{\Oone} \circ \Ide_{\Oone} \circ T_{\Oone} \circ S_{\Oone} \circ \Ide_{\Oone} \circ T_{\Oone} \ar[rr]^-{\id \circ K_{\Oone} \circ \id} &&
            S_{\Oone} \circ \Ide_{\Oone} \circ \Ide_{\Oone} \circ T_{\Oone} \ar[r]^-{\cong} & S_{\Oone} \circ \Ide_{\Oone} \circ T_{\Oone}
        }
    \end{equation}
    where $K_{\Oone}$ denotes the counit of the adjunction between $T_{\Oone}$ and $S_{\Oone}$, and the last equivalence is described in 
    \Cref{Lemma: Oone idempotent completion is idempotent}. In general, we do not expect the multiplication map defined in this way to be an equivalence.
\end{remark}

\subsection[$\rSpintwo$-idempotent completion]{$\textrm{\textbf{Spin}}\pmb{(2)^r}$-idempotent completion}
\label{Subsection: rSpintwo orbifold completion}

In this section we lift $\rSpintwo$-volutive structures from 2-categories~$\B$ to their idempotent completion, producing $\rSpintwo$-volutive 2-categories $I_{\rSpintwo}\B$. 
The main result of this section then concerns the case $r=1$: we show that conjugating with the operations $S_{\Sotwo}, T_{\Sotwo}$ of \Cref{Subsection: rspintwovolutive 2-categories} recovers the ``Euler completed orbifold completion'' of~$\B$. 

\medskip 

Throughout this section, let $r \in \Z_{\geqslant 0}$. 
We claim that any $\rSpintwo$-volutive structure $S$ on a locally idempotent complete 2-category $\B$ with right adjoints extends to its idempotent completion~$\Ide \B$ of \Cref{Definition: IB}. 
To see this, 
the following terminology will be helpful. Let $(a,A) \in \Ide \B$. The \emph{Nakayama$^r$ morphism of $A$} is the Frobenius algebra isomorphism 
\begin{equation}
    \widetilde\gamma_A^r := \widetilde{\gamma}_{A^{(\vee \vee)^{r-1}}} \circ \dotsc \circ \widetilde{\gamma}_A \colon A \xymatrix{\ar[r]&} A^{(\vee \vee)^r}
\end{equation}
where~$\widetilde\gamma_A$ is the Nakayama map of~\eqref{The Nakyama morphism and its inverse}. 

\begin{remark} \label{Remark: The Nakayama automorphism}
    Consider a 2-category $\B$ with right adjoints and an $\Sotwo$-dagger structure $S \colon \id \Rightarrow (-)^{\vee \vee}$ (recall that this is equivalently 
    a pivotal structure on $\B$). Let $(a,A) \in \Ide \B$. The composition $S_A \circ \widetilde\gamma_A^1 \colon A \to A$ is commonly known as the \emph{Nakayama automorphism of $A$}, see e.g.\ \cite[Equation (4.3)]{18}. 
    In the general case, for $r \geqslant 1$ and a $\rSpintwo$-dagger structure $S \colon \id \Rightarrow ((-)^{\vee \vee})^r$, we will call the composition $S_A \circ \widetilde\gamma_A^r \colon A \to A$ the \emph{Nakayama$^r$ automorphism of $A$}.
\end{remark}

The following result is obtained from the (graphical) definition of the Nakyama map. 

\begin{lemma} \label{Lemma: The r-th power of the Nakayama automorphism in the Sotwo case}
    Let $(\B,S)$ be an $\Sotwo$-dagger 2-category and $(a,A) \in \Ide \B$. Taking the $r$-fold composition of $S$ with itself, we obtain a $\rSpintwo$-dagger structure 
    $S^r \colon \id \Rightarrow ((-)^{\vee \vee})^r$. 
    In this case, the Nakayama$^r$ automorphism of $A$ is the $r$-fold power of the Nakayama morphism of $A$, that is, 
    $S_A^r \circ \widetilde\gamma_A^r = (S_A \circ \widetilde\gamma_A^1)^r$.
\end{lemma}

Recall from \Cref{Lemma: right adjoints in the idempotent completion} that the idempotent completion $\Ide \B$ 
of a 2-category $\B$ with right adjoints $X^\vee$ also has right adjoints $X^\star$.   

\begin{lemma} 
	\label{Lemma: rSpintwo volution on idempotent completion}
    Let $\B$ be a locally idempotent complete 2-category with right adjoints.
    Any $\rSpintwo$-volutive structure $S$ on $\B$ induces a $\rSpintwo$-volutive structure $S' \colon \id \Rightarrow ((-)^{\star \star})^r$ on $\Ide \B$.
\end{lemma}

\begin{proof}
    We define the 1-morphism component of $S'$ at an object $(a,A) \in \Ide \B$ to be the 1-morphism $A_{(\widetilde\gamma_A^r)^{-1}} \otimes S_a$ with $A$-$A$-bimodule structure 
    given by 
    \begin{equation}
        \xymatrix{
            A \otimes(A_{(\widetilde\gamma_A^r)^{-1}} \otimes S_a) \otimes A \ar[r]^-{\id \otimes S_A^{-1}} & 
            A \otimes A_{(\widetilde\gamma_A^r)^{-1}} \otimes A^{(\vee \vee)^r} \otimes S_a \ar[r]^-{} &
            A_{(\widetilde\gamma_A^r)^{-1}} \otimes S_a
        }
    \end{equation}
    where we used the canonical $A$-$A^{(\vee \vee)^r}$-action on $A_{(\widetilde\gamma_A^r)^{-1}}$ in the last step. An inverse of $S_{(a,A)}'$ is given by 
    the $A$-$A$-bimodule $S_a^{-1} \otimes {_{(\widetilde\gamma_A^r)^{-1}}A}$ with bimodule actions constructed similarly as for $S_{(a,A)}'$. 
    
    In order to define the 2-morphism component of $S'$, we note that for any $A$-$B$-bimodule $X$ and any Frobenius algebra isomorphisms $\phi \colon A' \to A$ and $\psi \colon B' \to B$ we have 
    \begin{equation}
        (_{\phi}X_{\psi})^\vee \cong {}_{\psi}(X^\vee)_{\phi^{\vee \vee}}
    \end{equation}
    as $B'$-$(A')^{\vee \vee}$-bimodules; to see this, one may use the graphical definition in~\eqref{Equation: Definition of the actions on the adjoint bimodule}. 
    Recall further that the right adjoint in $\Ide \B$ is given by $X^\star := (X^\vee)_{\tilde{\gamma}_B}$ and that 
    $\widetilde{\gamma}_A^{\vee \vee} = \widetilde{\gamma}_{A^{\vee \vee}}$ holds. 
    We define the 2-morphism component of $S'$ at $X$ to be the 2-isomorphism 
    \begin{equation*}
        \xymatrix{
            X^{(\star \star)^r} \otimes_A S_{(a,A)}'
            \cong {}_{\widetilde\gamma_B}{X^{(\vee \vee)^r}}_{\widetilde\gamma_A} \otimes_A (A_{(\widetilde\gamma_A^r)^{-1}} \otimes S_a) 
            \cong {}_{\widetilde\gamma_B}{X^{(\vee \vee)^r}} \otimes S_a 
            \ar[r]^-{S_X} &  {}_{\widetilde{\widetilde\gamma_B}}(S_b \otimes X)
            \cong (B_{(\widetilde\gamma_B^r)^{-1}} \otimes S_b) \otimes_B X
        }
    \end{equation*}
    which one checks to be an intertwiner of $B$-$A$-bimodules using the naturality of $S$. One moreover checks that this defines an $\Sotwo$-volutive structure on $\Ide \B$, 
    using the respective properties of $S$.
\end{proof}

\begin{remark}
    Recall from \Cref{Example: Symmetric monoidal 2-categories with duals and adjoints are rSpintwo volutive} that any symmetric monoidal 2-category $\B$ with duals, adjoints, 
    and idempotent complete Hom categories carries an $\Sotwo$-volutive structure given by the Serre automorphism~$\mathbb{S}$. 
    By \Cref{Lemma: rSpintwo volution on idempotent completion}, $\Ide \B$ inherits an $\Sotwo$-volutive structure from $\B$.
    Its 1-morphism component at $(a,A) \in \Ide \B$ is given by $S_{(a,A)} = A_{(\widetilde\gamma_A^1)^{-1}} \otimes S_a$. Recalling \Cref{Remark: The Nakayama automorphism} and using the 
    isomorphism described in \cite[Equation (4.11)]{75}, this reproduces the 1-morphism component of the Serre automorphism determined in \cite[Proposition 4.8]{75} 
    up to an erroneous appearence of ``$\otimes  A$''. We expect there to be an analogous commutative diagram as in 
    \Cref{Lemma: Oone volutive idempotent completion is compatible with symmetric monoidal origin}.
\end{remark}

\begin{remark}
    We expect that the assignment $(B,S) \mapsto (\Ide \B,S')$ extends to a 3-functor 
    \begin{equation}
        \Ide_{\rSpintwo} \colon \rSpintwovoltwocat \xymatrix{\ar[r]&} \rSpintwovoltwocat.
    \end{equation} 
    This should then give a sequence of 3-functors 
    \begin{equation}
    	\label{eq:SITSpin}
        \xymatrix{
            \strictrSpintwovoltwocat\ar[rr]^-{T_{\rSpintwo}} && \rSpintwovoltwocat \ar[rr]^-{\Ide_{\rSpintwo}} && \rSpintwovoltwocat \ar[rr]^-{S_{\rSpintwo}} && \strictrSpintwovoltwocat .
        }
    \end{equation}
\end{remark}

In the remainder of this section, we investigate aspects of the composition~\eqref{eq:SITSpin}. 
Given a $\rSpintwo$-dagger structure $S$ on a 2-category $\B$ with right adjoints and idempotent complete Hom categories, we want to compute the image 
of $(\B,S)$ under the (conjectural) 3-functor $S_{\rSpintwo} \circ I_{\rSpintwo} \circ T_{\rSpintwo}$. First, the $\rSpintwo$-volutive 2-category 
$T_{\rSpintwo}(\B,S)$ has the same underlying 2-category and $\rSpintwo$-volutive structure. Second, the 2-category 
underlying $(I_{\rSpintwo} \circ T_{\rSpintwo})(\B,S)$ is $\Ide \B$. The $\rSpintwo$-volutive structure $S'$ on $\Ide \B$ has 1-morphism components 
$S'_{(a,A)} = A_{(S_A \widetilde\gamma_A^r)^{-1}}$, as shown in the proof of \Cref{Lemma: rSpintwo volution on idempotent completion}. 
Finally, the 2-category underlying $(S_{\rSpintwo} \circ I_{\rSpintwo} \circ T_{\rSpintwo})(\B,S)$ has objects given by 
triples $(a,A,\lambda_{(a,A)})$ where $(a,A) \in \Ide \B$ and $\lambda_{(a,A)} \colon A_{(S_A \widetilde\gamma_A^r)^{-1}} \to A$ is a 2-isomorphism in $\Ide \B$, while 1- and 2-morphisms 
are those of $\Ide \B$. The $\rSpintwo$-dagger structure on this 2-category is then given by 
\begin{equation*}
    \xymatrix{
        X^{(\star \star)^r} \cong X^{(\star \star)^r} \otimes_A A \ar[rr]^-{\id \otimes_A \lambda_{(a,A)}^{-1}} && X^{(\star \star)^r} \otimes_A A_{(S_A \widetilde\gamma_A^r)^{-1}} \ar[r]^-{S_X'} & 
        B_{(S_B^{-1} \widetilde\gamma_B^r)^{-1}} \otimes_B X \ar[rr]^-{\lambda_{(b,B)} \otimes_B \id} && B \otimes_B X \cong X.
    }
\end{equation*}

\begin{definition}\label{Definition: Euler completion}
    Let $(\B,S)$ be an $\Sotwo$-dagger 2-category. 
    Its \emph{Euler completion} $(\B,S)_{\eu}$ consists of the 2-category whose objects are pairs $(a,\psi_a)$ where $\psi_a \colon \id_a \to \id_a$ is a 2-isomorphism in~$\B$, and whose 1- and 2-morphisms are those of $\B$, and 
    the $\Sotwo$-dagger structure given by 
    \begin{equation}
        \xymatrix{
            X^{\vee \vee} \cong X^{\vee \vee} \circ \id_a \ar[r]^-{\id \circ \psi_a^{-1}} & X^{\vee \vee} \circ \id_a \ar[r]^-{S_X} & 
            \id_b \circ X \ar[r]^-{\psi_b \circ \id} & \id_b \circ X \cong X.
        }
    \end{equation}
\end{definition}

\begin{remark}
    This definition is closely related to the Euler completion of a pivotal (i.e., $\Sotwo$-dagger) 2-category $(\B,S)$ described in \cite[Section 5.1.1]{carqueville2023orbifold} and coming from \cite{6}. 
    Changing the pivotal structure according to our prescription and then computing the left adjoints yields precisely the ones given in \cite[Equation (5.1)]{carqueville2023orbifold}. 
    The difference between the procedures is that we do not change the right adjoints, whereas loc.\ cit.\ does. 
    This distinction is not crucial, for the following reason. 
    Given two choices of right adjunction data, there is an invertible 2-transformation $\chi \colon (-)^{\vee} \Rightarrow (-)^{\tilde{\vee}} $ whose 1-morphism components are identities.
    The composition 
    \begin{equation}
        \xymatrix{
            \id  \ar[r]^-{S} & (-)^{\vee \vee} \ar[r]^-{\chi \circ \chi} & (-)^{\tilde{\vee}\tilde{\vee}}
        }
    \end{equation}
    is then again an $\Sotwo$-dagger structure, which a priori might be different from~$S$.
    However, if we change our right adjoints according to the prescription in \cite[Equation (5.1)]{carqueville2023orbifold}, one computes that 
    $(-)^{\tilde{\vee}\tilde{\vee}} = (-)^{\vee \vee}$ and hence the pivotal structure remains the same. 
\end{remark}

\begin{remark}
    Recall from \cite{55} that any two-dimensional oriented defect TQFT $\mathscr{Z}$ gives rise to an $\Sotwo$-dagger (i.e.\ pivotal) 2-category $D_{\mathscr{Z}}$.
    The Euler completion of an $\Sotwo$-dagger 2-category was introduced in order to describe the $\Sotwo$-dagger 2-category 
    associated to the Euler completion $\mathscr{Z}_{\eu}$ of the defect TQFT $\mathscr{Z}$; namely, one has $D_{\mathscr{Z}_{\eu}} \cong (D_{\mathscr{Z}})_{\eu}$; see 
    \cite[Section 5]{carqueville2023orbifold}. 
\end{remark}

The following is obtained by comparison of definitions, explaining how the notion of Euler completion fits into our picture. We note that our theory immediately
suggests a generalization to the case of $r \geqslant 2$. 

\begin{lemma} \label{Lemma: Euler completion from our perspective}
    Let $(\B,S)$ be an $\Sotwo$-dagger 2-category. 
    Then 
    \begin{equation}
    	(\B,S)_{\eu} \cong (S_{\Sotwo} \circ T_{\Sotwo})(\B,S) .
    \end{equation}
\end{lemma}

Let $(\B,S)$ be an $\Sotwo$-dagger 2-category with idempotent complete Hom categories. 
Recall that a Frobenius algebra $A$ is called \emph{symmetric} if 
$A \cong A_{S_A \circ \widetilde\gamma_A^1}$ or equivalently $S_A \circ \widetilde\gamma_A^1 = \id_A$, see \cite[Definition 17]{88} and \cite[Equation (4.11)]{75}. 
Recall also the idempotent completion $\Ide\B$ of \Cref{Definition: IB}, the induced $\rSpintwo$-volutive structure $S'$ on $\Ide \B$ of 
\Cref{Lemma: rSpintwo volution on idempotent completion} for a $\rSpintwo$-volutive 2-category $(\B,S)$, and the equivalent description of $\Sotwo$-dagger 
(that is, pivotal) structures in \Cref{Remark: Equivalent ways of talking about pivotal structures}.

\begin{definition}[\cite{18}]
    The \emph{orbifold completion} $(\B,S)_{\orb}$ of $(\B,S)$ consists of the full sub-2-category of $\Ide\B$ whose objects are pairs $(a,A)$ with $a \in \B$ and $A \in \B(a,a)$ a $\Delta$-separable symmetric Frobenius algebra, as well as 
    the $\Sotwo$-dagger structure $S'$ induced from~$(\B,S)$.
\end{definition}

\begin{remark}
    Consider again a two-dimensional oriented defect TQFT $\mathscr{Z}$ giving rise to an $\Sotwo$-dagger (i.e.\ pivotal) 2-category $D_{\mathscr{Z}}$. 
    The orbifold completion of an $\Sotwo$-dagger 2-category was introduced in order to describe the $\Sotwo$-dagger 2-category
    associated to the orbifolded two-dimensional oriented defect TQFT $\mathscr{Z}_{\orb}$; namely, one has $D_{\mathscr{Z}_{\orb}} \cong (D_\mathscr{Z})_{\orb}$; see 
    \cite{carqueville2023orbifoldstopologicalquantumfield} and references therein. 
\end{remark}

We have the following comparison result, explaining how the orbifold completion fits into our picture. 

\begin{theorem} \label{Theorem: comparison of Euler orbifold completion and SIT}
    Let $(\B,S)$ be an $\Sotwo$-dagger 2-category with idempotent complete Hom categories. 
    There is an equivalence of $\Sotwo$-dagger 2-categories 
    \begin{equation}
        (S_{\Sotwo} \circ I_{\Sotwo} \circ T_{\Sotwo})(\B,S) \cong ((\B,S)_{\orb})_{\eu}
    \end{equation}
    between our construction and the Euler completed orbifold completion. 
\end{theorem}
\begin{proof}
    We define a 2-functor from the right-hand side to the left-hand side. 
    Let $(a,A,\psi_{(a,A)})$ be an object in $(\B_{\orb})_{\eu}$, i.e.\ $a \in \B$, $A \in \Hom_{\B}(a,a)$ is a $\Delta$-separable symmetric Frobenius algebra, and $\psi_{(a,A)} \colon A \to A$ an invertible $A$-$A$-bimodule intertwiner. 
    Since $A$ is symmetric, we have $S_A \circ \widetilde\gamma_A^1 = \id_A$ and hence $\psi_{(a,A)}$ gives a 2-isomorphism $A_{(S_A \circ \widetilde\gamma_A^1)^{-1}} = A \to A$. 
    Thus we may send the triple $(a,A,\psi_{(a,A)})$ to $(a,A,\psi_{(a,A)})$. 
    On 1- and 2-morphism level the 2-functor is defined to be the identity.
    It is clear that this defines a 2-functor which is $\Sotwo$-dagger. 
    
    Conversely, if $(a,A,\lambda_{(a,A)})$ is an object in $(S_{\rSpintwo} \circ I_{\rSpintwo} \circ T_{\rSpintwo})(\B,S)$, then by existence of $\lambda_{(a,A)}$ the Frobenius algebra 
    $A$ is symmetric and $\lambda_{(a,A)} \colon A_{(S_A \circ \widetilde\gamma_A^1)^{-1}} = A \to A$ is an invertible intertwiner of $A$-$A$-bimodules. 
    Hence we may send the 
    triple $(a,A,\lambda_{(a,A)})$ to $(a,A,\lambda_{(a,A)})$. 
    The assignment on 1- and 2-morphism level is again trivial.
    This too is an $\Sotwo$-dagger 2-functor and clearly an inverse to the first $\Sotwo$-dagger 2-functor. 
\end{proof}

\subsection[$\Otwo$-idempotent completion]{$\textrm{\textbf{O}}\pmb{(2)}$-idempotent completion}
\label{Subsection: Otwo orbifold completion}

In this section we lift $\Oone$- and $\Otwo$-volutive structures from 2-categories~$\B$ to their idempotent completion, producing $\Oone$- and $\Otwo$-volutive 2-categories $I_{\Oone}\B$ and $I_{\Otwo}\B$, respectively. 
We also work out what conjugating $I_{\Otwo}$ with the operations $S_{\Otwo}, T_{\Otwo}$ of \Cref{Subsection: Otwovolutive 2-categories} does to $\Otwo$-dagger 2-categories, thus producing a candidate for ``$\Otwo$-idempotent completion''. 

\medskip 

We start with an $\Oone$-variant of \Cref{Lemma: rSpintwo volution on idempotent completion}: 

\begin{lemma}
	\label{Lemma: Oone volution on two idempotent completion}
    Let $\B$ be a locally idempotent complete 2-category. 
    Any $\Oone$-volutive structure $(d,\eta,\tau)$ on $\B$ induces an $\Oone$-volutive structure $(d',\eta',\tau')$ on $\Ide \B$.
\end{lemma}

\begin{proof}
    We start by defining the 2-functor $d' \colon \Ide \B \to \Ide \B^{\opp}$. 
    On object level, it assigns $(a,A) \in  \Ide \B$ to the pair $(d(a),d(A))$ where $d(A)$ 
    becomes a $\Delta$-separable Frobenius algebra by applying $d$ to the Frobenius algebra structure on $A$, e.g.\ $d(A) \otimes d(A) \cong d(A \otimes A) \to d(A)$ 
    where we used $d(\mu_A)$ in the second step. 
    On 1-morphism level, we assign to $X \colon (a,A) \to (b,B)$ the 
    1-morphism $d(X) \colon (d(b),d(B)) \to (d(a),d(A))$ where $d(X) \colon d(b) \to d(a)$ is equipped with the $d(A)$-$d(B)$-bimodule structure induced by the one of~$X$ upon applying the 2-functor~$d$. 
    On 2-morphism level, we assign to $f \colon X \to Y$ the map $d(f) \colon d(X) \to d(Y)$. 
    One checks that this assignment is well-defined and functorial due to the functoriality of~$d$. 
    
    Next, we define the invertible 2-transformation $\eta' \colon {d'}^{\opp} \circ d' \Rightarrow \id$. 
    Its 1-morphism component 
    at $(a,A) \in \Ide \B$ is given by the 1-morphism $A \otimes \eta_a \colon (d^{\opp} \circ d)(a) \to a$ equipped with the left $A$-module action 
    induced by $\mu_A$ and the right $d^2(A)$-module action given by 
    \begin{equation}
        \xymatrix{
            A \otimes \eta_a \otimes d^2(A) \ar[rr]^-{\id \otimes \eta_A^{-1}} && A \otimes A \otimes \eta_a \ar[rr]^-{\mu_A \otimes \id} && A \otimes \eta_a .
        }
    \end{equation}
    One deduces from the associativity of $\mu_A$ that $A \otimes \eta_a$ together with these actions becomes an $A$-$d^2(A)$-bimodule. 
    An inverse of this bimodule 
    is given by $\eta_a^{-1} \otimes A$ with actions defined analogously. 
    The 2-morphism component of $\eta'$ at $X \colon (a,A) \to (b,B)$ is defined as 
    \begin{equation}
        \xymatrix{
            X \otimes_A (A \otimes \eta_a) \ar[r]^-{\cong} & X \otimes \eta_a \ar[r]^-{\eta_X} & \eta_b \otimes d^2(X) \ar[r]^-{\cong} & 
            (\eta_b \otimes d^2(B)) \otimes_B d^2(X) \ar[r]^-{\eta_B^{-1} \otimes \id} & (B \otimes \eta_b) \otimes_B d^2(X)
        }
    \end{equation}
    which one checks to be an intertwiner. 
    Moreover, since $\eta$ is an invertible 2-transformation, so is $\eta'$. 
    Lastly, the invertible modification 
    $\tau' \colon ({\eta'}^{\opp})^{-1} \circ \id_{d'} \Rrightarrow \id_{d'} \circ \eta'$ has as its 2-morphism component at $(a,A) \in \Ide \B$ the intertwiner 
    \begin{equation}
        \xymatrix{
            \eta_{d'(a,A)}^{-1} = \eta_{d(a)}^{-1} \otimes d(A) \ar[r]^-{\tau_a} & d(\eta_a) \otimes d(A) \cong d(A \otimes \eta_a) = d'(\eta_{a,A}) .
        }
    \end{equation}
    The coherence condition of $\tau'$ follows from the respective coherence property of $\tau$. 
\end{proof}

\begin{remark} \label{Remark: Oone strictification is compatible with rigid structures}
    In the case that the $\Oone$-volutive structure is induced by a rigid symmetric monoidal structure as in 
    \Cref{Example: symmetric monoidal 2-categories with duals are Oone volutive}, this construction is compatible with the 
    construction of \cite[Section 4.1.4]{75}. 
    More specifically, loc.\ cit.\ describes duals on $\Ide \B$ induced by duals in $\B$, which yields an 
    $\Oone$-volution on $\Ide \B$ induced by the $\Oone$-volution on $\B$ determined by the rigid symmetric monoidal structure on $\B$. 
    The $\Oone$-volution induced 
    on $\Ide \B$ in this way is then recovered by \Cref{Lemma: Oone volution on two idempotent completion}.
    We expect that this extends to an analogue of \Cref{Lemma: Oone volutive idempotent completion is compatible with symmetric monoidal origin}.
\end{remark}

Combining Lemmas~\ref{Lemma: rSpintwo volution on idempotent completion} and~\ref{Lemma: Oone volution on two idempotent completion} allows us to lift $\Otwo$-volutive structures to the idempotent completion: 

\begin{proposition} \label{Proposition: Otwo volutive structures extend to the idempotent completion}
    Let $\B$ be a locally idempotent complete 2-category with right adjoints. 
    Any $\Otwo$-volutive structure $(d,\eta,\tau,S,\Gamma)$ on $\B$ induces an 
    $\Otwo$-volutive structure $(d',\eta',\tau',S',\Gamma')$ on $\Ide \B$.
\end{proposition}

\begin{proof}
    The $\Oone$-volutive structure $(d',\eta',\tau')$ is described in \Cref{Lemma: Oone volution on two idempotent completion} while the $\Sotwo$-volutive 
    structure is described in \Cref{Lemma: rSpintwo volution on idempotent completion}, recalling that we have $\Sotwo = \rSpintwo$ for $r=1$. 
    Next, we describe 
    the invertible 2-transformation $\epsilon' \colon {}^{\star \star}(-) \circ d' \Rightarrow d' \circ (-)^{\star\star}$. 
    Its 1-morphism component at $(a,A) \in \Ide \B$ 
    is given by $d(A^{\vee \vee}) \otimes \epsilon_a$ with $d(A)$-$d(A)$-bimodule actions constructed from the Nakayama map 
    $\widetilde{\gamma}_A \colon A \to A^{\vee \vee}$ (and its left analogue), and the 2-morphism component  $\epsilon_A \colon d(A^{\vee \vee}) \otimes \epsilon_a \to \epsilon_a \otimes {} ^{\vee \vee}d(A)$ of~$\epsilon$ at the 1-morphism $A \colon (a,A) \to (a,A)$. 
    The 2-morphism component of $\epsilon$ at $X \colon (a,A) \to (b,B)$ is given by 
    \begin{equation}
        \xymatrix{
            d(X^{\star \star}) \otimes_{d(A)} (d(A^{\vee \vee}) \otimes \epsilon_a) \ar[r]^-{\cong} & 
            d(X^{\vee \vee}) \otimes_{d(A^{\vee \vee})} d(A^{\vee \vee}) \otimes \epsilon_a  \ar[r]^-{\cong} &
            d(X^{\vee \vee}) \otimes \epsilon_a \ar[r]^-{\epsilon_X} & 
            \epsilon_b \otimes ^{\vee \vee}d(X) \ar[dlll]_-{\cong}  
            \\
            \epsilon_b \otimes ^{\vee \vee}d(B) \otimes_{^{\vee \vee}d(B)} {^{\vee \vee}d(X)} \ar[r]_-{\cong} & 
            (\epsilon_b \otimes ^{\vee \vee}d(B)) \otimes_{d(B)} {^{\star \star}d(X)} \ar[r]_-{\epsilon_B^{-1} \otimes \id} & 
            (d(B^{\vee \vee}) \otimes \epsilon_b) \otimes_{d(B)} {^{\star \star}d(X)} . &
        }
    \end{equation}
    One checks that this defines an invertible 2-transformation using the respective properties of $\epsilon$. 
    Finally, we define the invertible modification 
    $\Gamma' \colon (\id_{d'} \circ (S')^{-1}) \bullet \epsilon' \Rrightarrow S' \circ \id_{d'}$. 
    Its 2-morphism component at $(a,A) \in \Ide \B$ is given by 
    \begin{equation}
        \xymatrix{
            d(S_a^{-1} \otimes A_{(\widetilde\gamma_A^1)^{-1}}) \otimes_{d(A)} (d(A^{\vee \vee}) \otimes \epsilon_a) \ar[r]^-{\cong} &
            d(A)_{(\gamma_{d(A)}^1)^{-1}} \otimes d(S_a^{-1}) \otimes \epsilon_a \ar[r]^-{\id \otimes \Gamma_a} &
            d(A)_{(\gamma_{d(A)}^1)^{-1}} \otimes S_{d(a)}
        }
    \end{equation}
    which one checks to be an invertible modification using the respective properties of $\Gamma$. The coherence condition for $\Gamma'$ follows similarly.
\end{proof}

\begin{remark}
    We expect that the assignment $(\B,d,\eta,\tau,S,\Gamma) \mapsto (\Ide \B,d',\eta',\tau',S',\Gamma')$ extends to a 3-functor 
    \begin{equation}
        \Ide_{\Otwo} \colon \Otwovolbicat \xymatrix{\ar[r]&} \Otwovolbicat.
    \end{equation}
\end{remark}

\begin{remark}\label{Remark: The induced strict Otwo volutive structure on SIT}
Let $(\B,d,\eta,S)$ be an $\Otwo$-dagger 2-category. We want to spell out the $\Otwo$-dagger 2-category 
$(S_{\Otwo} \circ \Ide_{\Otwo} \circ T_{\Otwo})(\B,d,\eta,S)$. The $\Otwo$-volutive 2-category $T_{\Otwo}(\B,d,\eta,S)$ is $(\B,d,\eta,\id,S,\id)$, as explained in the paragraph containing~\eqref{eq:TO2}. 
Next, we describe the $\Otwo$-volutive structure on $(\Ide_{\Otwo} \circ T_{\Otwo})(\B,d,\eta,S)$. 
As we already described the arising $\Sotwo$-volutive 
structure on $\Ide \B$ in \Cref{Subsection: rSpintwo orbifold completion}, it remains to describe the $\Oone$-volutive structure as well as the modification $\Gamma$. 
The 2-functor $d' \colon \Ide \B \to \Ide \B^{\oneop}$ assigns $(a,A) \mapsto (a,d(A))$, $X \mapsto d(X)$, $f \mapsto d(f)$. 
The invertible 2-transformation 
$\eta'$ has 1-morphism component at $(a,A) \in \Ide \B$ given by the $A$-$d^2(A)$-bimodule $A$ where the right action is constructed using $\eta_A$. The 
invertible modification $\eta'$ is virtually trivial. Next, the 1-morphism component of the invertible 2-transformation~$\epsilon'$ at $(a,A) \in \Ide \B$ 
is given by the $d(A)$-$d(A)$-bimodule $d(A^{\vee \vee})$. Finally, the 2-morphism component of the invertible modification $\Gamma'$ at $(a,A) \in \Ide \B$ 
is also virtually trivial. 

We are now in a position to describe the $\Otwo$-volutive 2-category $(S_{\Otwo} \circ \Ide_{\Otwo} \circ T_{\Otwo})(\B,d,\eta,S)$. 
First, the underlying 2-category has objects given by tuples $(a,A,\theta_{(a,A)},\Pi_{(a,A)},\lambda_{(a,A)},\phi_{(a,A)})$ where $(a,A) \in \Ide \B$, $\theta_{(a,A)}
\colon (a,A) \to (a,d(A))$ is a 1-isomorphism in $\Ide \B$, and 
\begin{equation}
	\label{eq:BiGeneralizesInvolution}
	\Pi_{(a,A)} \colon A \otimes_{d^2(A)} d(\theta_{(a,A)})^{-1} \otimes_{d(A)} \theta_{(a,A)} \xymatrix{\ar[r]&} A
\end{equation}
is a 2-isomorphism in $\Ide \B$ satisfying a coherence condition. 
Moreover, $\lambda_{(a,A)} \colon A_{(S_A \widetilde\gamma_A^1)^{-1}} \to A$ is a 2-isomorphism in $\Ide \B$, and $\phi_{(a,A)} \colon 
\theta_{(a,A)}^{-1} \otimes d(A^{\vee \vee}) \otimes ^{\vee \vee} \theta_{(a,A)} \to A$ is a 2-isomorphism satisfying a coherence condition. The 1- and 2-morphisms 
of the underlying 2-category are the same as those of $\Ide \B$. The $\Otwo$-dagger structure on this 2-category is then obtained by carrying out the construction 
described in \Cref{Subsection: Otwovolutive 2-categories}. 
\end{remark}

We explained in \Cref{Subsection: Otwovolutive 2-categories and Oonevolutive 1-categories} that the Hom categories of any $\Otwo$-dagger 2-category admit 
$\Oone$-volutions. Moreover, we were able to construct a new $\Otwo$-dagger 2-category whose Hom categories carry compatible $\Oone$-dagger 
structures (that is, dagger structures) so that the arising 2-category carries a bi-$\dagger$ structure. One may in particular apply this construction to the 
$\Otwo$-dagger 2-category $(S_{\Otwo} \circ \Ide_{\Otwo} \circ T_{\Otwo})(\B,d,\eta,S)$. The objects of the resulting 2-category are the same as those 
of the 2-category underlying $(S_{\Otwo} \circ \Ide_{\Otwo} \circ T_{\Otwo})(\B,d,\eta,S)$, while its 1-morphisms additionally carry 2-isomorphisms enforcing 
compatibility with the $\Otwo$-volutive structure as explained in \Cref{Subsection: Otwovolutive 2-categories and Oonevolutive 1-categories}.

\begin{remark}\label{Remark: simpler objects in the Otwo SIT construction}
    Recall that the 2-category $\Ide \B$ descends from a double category, meaning in particular that any morphism of Frobenius algebras 
    $\phi \colon A \to B$ canonically induces an $A$-$B$-bimodule obtained by twisting the left action of the $B$-$B$-bimodule $B$ along $\phi$. 
    In consequence, there may be objects and 1-morphisms in the 2-category underlying $(S_{\Otwo} \circ \Ide_{\Otwo} \circ T_{\Otwo})(\B,d,\eta,S)$ of 
    a particularly simple form; by that we mean that for instance a 1-isomorphism $\theta_{(a,A)} \colon (a,A) \to (a,d(A))$ in $\Ide \B$, which is by definition an 
    invertible $d(A)$-$A$-bimodule, may come from an invertible Frobenius algebra morphism $A \to d(A)$. 
\end{remark}

\begin{remark}
    Analogously to the oriented case described in \Cref{Definition: Euler completion} and \Cref{Lemma: Euler completion from our perspective}, 
    we propose the composition $S_{\Otwo} \circ T_{\Otwo}$ as a generalization of ``two-dimensional unoriented Euler\footnote{In general, we expect (the TQFT analogue of) 
    this generalization of ``Euler'' completion to add more than the invariants constructed from the Euler characteristic; a more precise statement would require a rigorous 
    definition of unoriented two-dimensional defect TQFTs.} completion''. 
    To see that this is indeed a completion, it would be sufficient to have a result analogous to 
    \Cref{Theorem: The equivalence of strict and non-strict Oone volutive categories with restrictions} in the $\Otwo$-case; if given such a theorem, 
    it would automatically follow that $S_{\Otwo} \circ I_{\Otwo} \circ T_{\Otwo}$ is complete with respect to $S_{\Otwo}\circ T_{\Otwo}$.
    We leave this for future work.
\end{remark}

\section{Examples and applications}
\label{sec:ExamplesApplications}

In this section we study $G$-volutive and $G$-dagger structures on several 2-categories, recovering and generalizing several known results from the literature as well 
as producing new ones. In \Cref{subsec:StateSumModels}, we show that our $G$-idempotent completion construction $S_G\circ I_G\circ T_G$ of Sections~\ref{sec:DaggerStructures} 
and~\ref{sec:AlgebraicDescriptionOrbifolds} for $G \in \{ \Sotwo, \rSpintwo, \Otwo\}$ recovers the algebraic input data and structural results 
for ``$G$-type'' state sum models (\Cref{subsec:StateSumModels}). In the remaining sections, we study $G$-dagger structures on the 2-categories of bundle gerbes 
(\Cref{subsec:BundleGerbes}), Landau--Ginzburg models (\Cref{Subsection: Landau--Ginzburg theory}), and truncated affine Rozansky--Witten models (\Cref{subsec:RWmodels}). 

\subsection{State sum models}
\label{subsec:StateSumModels}

In this section we apply the theory developed in Sections~\ref{sec:DaggerStructures} and~\ref{sec:AlgebraicDescriptionOrbifolds} to the 2-category which is the delooping of finite-dimensional (super) vector spaces. 
The examples of $G$-dagger structures this gives rise to turn out to reproduce the algebraic input data of two-dimensional state sum models for $G$-tangential structure, for $G \in \{ \Sotwo, \rSpintwo, \Otwo\}$. 
This is a non-trivial consistency check of the expectation that our higher idempotent constructions actually give rise to TQFTs. 

\medskip 

Recall the symmetric monoidal rigid category of finite-dimensional vector spaces $\vect$ over some field~$\Bbbk$. 
The delooping $\textrm{B}\vect$ is a symmetric monoidal 2-category with duals and adjoints. 
Its adjoints correspond to the dual objects of $\vect$, while the dual of the unique object $*\in \textrm{B}\vect$ is trivially~$*$. 

We first note that by \Cref{Example: The delooping of vector spaces is strict rSpintwo volutive} and \Cref{Example: The delooping of vector spaces is strict Otwo volutive}, 
$\textrm{B}\vect$ carries both a $\rSpintwo$-dagger structure~$S$ (with 2-morphism components $S_V \colon V^{**} \to V$ the canonical isomorphisms for finite-dimensional vector spaces~$V$) and an $\Otwo$-dagger structure $(d,\eta,S)$. 
We stress that here $d,\eta$ are trivial, hence below we write $(\id,\id,S)$ for $(d,\eta,S)$.
Secondly, by \Cref{Proposition: Otwo volutive structures extend to the idempotent completion} and \Cref{Lemma: rSpintwo volution on idempotent completion}, the idempotent completion 
$\Ide \textrm{B}\vect$ then also carries an $\Otwo$-volutive structure and a 
$\rSpintwo$-volutive structure, see the respective proofs for explicit details. 
Alternatively, this can be seen by noting that $\Ide \textrm{B}\vect$ too is symmetric monoidal and has duals and adjoints. 
Hence we obtain the $\rSpintwo$-dagger 2-category 
\begin{equation}
	\big(S_{\rSpintwo} \circ \Ide_{\rSpintwo} \circ T_{\rSpintwo}\big)(\textrm{B}\vect,S)
\end{equation}
and the $\Otwo$-dagger 2-category 
\begin{equation}
	\big(S_{\Otwo} \circ \Ide_{\Otwo} \circ T_{\Otwo}\big)(\textrm{B}\vect,\id,\id,S) . 
\end{equation}
We will subsequently compare both of these to the literature, starting with the former.

\subsubsection*{The case of $\textrm{\textbf{SO}}\pmb{(2)}$}

We first consider the special case of $r=1$, that is, the case of an $\Sotwo$-dagger structure on $\textrm{B}\vect$. 
Thanks to \Cref{Theorem: comparison of Euler orbifold completion and SIT}, we have 
$(S_{\Sotwo} \circ \Ide_{\Sotwo} \circ T_{\Sotwo})(\textrm{B}\vect,S) \cong ((\textrm{B}\vect,S)_{\orb})_{\eu}$. 
On the other hand, the right-hand side is equivalent 
to the 2-category $\textrm{SepSymFrobAlg}_\Bbbk$ of {separable} symmetric Frobenius algebras, cf.\ \cite[Section 5.1.1]{carqueville2023orbifold}. 
Thus we find
\begin{equation}
	\big(S_{\Sotwo} \circ \Ide_{\Sotwo} \circ T_{\Sotwo}\big)(\textrm{B}\vect,S)
	\;\cong\;
	\textrm{SepSymFrobAlg}_\Bbbk . 
\end{equation}
This means that our construction specializes to precisely reproduce the symmetric monoidal pivotal 2-category which governs two-dimensional oriented state sum models with defects, see \cite{Fukuma1994} for the closed case, 
\cite{LAUDA2007} for the open-closed case, and \cite{carqueville2023orbifoldstopologicalquantumfield} for the general defect case.

\subsubsection*{The case of $\textrm{\textbf{Spin}}\pmb{(2)^r}$}

We turn to the case of $\rSpintwo$ with $r \geqslant 2$. As explained in \Cref{Subsection: rSpintwo orbifold completion}, the objects in $(S_{\Sotwo} \circ \Ide_{\Sotwo} \circ 
T_{\Sotwo})(\textrm{B}\vect,S)$ are separable Frobenius $\Bbbk$-algebras~$A$ for which the Nakayama$^r$ automorphism is trivial. More generally, one may consider any 
symmetric monoidal rigid category $\C$ rather than $\vect$ and carry out the same construction; of particular interest here is $\C = \mathrm{s}\hspace{-0.05cm}\vect$, the category of super 
vector spaces. 
Such algebras are closely related to the ones of \cite{85,Novak2015} describing  two-dimensional $r$-spin state sum models. 
The main difference to our approach is that our theory is set in the context 
of $\rSpintwo$-dagger 2-categories, generalizing the well-understood pivotal ($r$=1) case, whereas the authors of \cite{85,Novak2015} study separable Frobenius algebras in 
symmetric monoidal categories and ask the honest $r$-th power of the Nakayama automorphism to be trivial; see \cite[Remark 4.2]{75} for a comparison of different conventions of Nakayama automorphisms. 
Here we recall from \Cref{Lemma: The r-th power of the Nakayama automorphism in the Sotwo case} that in case the $\rSpintwo$-dagger structure is induced from an $\Sotwo$-dagger structure (which in turn is induced from the symmetric monoidal structure on the category), the $r$-th power of the Nakayama automorphism coincides with the Nakayama$^r$ automorphism. 
Hence our construction precisely reproduces the algebraic input data of two-dimensional $r$-spin state sum models, as expected.

\subsubsection*{The case of $\O\textrm{\textbf{O}}\pmb{(2)}$}

Finally, we consider the case of $\Otwo$. Considering the special case of $(\textrm{B}\vect,\id,\id,S)$ in \Cref{Remark: The induced strict Otwo volutive structure on SIT},
we find that the objects in $(S_{\Otwo} \circ \Ide_{\Otwo} \circ T_{\Otwo})(\textrm{B}\vect,\id,\id,S)$ are essentially separable symmetric Frobenius algebras~$A$ that additionally carry a (generalized) involution. 
We already explained why these Frobenius algebras are separable and symmetric in the $\Sotwo$-case, so it remains to discuss the involution on $A$. First, we recall from 
\Cref{Lemma: Oone volution on two idempotent completion} that the $\Oone$-volutive structure $(\id,\id)$ on $\textrm{B}\vect$ induces an $\Oone$-volutive structure on $\Ide\textrm{B}\vect$. The 
underlying functor $d \colon \Ide\textrm{B}\vect \to \Ide\textrm{B}\vect^{\oneop}$ of this $\Oone$-volutive structure assigns an object $(\star,A) \in \Ide\textrm{B}\vect$ to the object 
$(\star,A^{\opp}) \in \Ide\textrm{B}\vect^{\oneop}$ where $A^{\opp}$ denotes the \emph{opposite algebra} of $A$; this can be seen either by spelling out the construction in 
\Cref{Lemma: Oone volution on two idempotent completion} for the case at hand or equivalently by combining \Cref{Remark: Oone strictification is compatible with rigid structures} 
with the explicit discussion of duals in the idempotent completion (there denoted by $(-)_{\operatorname{eq}})$ of an arbitrary rigid monoidal 2-category in \cite[Section 4.1.4]{75} and 
then again specializing to $\textrm{B}\vect$. We note that the underlying 2-transformation $\eta \colon d^{\oneop} \circ d \to \id$ of the $\Oone$-volutive structure on $\Ide\textrm{B}\vect$ is trivial, since 
$(A^{\opp})^{\opp} = A$. 

A particular special case of the data $\theta_{(a,A)},\Pi_{(a,A)}$ described in \Cref{Remark: The induced strict Otwo volutive structure on SIT} 
for the case at hand arises when $\theta_{(a,A)}$ is induced by an involutive isomorphism of Frobenius algebras $\theta \colon A \to A^{\opp}$ (see \Cref{Remark: simpler objects in the Otwo SIT construction}) such that 
$\theta^2 = \id_A$ and $\Pi_{(a,A)}$ is trivial, which is why we speak of \emph{involutions}. 
These algebraic structures are closely related to the stellar algebras of \cite{26}, and to the ones of \cite{tavares2014spinstatesummodels,Barrett2013twodimensionalstatesummodels} following \cite{Karimipour1997}, describing two-dimensional unoriented state sum models. 
There are minor differences between our algebraic structures and those of loc.\ cit., however below we shall explain that the latter are equivalent to a special case of ours. 

One difference between our construction and that of \cite{Barrett2013twodimensionalstatesummodels} is that the latter considers special Frobenius algebras. 
Such algebras are always separable and hence included in our setting. 
Moreover the involutions $A \to A^{\opp}$ of loc.\ cit.\ are algebra isomorphisms rather than bimodules, which are included in our setting via \Cref{Remark: simpler objects in the Otwo SIT construction} and the discussion above.
We expect that our slightly more general definitions also carry physical meaning in the sense that all objects of $(S_{\Otwo} \circ \Ide_{\Otwo} \circ T_{\Otwo})(\textrm{B}\vect,\id,\id,S)$ can be used as the algebraic input data for unoriented closed TQFTs of state sum type. 
In fact it is natural to expect that the whole of $(S_{\Otwo} \circ \Ide_{\Otwo} \circ 
T_{\Otwo})(\textrm{B}\vect,\id,\id,S)$ can be used as the algebraic input data of an unoriented \textsl{defect} TQFT, analogously to the well-studied oriented case discussed above. 
In particular, our discussion of the $\Oone$-volutions on Hom categories in \Cref{Subsection: Otwo orbifold completion} should yield the correct algebraic notion of unoriented line defects.

\subsection{Bundle gerbes}
\label{subsec:BundleGerbes}

In this section we discuss volutive structures on the 2-category of bundle gerbes over a fixed manifold. 

\medskip 

Let $M$ be a smooth manifold. 
We consider the 2-category of bundle gerbes $\Grb(M)$ over $M$ described in \cite{waldorf2007morphismsbundlegerbes}. 
Here we recall only some of the relevant definitions. 
Given surjective submersions $\pi \colon Y \to M$ and $\pi' \colon Y' \to M$, we denote their fiber product by 
$Y {}_{\pi}\!\times_{\pi'} Y' := \{(y,y') \,|\, \pi(y) = \pi'(y') \}$. 
We denote the $k$-fold fiber product of $Y$ with itself by 
\begin{equation}
	Y^{[k]} := Y {}_{\pi}\!\times_\pi Y {}_{\pi}\!\times_\pi \cdots {}_{\pi}\!\times_\pi Y \,.
\end{equation}

\begin{definition}
    A \emph{bundle gerbe} $\mathscr{G}$ over $M$ consists of a surjective submersion $\pi: Y \to M$, a line bundle $L$ over $Y^{[2]}$ and an associative
    isomorphism $\mu : \pi_{23}^* L \otimes \pi_{12}^* L \to \pi_{13}^* L$ of line bundles over $Y^{[3]}$. Here associativity amounts to a coherence condition 
    over $Y^{[4]}$.
\end{definition}

\begin{example}
    The \emph{trivial} bundle gerbe $\mathscr{I}$ over $M$ consists of the surjective submersion $\id \colon M \to M$, the line bundle $M \times \mathbb{C}$, and 
    the isomorphism $\mu$ induced by multiplication of complex numbers. 
\end{example}

\begin{remark}
    Ignoring the surjective submersion $\pi$, it is helpful to think of bundle gerbes as structures resembling algebras, which suggests natural notions of 1- and 2-morphisms 
    as bimodule and intertwiners, respectively. 
    This leads to the construction of a 2-category $\Grb(M)$, see  \cite{waldorf2007morphismsbundlegerbes} for a careful treatment of the involved structures. 
\end{remark}

The 2-category $\Grb(M)$ can be endowed with a symmetric monoidal structure such that it is a categorification of the symmetric monoidal category $\Vectbdl(M)$ of vector bundles over~$M$, i.e.\ $\Grb(M)(\mathscr{I},\mathscr{I}) \cong \Vectbdl(M)$. 
Moreover, one shows that $\Grb(M)$ admits duals and adjoints. In particular, the dual of a bundle gerbe $\mathscr{G}=(Y,\pi,L,\mu)$ is the bundle gerbe $(Y,\pi,L^*,(\mu^*)^{-1})$ where $L^*$ denotes the dual line bundle.

Again we refer to \cite{waldorf2007morphismsbundlegerbes} and \cite{kristel20222vectorbundles} for details. 
By Examples~\ref{Example: symmetric monoidal 2-categories with duals are Otwo volutive} and~\ref{Example: Symmetric monoidal 2-categories with duals and adjoints are rSpintwo volutive}, 
the 2-category $\Grb(M)$ hence admits an $\Otwo$-volutive structure as well as a $\rSpintwo$-volutive structure for any $r \in \Z_{\geqslant 0}$. 
We saw that the strictification procedures $S_{\Otwo}$ and $S_{\rSpintwo}$ essentially amount to finding certain trivializations of structures internal to $\Grb(M)$, 
such as the duals and (higher) powers of the Serre automorphism. 
In the following we want to determine which objects admits such trivializations, that is, we wish to determine 
$S_{\Otwo}\Grb(M)$ and $S_{\rSpintwo}\Grb(M)$, starting with the latter. 

\medskip 

Let $\mathscr{G}=(Y,\pi,L,\mu)$ be a bundle gerbe over $M$. The evaluation 1-morphism $\ev \colon \mathscr{G} \otimes \mathscr{G}^* \to \mathscr{I}$ is given by the 
trivial line bundle together with the line bundle isomorphism which fiberwise reads $L_{(y_1,y_2)} \otimes L_{(y_1,y_2)}^* \to \mathbb{C}$. In other words, this 1-isomorphism 
is fiberwise induced by the ordinary evaluation map for vector spaces. Here we have tacitly used the canonical refinement along the diagonal map $\Delta \colon Y \to Y^{[2]}$.
Since $\ev$ is invertible, we may take its inverse to be its adjoint. The braiding $\beta_{\mathscr{G},\mathscr{G}} \colon \mathscr{G} \otimes \mathscr{G} \to 
\mathscr{G} \otimes \mathscr{G}$ consists again of the trivial line bundle together with the isomorphism constructed fiberwise from the braiding in the category 
of vector spaces, that is, $\beta_{ L_{(y_1,y_2)},L_{(y_1,y_2)}} \colon L_{(y_1,y_2)} \otimes L_{(y_1,y_2)} \cong L_{(y_1,y_2)} \otimes L_{(y_1,y_2)}$. 
Combining these structures, we find that the Serre automorphism of $\mathscr{G}$ consists of the line bundle~$L$ over $Y^{[2]}$, and the isomorphism of line bundles 
\begin{equation}
    \xymatrix{
        L \ar[rr]^-{\id \otimes \ev_L^{\vee}} && L \otimes L \otimes L^*  \ar[rr]^-{\beta_{L,L} \otimes \id} && L \otimes L \otimes L^* \ar[rr]^-{\id \otimes \ev_L} && L
    }
\end{equation}
which is straightforwardly checked to be trivial too. 
Indeed, working fiberwise over a chosen point, the amounts to the standard computation
\begin{align*}
    v \longmapsto v \otimes \left(\sum_i b_i \otimes b_i^* \right) 
      \longmapsto \sum_i b_i \otimes v \otimes b_i^*
      = \sum_i b_i \otimes (v_i \cdot b_i) \otimes b_i^*
      = \sum_i (v_i \cdot b_i) \otimes b_i \otimes b_i^*
      \longmapsto \sum_i (v_i \cdot b_i) 
      = v
\end{align*}
in $\operatorname{vect}_{\mathbb{C}}$, where $b_i$ is a basis of the fiber of $L$ over the chosen point. 
Hence: 

\begin{lemma}
	\label{lem:SerreTrivialGrbM}
    The Serre automorphism of any object $\mathscr{G} \in \Grb(M)$ is trivializable. 
\end{lemma}
Recall from \cite[Theorem 2.5.4]{Waldorf2008phd} that the automorphism 2-group of every bundle gerbe $\mathscr{G} \in \Grb(M)$ is equivalent to the groupoid 
$\Linebdl(M)^\times$ of line bundles over $M$. In particular, trivializations of the Serre automorphism of $\mathscr{G}$ are the same as automorphisms of the trivial line bundle 
over $M$; in other words, these are smooth functions $M \to \mathbb{C}^\times$.

\begin{corollary}\label{Corollary: Bundle gerbes define extended spin TQFTs}
    For $r \in \Z_{\geqslant 0}$, any object $\mathscr{G} \in \Grb(M)$ gives rise to a fully extended $r$-spin TQFT 
    \begin{equation}
        \mathscr{Z}_{\mathscr{G}} \colon \Bord_{2,1,0}^{\rSpintwo} \xymatrix{\ar[r]&} \Grb(M) \, .
    \end{equation} 
\end{corollary}

\begin{proof}
    This follows directly from the two-dimensional $r$-spin cobordism hypothesis proved in \cite{84,82,75} and \Cref{lem:SerreTrivialGrbM}.
\end{proof}
\begin{remark}
    Every fully extended $r$-spin TQFT for $r \geqslant 2$ described in \Cref{Corollary: Bundle gerbes define extended spin TQFTs} is induced by a fully extended oriented 
    TQFT, hence, cannot distinguish different spin structures. We expect that, in passing from ordinary bundle gerbes to super bundle gerbes, one obtains 
    honest spin theories. 
\end{remark}

\begin{remark}
	\label{rem:GrbMSpinVolutive}
    Having proved that any object in $\Grb(M)$ has trivializable Serre automorphism, we find that $\Grb(M)$ embeds into the 2-category 
    underlying the $\rSpintwo$-dagger 2-category $S_{\rSpintwo}\Grb(M)$. In other words, $\Grb(M)$ admits a $\rSpintwo$-dagger structure. 
\end{remark}

Next, we wish to determine the self-dual objects in $\Grb(M)$, that is, those objects $\mathscr{G} \in \Grb(M)$ for which $\mathscr{G} \cong \mathscr{G}^*$ holds. 
Recall that the dual bundle gerbe is automatically an inverse; this follows essentially from the classification of bundle gerbes (up to isomorphism) in terms of 
their Dixmier--Douady class $\DD(\mathscr{G}) \in H^3(M,\mathbb{Z})$, see \cite[Section 4.3]{murray2008introductionbundlegerbes}. We conclude the following.

\begin{lemma}
    Every object in $S_{\Otwo}\Grb(M)$ has an underlying bundle gerbe $\mathscr{G}$ which is 2-torsion, i.e.\ $\mathscr{G} \otimes \mathscr{G} \cong \mathscr{I}$. 
\end{lemma}
Determining which bundle gerbes are coherently self-dual in the sense of \Cref{Remark on coherently self-dual objects} is more difficult; we will leave this for future research.
As explained in \Cref{Subsection: Otwovolutive 2-categories and Oonevolutive 1-categories}, the Hom categories in any $\Otwo$-dagger 2-category 
carry $\Oone$-volutions. 
This holds in particular for $S_{\Otwo}\Grb(M)$, whose Hom categories are those of $\Grb(M)$. Recall further that there is a fully faithful functor 
$\mathrm{B}\Vectbdl(M) \to \Grb(M)$ which sends the trivial object of the delooping of the monoidal category of vector bundles over $M$ to the trivial bundle gerbe over $M$. 
This functor extends to a fully faithful functor $\mathrm{B}\Vectbdl(M) \to S_{\Otwo}\Grb(M)$, so that $\mathrm{B}\Vectbdl(M)$ inherits an $\Oone$-volutive structure 
from the $\Otwo$-dagger structure on $S_{\Otwo}\Grb(M)$. The discussion in \Cref{Example: comparing Oone volutions on deloopings} implies the following. 
\begin{proposition}
    The $\Oone$-volution on $\Vectbdl(M)$ induced by the $\Otwo$-dagger structure on $S_{\Otwo}\Grb(M)$ according to the prescription in 
    \Cref{Subsection: Otwovolutive 2-categories and Oonevolutive 1-categories} coincides with the $\Oone$-volution induced by 
    \Cref{Example: symmetric monoidal categories with duals are Oone volutive}. 
\end{proposition}

\begin{remark}
    If one considers the $\Oone$-volution on $\Vectbdl(M)$ obtained by combining the one of \Cref{Example: symmetric monoidal categories with duals are Oone volutive} 
    with the complex conjugation functor along the same lines as \Cref{Example: complex vector spaces with complex conjugation and dualization}, one can 
    apply $S_{\Oone}$ to this $\Oone$-volutive category to obtain the dagger category of hermitian vector bundles. 
    This is analogous to the dagger structure on finite-dimensional hermitian vector spaces, cf.\ Example~\ref{exa:HermVect}. 
\end{remark}

\subsection{Landau--Ginzburg models}
\label{Subsection: Landau--Ginzburg theory}

In this section we consider the 2-category $\LG$ of Landau--Ginzburg models. 
We show that the standard pivotal structure of a certain full sub-2-category can be recovered from the general theory of Section~\ref{sec:DaggerStructures} by constructing the $\Sotwo$-dagger structure coming from the fully dualizable symmetric monoidal structure of $\LG$. 
Analogously working out the $\Oone$-volutive structure connects to the physics literature of orientifold branes. 

\medskip 

We start by sketching the symmetric monoidal structure with duals and adjoints for the 2-category $\LG$ of Landau--Ginzburg models, see \cite{3,CarquevilleMurfet2016} for details. 
Objects of $\LG$ are pairs $(\Bbbk[x_1,\dots,x_n],W)$ 
where $n \in \mathbb{Z}_{\geqslant 0}$ and $W \in \Bbbk[x_1,\dots,x_n]$ is a polynomial such that the Jacobi algebra
$
    \Bbbk[x_1,\dots,x_n]/(\partial_{x_1}W,\dots,\partial_{x_n}W)
$
is finite-dimensional over $\Bbbk$. 
We often abbreviate $(\Bbbk[x_1,\dots,x_n],W)$ to $(\Bbbk[{\xu}],W)$ or just~$W$. 
Given two objects 
$(\Bbbk[{\xu}],W)$, $(\Bbbk[{\zu}],V)$, the Hom category is 
\begin{equation*} \label{Definition: Hom categories of Landau-Ginzburg}
    \LG\big((\Bbbk[{\xu}],W),(\Bbbk[{\zu}],V)\big) = \Ide \hmf(\Bbbk[{\xu},{\zu}], V - W)
\end{equation*}
where the right-hand side is the idempotent completion of the homotopy category 
of finite-rank matrix factorizations $(X,d_X)$ of $V-W$ over $\Bbbk[{\xu},{\zu}]$. 
This means that $X=X^0 \oplus X^1$ is a $\mathbb{Z}$-graded free finite-rank $\Bbbk[{\xu},{\zu}]$-module together with an odd endomorphism~$d_X$ such that $d_X^2 = (V-W)\cdot \id$, while morphisms $(X,d_X)\to(Y,d_Y)$ in $\hmf(\Bbbk[{\xu},{\zu}], V - W)$ are even cohomology classes of the differential on $\operatorname{Hom}_{\Bbbk[\xu,\zu]}(X,Y)$ given by $\Phi \mapsto d_Y\Phi - (-1)^{|\Phi|}\Phi d_X$. 
The identity 1-morphism on $(\Bbbk[{\xu}],W)$ is denoted $(I_W,d_{I_W})$ and given by a deformation of the Koszul complex associated to the differences $x_i-x'_i \in \Bbbk[\xu,\xu']$ for $i\leqslant n$.

The monoidal structure on $\LG$ is given on objects by  
\begin{equation*}
    (\Bbbk[{\xu}],W) \otimes (\Bbbk[{\zu}],V) := (\Bbbk[{\xu},{\zu}], W + V)
\end{equation*}
while on morphisms it is the tensor product over~$\Bbbk$. 
The monoidal unit is $\mathbb{I} = (\Bbbk,0)$, and the symmetric braiding is essentially given by identity and coherence morphisms. 
The dual of an object $(\Bbbk[\xu],W)$ is $(\Bbbk[\xu],W)^* := (\Bbbk[x_1,\dots,x_n],-W)$. 
To describe adjoints of 1-morphisms $X (\Bbbk[x_1,\dots,x_n],W) \to (\Bbbk[z_1,\dots,z_m],V)$, we first consider $X^{\#} := \operatorname{Hom}_{\Bbbk[{\xu},{\zu}]}(X,\Bbbk[{\xu},{\zu}])$ together with 
twisted differential given by $d_{X^{\#}}(\varphi) = (-1)^{|\varphi| +1} \phi \circ d_X$ for homogenous $\varphi \in X^{\#}$. 
The right and left adjoints then are $X^\vee = X^{\#}[n]$ and $^{\vee} X = X^{\#}[m]$ together with their shifted twisted differentials, where~$n$ and~$m$ are the numbers of variables on which the source and target polynomials of~$X$ depend. 
In particular, every object in $\LG$ is fully dualizible. 

\medskip 

Since $\LG$ is symmetric monoidal with all duals and adjoints, Examples~\ref{Example: symmetric monoidal 2-categories with duals are Otwo volutive} and~\ref{Example: Symmetric monoidal 2-categories with duals and adjoints are rSpintwo volutive} imply that $\LG$ carries an $\Otwo$-volutive structure and a $\rSpintwo$-volutive structure for any $r \in \Z_{\geqslant 0}$. 
We then obtain an $\Otwo$-dagger 2-category $S_{\Otwo}\LG$ and a $\rSpintwo$-dagger 2-category $S_{\rSpintwo}\LG$ for any $r \in \Z_{\geqslant 0}$. In the rest of this section, we compare these structures to the current literature on higher dagger structures in the context of Landau--Ginzburg models, starting with the $\rSpintwo$-case. 

First, it was determined in \cite[Lemma 3.8]{3} that the Serre automorphism of an object satisfies $S_{(\Bbbk[x_1,\dotsc, x_n],W)} \cong I_W[n]$; from this fact one deduces that the Serre automorphism is trivializable if and only if $n$ is even. 
In particular, the full sub-2-category $\LG^{\operatorname{even}} \subset \LG$
whose objects depend on an even number of variables embeds into $S_{\Sotwo}\LG$ and hence admits an $\Sotwo$-dagger structure, that is, a pivotal structure. 
The $\Sotwo$-dagger structure on $\LG^{\operatorname{even}}$ depends on the embedding of $\LG^{\operatorname{even}}$ into $S_{\Sotwo}\LG$, that is, on the choice of trivialization of the Serre automorphism 
for a given object. Noticing that the Serre automorphism is really trivial in the present case, we may choose the trivial trivialization for each object, which yields 
the following comparison result, whose proof consists of a translation between different conventions for $\Sotwo$-dagger structures and will be omitted here. 

\begin{proposition} \label{Proposition: comparison of pivotal structures on LG}
    The $\Sotwo$-dagger structure induced on $\LG^{\operatorname{even}}$ as explained in the previous paragraph coincides with the one described in 
    \cite[Section 7]{CarquevilleMurfet2016}. 
\end{proposition}
We remark that the 1-morphism $\Omega_W$ in \cite[Section 7]{CarquevilleMurfet2016}. is in fact the Serre automorphisms of $(\Bbbk[\xu],W)$. 
Regarding the case of $r = 2$, the 2-isomorphism $\mu_W \colon \Omega_W \otimes \Omega_W \to I_W$ of loc.\ cit.\ is then a trivialization of the square of the Serre automorphism, which always exists. 
In other words, the 2-category $\LG$ embeds into $S_{\rSpintwo}\LG$ for $r=2$, hence 
admits a $\rSpintwo$-dagger structure for $r=2$. 
For general $r \in \Z_{\geqslant 1}$, we deduce from the periodicity of the shift functor, $[2] = \id$, that $\LG$ admits a $\rSpintwo$-dagger structure if and only if $r$ is even and $\LG^{\operatorname{even}}$ admits a $\rSpintwo$-dagger structure for any~$r$. 
Again, these structures depend on the respective chosen embeddings into $S_{\rSpintwo}\LG$.

\begin{remark}
    We stress that the technical difficulties of determining the $\Sotwo$-dagger structure in \cite[Section 7]{CarquevilleMurfet2016} are circumvented by our method.
    Indeed, we deduce the existence of the $\Sotwo$-dagger structure from the symmetric monoidal structure and the explicit description of duals and adjoints in loc. 
    cit. and \cite{3}. Moreover, without any additional work, we obtain the respective results for $r \geqslant 2$, which have not been described so far.
\end{remark}

\medskip 

In the remainder of this section, we study the $\Otwo$-dagger 2-category $S_{\Otwo}\LG$ as well as the induced $\Oone$-volutions on its Hom categories, applying the general discussion \Cref{Subsection: Otwovolutive 2-categories and Oonevolutive 1-categories}. It turns out that the latter are closely related to the constructions of 
\cite{Hori2008}. 

Determining the objects in $S_{\Otwo}\LG$ is a difficult problem in general, due to the non-semisimplicity of the Hom categories of $\LG$. 
There are however special objects of a simpler form, which also naturally arise in the study of Landau--Ginzburg orientifolds in \cite{Hori2008}: 

\begin{lemma}\label{Lemma: Hori Walcher objects in our setting}
    Let $(\Bbbk[{\xu}],W) \in \LG$ and let $\tau \colon \Bbbk[{\xu}] \to \Bbbk[{\xu}]$ be a ring isomorphism such that $\tau(W) = -W$ and $\tau^2 = \id$. 
    Then $(\Bbbk[{\xu}],W)$ is self-dual, i.e.\ $(\Bbbk[{\xu}],W) \cong (\Bbbk[{\xu}],W)^*$. 
    Moreover, $(\Bbbk[{\xu}],W)$ is coherently self-dual in the sense that the tuple $((\Bbbk[{\xu}],W),\tau)$ can be completed to an object $S_{\Oone}\LG$, see 
    \Cref{Remark on coherently self-dual objects}.
\end{lemma}
\begin{proof}
	We work with the fibrant double category $((\LG)_0,(\LG)_1)$ whose horizontal 2-category is $\LG$, see \cite{McNamee2009} and \cite{25} for details. 
    We recall that morphisms $(\Bbbk[{\xu}],W) \to (\Bbbk[{\zu}],V)$ in $(\LG)_0$ are ring isomorphisms $f \colon \Bbbk[{\xu}] \to \Bbbk[{\zu}]$ satisfying $f(W) = V$. 
    Since the double category $((\LG)_0,(\LG)_1)$ has companions, any morphism $f \colon (\Bbbk[{\xu}],W) \to (\Bbbk[{\zu}],V)$ in $(\LG)_0$ induces a 1-morphism 
    $(\Bbbk[{\xu}],W) \to (\Bbbk[{\zu}],V)$ in $\LG$ by twisting the identity 1-morphism $(I_V,d_{I_V})$ by $f$. In the situation at hand, $\tau \colon (\Bbbk[{\xu}],W) \to (\Bbbk[{\xu}],-W)$ is an isomorphism in $(\LG)_0$ and hence 
    induces a 1-isomorphism $(\Bbbk[{\xu}],W) \to (\Bbbk[{\xu}],-W)$ in $\LG$, namely the companion of $\tau$ given by the fibrant structure of the double category 
    $((\LG)_0,(\LG)_1)$. 
    The first claim follows directly from $(\Bbbk[{\xu}],W)^* = (\Bbbk[{\xu}],-W)$. 
    The second claim follows from $\tau^2 = \id$ by noting that the 1-isomorphism $(\Bbbk[{\xu}],W)^{**} \cong (\Bbbk[{\xu}],W)$, which is part of the data of the $\Oone$-volutive structure on $\LG$, is trivial. 
\end{proof}

\begin{remark}
    In general, not every object in $S_{\Oone}\LG$ has to be of the form described in \Cref{Lemma: Hori Walcher objects in our setting}. The restriction we made is that the 
    structure expressing the coherent self-duality of an object is induced by a morphism in $(\LG)_0$, which is not the most general case. However, determining whether an 
    arbitrary object in $\LG$ admits coherent self-duality data in the general sense is a difficult problem. 
\end{remark}

Combining the results of \Cref{Lemma: Hori Walcher objects in our setting} and the discussion in the $\Sotwo$-case, we obtain the following. 

\begin{theorem}
    The full sub-2-category $\widehat{\LG} \subseteq \LG$ consisting of objects $(\Bbbk[x_1,\dotsc,x_{2n}],W)$ depending on an even number of variables and admitting an isomorphism~$\tau$ as in \Cref{Lemma: Hori Walcher objects in our setting} embeds into $S_{\Otwo}\LG$, hence admits an $\Otwo$-dagger structure.
\end{theorem}

By the discussion in \Cref{Subsection: Otwovolutive 2-categories and Oonevolutive 1-categories}, the Hom categories of $S_{\Otwo}\LG$ carry $\Oone$-volutions. 
Given an object $(\Bbbk[{\xu}],W,\tau) \in S_{\Otwo}\LG$ (here we have slightly shortended the notation by omitting some of the trivialization data, 
cf. \Cref{Subsection: Otwovolutive 2-categories}), we consider the category $\MF(W) := \Hom_{S_{\Otwo}\LG}(\mathbb{I},(\Bbbk[{\xu}],W,\tau))$ for which we want to describe the $\Oone$-volutive structure explicitly. Recalling our discussion of duals and adjoints in $\LG$, we find the following (under some tacit identifications). First, the functor 
$\hat{d} \colon \MF(W) \to \MF(W)^{\opp}$ assigns a matrix factorization $(X,d_X) \in \MF(W)$ to the matrix factorization $(X,d_X)_{\tau}^\vee$ where $(-)_{\tau}$ 
denotes the twisting by $\tau$ (which amounts to changing the module action of $X$ by composition with $\tau$ as explained in \cite{McNamee2009}), and $(-)^\vee$ denotes the 
right adjoint (which amounts to considering the module $X^{\#} = \operatorname{Hom}_{\Bbbk[{\xu}]}(X,\Bbbk[{\xu}])$ with the corresponding differential $d_{X^{\#}}$ as explained 
above). Recalling that $\tau^2= \id$, the component of the natural isomorphism $\hat{\eta} \colon \hat{d}^{\opp} \circ \hat{d} \to \id_{\MF(W)}$ at $(X,d_X)$ is given by 
the canonical isomorphism $(X^{\# \#},d_{X^{\# \#}}) \cong (X,d_X)$. 

The $\Oone$-volutive structure on $\MF(W)$ we described in the previous paragraph coincides (up to conventions, e.g.\ application of the shift functor [1] described in \cite{3}) 
with the \emph{parity functor} $\mathcal{P}$ described in \cite[Equation (3.6)]{Hori2008}. Moreover, the dagger category $S_{\Oone}(\MF(W),\hat{d},\hat{\eta})$ is 
(again, up to conventions) the category $\MF^{\epsilon=1}_{\mathcal{P}}(W)$ given in \cite[Section 4.1]{Hori2008} with dagger structure $P$ described in \cite[Equation (4.6)]{Hori2008} for 
their case\footnote{The case of $\epsilon =-1$ considered in \cite{Hori2008} is not covered by our theory; this would require a version of our theory in the setting of linear categories.
In \cite[Section 4.1]{Hori2008}, $\epsilon$ is an additional parameter appearing in a strification procedure similar to the one we presented as $S_{\Oone}$.}  of $\epsilon = +1$, i.e.\  
\begin{equation}
    S_{\Oone}(\MF(W),\hat{d},\hat{\eta}) \cong \MF^{\epsilon=1}_{\mathcal{P}}(W).
\end{equation}
Finally, we note that the \emph{antibrane functor} of \cite[Section 3.2]{Hori2008} is (up to conventions) the shift functor [1] and objects in $S_{\Oone}(\MF(W),\hat{d},\hat{\eta})$ 
are the \emph{invariant branes} of \cite{Hori2008}.

We summarize this comparison by remarking that the dagger categories of invariant branes constructed in \cite{Hori2008} are essentially the $\Oone$-strictified Hom categories 
of the $\Otwo$-dagger 2-category $S_{\Otwo}\LG$ extracted from the symmetric monoidal rigid 2-category $\LG$. As such, our discussion 
reveals the 2-categorical origin of the constructions in \cite{Hori2008}. 

\begin{remark}
    Given the $\Otwo$-dagger 2-category $S_{\Otwo}\LG$, it would be interesting to compare the discussion of Landau--Ginzburg orientifolds in \cite[Sections 3.3 \& 4.2]{Hori2008}
    to the structures contained in $(S_{\Otwo} I_{\Otwo} T_{\Otwo}) (S_{\Otwo}\LG)$. In the framed/oriented-case, such a comparison has been carried out in 
    \cite[Section 7]{18}.
\end{remark}

\subsection{Truncated affine Rozansky--Witten models}
\label{subsec:RWmodels}

Recall from \cite{Kapustin_2009, kapustin2009threedimensionaltopologicalfieldtheory, 1,Brunner2023TruncatedAR} that there is a\footnote{More precisely, for simplicity here we consider the 2-category of loc.\ cit.\ ``without R-charges'', i.e.\ without extra $\mathbb{Q}$-gradings.} symmetric monoidal 2-category $\RW$ with duals and adjoints which can be thought of as the homotopy 2-category of a 3-category of affine Rozansky--Witten models, see \cite{riva2024highercategoriespushpullspans}.
Roughly, objects of this 2-category are (possibly empty) lists of variables $\xu = (x_1,...,x_n)$ corresponding to Rozansky--Witten models with target space $T^*\mathbb{C}^n$ for some $n\in\Z_{\geqslant 0}$, 1-morphisms ${\xu}\to{\yu}$ are polynomials $W \in \mathbb{C}[{\xu},{\yu},\zu]$ where~$\zu$ is an extra list of variables, and 2-morphisms $W \to V$ are isomorphism classes in 
the category $\text{hmf}(V-W)$, where we have omitted the lists of variables for readability, see \Cref{Subsection: Landau--Ginzburg theory}.
Vertical composition in $\RW$ is analogous to horizontal composition in $\mathscr{LG}_{\mathbb C}$, while horizontal composition in $\RW$ is basically addition of polynomials. 
The monoidal structure on $\RW$ amounts to concatenations of lists, addition of polynomials and tensor product over~$\C$ on the level of objects, 1- and 2-morphisms, respectively. 
For details we refer to \cite{1,Brunner2023TruncatedAR}.

The main technical result relevant for us is that $\RW$ admits a (trivial) pivotal (i.e.\ $\Sotwo$-dagger) structure. 
The left and right adjoint of a 1-morphism~$W$ is simply $-W$ (which allows for a trivial map from~$W$ to its double adjoint), but the adjunction 2-morphisms are non-trivial, see \cite[Theorem 1.1 \& Section 4.2.2]{Brunner2023TruncatedAR}. 
In the present section, we study $\rSpintwo$-dagger structure on $\RW$ from the purely algebraic perspective of the present paper. 
This reproduces and refines the results of \cite{Brunner2023TruncatedAR}, adding another consistency check of our approach. 

\medskip 

We start by considering the $\rSpintwo$-case, for $r \in \Z_{\geqslant 1}$. 
As an application of \Cref{Example: Symmetric monoidal 2-categories with duals and adjoints are rSpintwo volutive}, the 2-category $\RW$ carries a $\rSpintwo$-volutive structure, from which we may construct the $\rSpintwo$-dagger 2-category $S_{\rSpintwo}\RW$. 
By definition, the objects in the latter are pairs consisting of an object ${\xu} \in \RW$ together with a trivializations $\lambda_{{\xu}} \colon S_{{\xu}}^r \to \id_{{\xu}}$ of the $r$-th power of the Serre automorphism. 
The Serre automorphism~$S_{{\xu}}$ allows precisely two trivializations for every object~${\xu}$, with one represented (up to unitor 2-morphisms, cf.\ the discussion around \cite[(2.59)]{1}) by the identity matrix factorization $I_{\id_{{\xu}}}$, and the other by its $\Z_2$-shift $I_{\id_{{\xu}}}[1]$; we remark that 
both these trivializations lead to isomorphic fully extended oriented TQFTs, see \cite[Theorem 3.4]{1}.
In particular, $\RW$ embeds into the 2-category underlying $S_{\rSpintwo}\RW$ and hence admits a $\rSpintwo$-dagger structure. 
For the following result we restrict to $r=1$ and choose the embedding that assigns each object ${\xu} \in \RW$ to the pair $(x,[I_{\id_{{\xu}}}]) \in S_{\rSpintwo}\RW$.

\begin{proposition} \label{Proposition: comparison of pivotal structures on RW}
    For $r=1$, the $\Sotwo$-dagger structure induced on $\RW$ as explained in the previous paragraph coincides with the one described in 
    \cite[Theorem 1.1]{Brunner2023TruncatedAR}. 
\end{proposition}

\begin{proof}
    It is straightforward to check that the 2-morphism components of the Serre automorphism on $\RW$ are trivial by using the explicit description of right adjoints 
    in $\RW$. 
    One then translates our convention for $\Sotwo$-dagger structures into the one of loc.\ cit.\ by constructing left adjoints from the right adjoints 
    and the $\Sotwo$-dagger structure as explained in \Cref{Remark: Equivalent ways of talking about pivotal structures} for $r=1$. 
    The claim then follows.
\end{proof}

\begin{remark}\label{Remark: usefulness of construction in RW}
    \Cref{Proposition: comparison of pivotal structures on RW} illustrates the usefulness of our approach to $\Sotwo$-dagger structures: 
    while in \cite{Brunner2023TruncatedAR} the $\Sotwo$-dagger structure was found by an educated guess, our construction is canonical once the symmetric monoidal rigid structure on the 2-category $\RW$ with adjoints is established. 
    Moreover, we deduce the coherences one has to check directly from the properties of the Serre automorphism. 
\end{remark}

We note that our approach automatically provides us with a $\rSpintwo$-dagger structure on $\RW$ for $ r > 1$, though in this particular case the latter contains essentially 
the same information as the $\Sotwo$-dagger structure, as every object has trivializable (first power of the) Serre automorphism. 

\medskip 

Having discussed the known higher dagger structures on $\RW$, we now apply our theory to the $\Otwo$-case, which had not been considered previously. 
By \Cref{Example: symmetric monoidal 2-categories with duals are Otwo volutive}, we obtain an $\Otwo$-dagger 2-category $S_{\Otwo}\RW$ from the symmetric monoidal rigid 2-category $\RW$. 
We claim that any object in $\RW$ can be completed into an object of $S_{\Otwo}\RW$. 
Indeed, this follows essentially from the explicit description of duals 
in \cite[Section 4.2.3]{Brunner2023TruncatedAR} implying that any object ${\xu} \in \RW$ is self-dual, i.e.\ ${\xu}^{*} = {\xu}$. Combining this with the discussion of the $\Sotwo$-case, 
any object in $\RW$ can be completed into an object in the 2-category underlying $S_{\Otwo}\RW$. Opposed to considering the $\Otwo$-dagger 2-category $S_{\Otwo}\RW$ which 
keeps track of all trivialization data at once, one may also consider the 2-category $\RW$ with a specific choice of trivialization data for each object, equipping it with a (non-canonical) 
$\Otwo$-dagger structure; this point of view is often taken in the older literature on 2-categories with pivotal structures. 
Spelling this $\Otwo$-dagger structure out is a straightforward exercise and essentially amounts to spelling out the Serre automorphism (as done previously) and the 
$\Oone$-dagger structure on $\RW$ based on the dualization 2-functor, whose action on object level is trivial while its action on 1- and 2-morphisms is described in 
\cite[Section 4.2.3]{Brunner2023TruncatedAR}. The higher coherence data of the $\Oone$-dagger structure is trivial.

\medskip

Finally, we recall from \Cref{Subsection: Otwovolutive 2-categories and Oonevolutive 1-categories} that the Hom categories of $S_{\Otwo}\RW$ carry 
$\Oone$-volutions, to which one may apply $S_{\Oone}$ resulting in an $\Otwo$-dagger 2-category whose Hom categories carry compatible $\Oone$-dagger 
structures. 
Determining the objects in $S_{\Oone}\Hom({\xu},{\yu})$ is more difficult than the previous discussion; a partial description can be 
inferred from the discussion in \Cref{Subsection: Landau--Ginzburg theory}, noticing that the Hom categories of $\RW$ are closely related to the homotopy 1-category of $\LG$, on 
which we have discussed dagger structures in \Cref{Subsection: Landau--Ginzburg theory}.

\begin{remark}
    We stress that, while the $\Sotwo$-dagger structure we reproduced in \Cref{Proposition: comparison of pivotal structures on RW} was previously known, the 
    $\Otwo$-dagger structure is described here for the first time. This further illustrates the argument in \Cref{Remark: usefulness of construction in RW} on the usefulness 
    of our approach to higher dagger structures.
\end{remark}

\addcontentsline{toc}{section}{References}
\printbibliography

\medskip 

\noindent
Universität Wien, Fakultät für Physik, Boltzmanngasse 5, 1090 Wien, Österreich
\\
\href{mailto:nils.carqueville@univie.ac.at}{nils.carqueville@univie.ac.at}, 
\href{mailto:tim.lueders@univie.ac.at}{tim.lueders@univie.ac.at}
\end{document}